\documentclass[11pt]{article}
\usepackage{graphicx}
\usepackage{amssymb}

\usepackage{amsmath}

\usepackage{hyperref}

\usepackage{caption}

\def\ie{{\rm i.e.,\/}\ }

\def\etc{{\rm etc.\/}\ }

\newcommand{\nco}{\newcommand}
\nco{\munite}{\ensuremath{\,\,\mathrm{l}\!\!\!1}} 
\nco{\ZZ}{\mathbb{Z}}
\nco{\QQ}{\mathbb{Q}}
\nco{\CC}{\mathbb{C}}

\title{Quantum symmetries and nimreps for an exotic exceptional modular invariant of $F_4$ at level $6$}

\author{R. Coquereaux$^1$}
\date{%
    $^1$ Aix Marseille Univ, CNRS, I2M, Marseille, France\\[2ex]%
}

\begin{document}
\maketitle
\begin{abstract} 
We revisit an old non-diagonal modular invariant of $F_4$ at level $6$, discovered in 1989. 
At the time of its discovery, it was noticed that it could not be obtained from a conformal embedding of affine Lie algebras.
Later, it was shown that it can be obtained from an extension of a vertex operator algebra which is not of Lie type. 
To this modular invariant, one can attach a full CFT defined (or formulated) as a module category $\mathcal M$ for the fusion category $\mathcal A = \mathcal C(F_4,6)$. More recently, it was recognized that $\mathcal M$ is associated with an étale algebra of $\mathcal A$. In other words, the resulting quantum module is a quantum subgroup (it is also said to have self-fusion, to be of type I, or to be flat).
We study the dual $\mathcal O$ of $\mathcal A$ with respect to $\mathcal M$. In particular, we describe the algebra of quantum symmetries (we obtain the chiral generators and display the Ocneanu graphs for the fundamental irreps of $\mathcal A$).
The chiral graphs are not connected, each one containing two pairs of isomorphic components. The vertex sets of these components are identified with the simple objects, up to equivalence, of two module categories $\mathcal M$ and $\mathcal M^\prime$ sharing the same modular invariant.
$\mathcal M$ enjoys self-fusion, whereas $\mathcal M^\prime$ does not.
The restriction of the action of the Verlinde algebra $A$ to the associated abelian groups
$M$ and $M^\prime$ defines two families of matrices with non-negative integer entries (nimreps) and their adjacency graphs. 
As a by-product, we obtain the matrices describing the defect lines of the theory, the integer and rational fusion polynomials attached to $A$ and its modules, and the relevant numerical invariants, quantum dimensions, quantum masses, and induction-restriction rules.
\end{abstract}

\section{Introduction}
    
\subsection{Extended abstract}

We revisit an old non-diagonal modular invariant of $F_4$ at level $6$ discovered in 1989 by A.N. Schellekens and S. Yankielowicz \cite{Schellekens-Yankielowicz}.  
More recently, it was shown \cite{LamAndLin} that it can be obtained from an extension of a vertex operator algebra which is not of Lie type.
Because it was neither diagonal nor obtained by a simple-current extension, it was called exceptional.
Most exceptional modular invariants can be obtained from a conformal embedding of affine Lie algebras, possibly followed by a contraction procedure.
At the time of its discovery, it was noticed that this particular invariant could not be obtained in this manner. 
For this reason we call it exotic.       
Other authors \cite{Gannon.ExoticQuantumSubgroups1} use the word exotic in the above sense of exceptional, and refer to the latter invariant as an example of an exotic invariant which is not of Lie type.
   
To this modular invariant, one can attach a full CFT defined (or formulated) as a module category $\mathcal M$ for the fusion category $\mathcal A = \mathcal C(F_4,6)$. Even more recently \cite{Gannon.ExoticQuantumSubgroups1}, it was recognized that this module category is associated with an étale algebra of $\mathcal A$. In other words, the resulting quantum module is a quantum subgroup --- it is also said to have self-fusion, or to be of type I, or to be flat.
We also consider $\mathcal O = End_{\mathcal A} (\mathcal M)= {\mathcal F}un_{\mathcal A}(\mathcal M,\mathcal M)$, the endomorphism category of $\mathcal M$, also called the dual of $\mathcal A$ with respect to $\mathcal M$. 

In the sequel we shall most of the time avoid the category-theoretic language \cite{Ostrik} and stay at the level of the (Grothendieck) groups, or rings, of the relevant categories $\mathcal A$, $\mathcal M$, and $\mathcal O$, that we call $A$, $M$ and $O$. The ring $A$ (and its various field extensions) is called the fusion algebra, or the Verlinde algebra, the module $M$ is often called a nimrep (stricto sensu), and the ring $O$, with its various field extensions, is often called the algebra of quantum symmetries of $M$, or Ocneanu algebra of $M$.
One has to remember that $M$ is a module over $A$, and over $O$.
The ring $O$ itself is a bimodule over the Verlinde algebra $A$ (with two actions called left and right).

Starting from the modular invariant and the knowledge of the fusion matrices, we solve the so-called modular splitting equations
as well as a family of intertwining equations that allow us to obtain matrices describing the bimodule structure of $O$ over $A$.
In particular, for all simple objects (irreps) of $\mathcal A$, we obtain matrices (chiral generators) that describe the left and right actions of $A$ on $O$.
To any integrable irrep of $(F_4,6)$, in particular to its four fundamental irreps (since they happen to be integrable at that level), 
one can attach a non-connected graph defined as the union of the adjacency graphs of the two chiral generators. This is the Ocneanu graph of the chosen irreducible representation.

The resulting chiral graphs are themselves non-connected, each one containing two isomorphic copies of a first graph $M$ and two isomorphic copies of another graph $M^\prime$. 
The restriction of the action (the left one, for instance) of $A$ on $O$ to the first copy of $M$ (we take the connected component of the unit of $O$) defines a representation of the fusion algebra $A$ in terms of a family of matrices with non-negative integer entries (nimreps). The action of $A$ on $M^\prime$ defines another nimrep. 
The simple objects of the module category $\mathcal M$ are identified with the vertices of $M$, and the vertices of $M^\prime$ are associated with the simple objects of another module category $\mathcal M^\prime$.
By construction, these two inequivalent quantum modules (the first being a quantum subgroup) share the same modular invariant. 

The restriction of the ring structure of $O$ to the connected component of its unit element defines a self-fusion on $M$, and the unit itself selects a particular vertex of $M$ which, at a later stage, defines a commutative algebra $\mathcal F$ associated by induction with a non-simple object of $\mathcal A$ (étale structure).
The same procedure does not work for $M^\prime$ (no self-fusion).

Going back to the general case, let $\mathcal M$ be a module category over the modular monoidal (or fusion) category $\mathcal A = \mathcal A(\frak g,k)$ defined in terms of the integrable weights of an affine Lie algebra $\frak g$ at level $k$, or
as a semi-simplified category of representations of an associated quantum group at roots of unity.
In general, $\mathcal M$ is neither modular nor monoidal. $\mathcal O$ is monoidal, but in general it is not modular, although some subcategory of $\mathcal O$, called $\mathcal J$, is. 
The Grothendieck ring of $\mathcal J$ is a subring $J$ of $O$ sometimes called its ``ambichiral'' part. In categorical parlance, $\mathcal J$ is called the dyslectic part of $\mathcal O$.
In the present particular case, (the graph of) $M$ is identified as a special connected component of (the graph of) $O$ which contains all the ambichiral vertices. This means that the action of the modular group is defined not only for a particular collection of vertices of $O$ (the ambichiral ones), but also for a particular set of vertices of $M$ (a set that we shall again call $J$), sometimes called modular vertices.

A reader familiar with the historical examples defined by the ADE graphs of the Cappelli-Itzykson-Zuber classification \cite{CIZ_classification} may remember that the module category defined by the $E_6$ graph has three modular points, and that it can indeed be obtained as the connected component of the identity of the Ocneanu graph defined by the choice of the fundamental representation of $SU(2)$ ---the left graph, for instance.
The resulting quantum subgroup is exceptional, obtained from a conformal embedding of $(SU(2), 10)$ into the affine Lie algebra $(B_2,1)$, with the ambichiral part of $O$ (or the modular points of $M$) identified with the fusion algebra $J$ of $(B_2,1)$. In our present case, there is no conformal embedding into any affine Lie algebra, but the structure is similar.

Most techniques used in the forthcoming sections are certainly known to many people and have been discovered or rediscovered several times.
As for the modular splitting equations and subsequent intertwining equations, which give rise to the family of toric matrices (twisted partition functions) and to the chiral generators described by Ocneanu graphs, we have followed the methods developed and described in \cite{RCGSsu4, CTR, CH,IsasiSchieber}, drawing on, and in some cases rediscovering, ideas only briefly sketched in \cite{OcGoto, OcBariloche}.
   
\subsection{Notations}
 
From now on we often denote the categories and their Grothendieck groups, or rings, by the same capital letter. 
The context should be, most of the time, sufficient to understand which is which.
Here we consider a general situation but we shall return to our specific case in the next section.

The simple objects of $\mathcal A$ are denoted $m,n,p, \ldots$ They are also called irreps because they correspond to those irreducible representations of $\frak g$ that exist at level $k$.
The simple objects of $\mathcal M$ are denoted $a,b,c, \ldots$ In CFT (or rather BCFT: Boundary Conformal Field Theory) they are the boundary-condition labels.
The simple objects of $\mathcal O$ are denoted $x,y,z, \ldots$ In CFT, they are the defect-line labels; they are also called seams.
The matrices $N_m = (N_m)_{np}$ describing the monoidal structure of $\mathcal A$ are called fusion matrices.
The matrices $F_m = (F_m)_{ab}$ describing the $A$-module structure of $M$ are often called the ``nimreps'' of $M$ in the literature, or matrices of the nimreps, or fused adjacency matrices, and sometimes ---in particular by the author--- ``annular matrices''.
The matrices $S_x = (S_x)_{ab}$ describing the $O$-module structure of $M$ are called ``dual annular matrices''.
The matrices $O_x = (O_x)_{yz}$ describing the monoidal structure of $O$ are called generators of quantum symmetries.

The matrices $W_{x,y} = (W_{x,y})_{mn}$ describing the bimodule structure of $O$ over $A$ are defined by $m \, x \, n = \sum_y (W_{x,y})_{mn} \, y$.
They are called ``toric matrices'', actually ``toric matrices with two twists'' since the $W_{x,1}$ and the $W_{1,x}$ are called ``toric matrices with one twist''. 
The special $Z=W_{1,1}$ matrix is called ``the modular invariant'' because it commutes with the elements of the modular group $SL(2,\ZZ)$; it encodes the partition function in the BCFT defined by $M$.
The toric matrices with one or two twists can be interpreted in terms of partition functions for a theory with boundary conditions twisted by one or two defect lines, see the work of \cite{PetkovaZuber_defects}.
In general they are not invariant under the modular group, only $W_{1,1}$ is.
One can set $V_{mn} = (V_{mn})_{xy} := (W_{x,y})_{mn}$. The matrices $V_{mn}$ are called ``bimodule intertwiners'', and sometimes ``double fusion matrices''.
They describe the full bimodule structure of $O$ over $A$.
It is usual to set $V^L_m = V_{m1}$ and $V^R_m = V_{1m}$ for the chiral (left and right) bimodule intertwiners.

Call $f_i$ the simple fundamental objects of $\mathcal A$ that exist at the chosen level; so the $f_i$'s are those particular irreps $m, n, p, \ldots$ that are fundamental as irreps of the underlying classical Lie algebra, and exist at level $k$ (\ie are integrable as irreps of the affine Lie algebra at level $k$).
In the present example $(F_4,6)$, the four fundamental irreps exist at level $6$.
The special annular matrices $F_f$, where $f$ is a fundamental object, are also called ``fundamental''.
The special intertwiners $V^L_f$ and $V^R_f$ are called chiral (left and right) fundamental intertwiners.

The adjacency graph (with multilines and self-loops) associated with $V^L_f$ (resp. $V^R_f$) is the left (resp. right) Ocneanu graph for the fundamental irrep $f$.
The union of these two chiral graphs is the Ocneanu graph associated with the fundamental irrep $f$.
The intersection of the two chiral graphs is the graph of $J$, the ambichiral part of $O$.
A chiral graph (left or right) has, in general, several connected components.

The induction-restriction rules between $A$ and $M$ are also encoded by the (square) $F_n$ matrices, but they are better displayed in terms of the rectangular matrices $Ess_a=(Ess_a)_{n,b}=(F_n)_{a,b}$, called essential matrices (from $a$).

\underline {Case of self-fusion}. \\
We now assume that $M$ is of type I, i.e., that it is a quantum subgroup, not only a quantum module.

The finite group and vector bundle metaphor: think of $A$ as a (quantized) finite group, $M$ as a finite subgroup, $A/M$ as a homogeneous space, and $A$ as a principal bundle with fibres of type $M$ above $A/M$.
For each irreducible representation $a$ of $M$ one considers the vector bundle over $A/M$ with a typical fibre that (classically) is a representation space for $a$, and with a space of sections that we call $\Gamma_a$.
In particular, if $a$ is the identity representation, the space of sections $\mathcal F = \Gamma_1$ is the commutative algebra of functions over $A/M$. By the analogue of Peter-Weyl (induction), the spaces $\Gamma_a$, and in particular the algebra $\mathcal F$ itself, are representation spaces for $A$ and can be decomposed as sums of irreps of $A$.
The decomposition of the space of sections $\Gamma_a$, an $\mathcal A$-representation induced from a simple object $a$ of $\mathcal M$, into a sum of irreps of $A$ is encoded by the $a$-column of the essential matrix $Ess_1$.

In the case of self-fusion, the abelian group generated by the irreps of $M$ becomes a ring compatible with the action of $A$.
Using the above finite-group analogy, it is the tensor product of irreps of $M$.
One introduces matrices $g_a = (g_a)_{bc}$, sometimes called Pasquier graph matrices \cite{Pasquiergraphs}, that encode this ring multiplication.

The ambichiral vertices of (the graph of) ${O}$ are labelled $\hat x, \hat y,\ldots$ They generate a subring $J$ of $O$.
They are identified with special vertices $\hat a, \hat b,\ldots$ of (the graph of) $M$ called modular points, or modular vertices.
By induction, the latter determine particular spaces of sections $\Gamma_{\hat a}$ called dyslectic modules ---they carry a representation of $A$ but some of them can be equivalent.
 
${\mathcal F} ={ \mathcal F}(A/M) = \Gamma_1$, where $1$ is the unit of $M$, is an étale (commutative) algebra, an object ---not simple in general--- in the category $\mathcal A$. It is the space of would-be functions over the quantum homogeneous space $A/M$. This commutative algebra structure attached to $\Gamma_1$ should not be confused with the self-fusion of $M$. 

\subsection{A selection of results}
   
The fusion ring $A$ of $\mathcal A(F_4,6)$ has $r_{\mathcal A}=39$ simple objects. They are given in the text together with their quantum dimensions, conformal weights, etc.
The generators $S$ and $T$ of the modular group are obtained by the Kac-Peterson formula. 
The fusion matrices are then obtained by using the Verlinde formula.

The exceptional and exotic invariant of Schellekens-Yankielowicz (1989) is compatible with the existence of a module category.
We actually find two of them: one, $\mathcal M$, is a quantum subgroup (it has self-fusion and one of its simple objects ---the unit--- determines a connected étale algebra structure); the other, $\mathcal M^\prime$, is only a quantum module. 
Both have $r_M=16$ simple objects. The algebra of quantum symmetries $O$ is the same for $\mathcal M$ and $\mathcal M^\prime$.

$\mathcal O$ possesses $r_{\mathcal O}=64$ simple objects $x$, but the rank of the family generated by the left (or right) representation of $A$ on $O$ is $48$ ---the rank of the family of matrices $V^L_m$ where $m$ runs from $1$ to $39$. The rank of the family of toric matrices $W_x$, with one twist, is actually $40$.

For each fundamental $f_i$ of $A$, with $i$ in $(1,2,3,4)$, the graph of $V^L_f$ has four connected components forming two pairs of isomorphic components. Each connected component contains $16$ vertices.
$\mathcal M$ possesses $r_M=16$ simple objects, but the rank of the family generated by the representation of $A$ on $M$ is $12$.

The graph of $V^L_f$ (for each $f$) contains two isomorphic copies of $M$, and also two isomorphic copies of $M^\prime$. 
The graph(s) of $V^R_f$, for each $f$, can be deduced from the graph(s) of $V^L_f$ by a permutation acting on the $64$ vertices. So, both left and right graphs share the same features.

The restriction of the bimodule action of $A$ on $O$ (choose one) determines the action of the fusion algebra $A$ on the two copies of $M$ and of $M^\prime$. For definiteness, $M$ is defined as the connected component that contains the unit of $O$ in the left chiral graph (choosing the right graph leads to an equivalent description of the same nimrep).

The self-fusion on $M$ is inherited from the multiplication in $O$.
$M^\prime$ does not enjoy self-fusion because the restriction of the $O$-multiplication is not inner in $M^\prime$.

Once $V^L_f$ and $V^R_f$ are determined for each of the four fundamental $f$, one can obtain the other $V^L_m$'s and the other $V^R_m$'s for arbitrary integrable irreps $m$, because $m$ can be expressed polynomially in terms of the fundamental irreps, using $\ZZ$-polynomials. Actually, the number of generators can be reduced if one uses $\QQ$-polynomials and primitive elements that are appropriate sums of fundamental irreps.

For various choices of $f$, the resulting annular matrices (nimreps) for $M$ and $M^\prime$ and their adjacency graphs are explicitly described in the text.

For each fundamental $f$, the intersection of $V^L_f$ and $V^R_f$ is a subgraph with $r_J=6$ vertices ---the ambichiral vertices labelled $\{1, 18, 32, 33, 60, 61\}$ in the graphs of $O$ given in sec.~\ref{OcneanuGraphs}, in particular in fig.~\ref{FLR1000M} and fig.~\ref{FLR0001M}.
The associated modular vertices of $M$, also called $J$, are labelled $\{1, 6, 8, 9, 15, 16\}$ in the graphs given in sec.~\ref{qModuleM} and correspond, by induction, to the collection of dyslectic modules (it turns out that only four are inequivalent).

One might naively expect $O$ to behave like $M \otimes_J M$, which at the level of ranks would suggest $r_O \sim r_M^2/r_J$, but in the present case this naive counting cannot be correct ($16^2/6$ is not an integer!). The correct result is obtained by writing $64 = 16 \times 16/4$. This phenomenon can be traced back to a degeneracy involving two pairs of ambichiral vertices (or modular points of $M$), those called $(8, 9)$ and $(15, 16)$ in the graphs of $M$,
see the details in sec.~\ref{dyslecticdecomp}, and the number of inequivalent dyslectic modules is $4$, not $6$.

Ultimately, one can also trace this to the existence of a ``module-over-category'' associated with the fusion algebra $F_4 \times A_2$ at level $(6,2)$, with central charge $24$. Its $(A_2,2)$ simple component has $6$ simple objects that correspond to the classical irreps $00, 01, 10, 20, 11, 02$ of $A_2$, and the exotic example that we study is obtained from the latter as a quotient that identifies the pairs $(01,10)$ and $(02, 20)$ ---each pair giving rise to real irreps. This can also be seen from the presence of entries of multiplicity $2$ in the modular invariant matrix. We shall come back to this discussion in sec.~\ref{modularZ}.

\subsection{Structure of the paper}
   
The fusion algebra alone is studied in sec.~\ref{F4level6}.

Several properties of the exotic exceptional invariant of Schellekens-Yankielowicz \cite{Schellekens-Yankielowicz} are discussed in sec.~\ref{ModularConsiderations}, in relation with later considerations ($F_4 \times A_2$, VOAs, $\ldots$). These comments are not used in the sequel, but we include them for the sake of completeness.

The modular splitting equation and the method used to solve it are studied in sec.~\ref{modularsplittingeq}. 
The output of this section is the list of toric matrices with one twist $W_x$ that encode part of the bimodule structure of $O$.
 
The intertwining equations for $O$ are studied in sec.~\ref{Intertwiners}. There we determine the left and right actions of $A$ on $O$, in terms of chiral generators $V^L_m$ and $V^R_m$, and display their Ocneanu graphs for various choices of fundamental representations.
The determination of the $V^{L,R}_m$ is first carried out for selected fundamental irreps ($0001$ and $0010$), then generalized to arbitrary irreps $m$ by using integral or rational fusion polynomials. 
In the same section we consider the restriction of these chiral generators to $M$ and to $M^\prime$.
 
Section \ref{qModuleM} is devoted to a detailed study of $M$: induction-restriction rules, 
decompositions of the spaces of sections $\Gamma_a$, in particular for the commutative algebra $\mathcal F=\Gamma_1$, 
quantum dimensions, masses, reduced rational polynomials for $M$, etc.
 
Section \ref{qModuleMprime} is devoted to $M^\prime$.
  
Several other results of interest, or open questions, are gathered in section \ref{varia}.
  
 \section{The fusion algebra of $F_4$ at level six}
   \label{F4level6}
  \subsection{Integrable representations}
  \subsubsection{Simple objects}
  \label{simpleobjects}
 The fusion category $\mathcal A(F_4,6)$ has 39 simple objects (irreps). More generally, denoting by $r_{\mathcal A}(k)$ the rank of $\mathcal A(F_4,k)$, one can show that $r_{\mathcal A}(k)=1,2,5,9,16,25,39,56,\ldots$, for $k=0,1,2,3,4,5,6,7,\ldots$.
 Each simple object is associated with some classical irreducible representation (irrep), and to a classical highest weight.
 Since we shall mostly consider the case $k=6$, we just set $r_{\mathcal A}=39$.
 Each simple object can be labelled by the components\footnote{The list of components $(i,j,k,l)$ is often written ijkl, to save space, whenever the labels are smaller than $10$.} of the highest weight (hw) of its classical counterpart in the basis of fundamental weights (the so-called Dynkin basis).
 The ordering chosen for the Dynkin diagram of $F_4$ is such that the dimensions of the classical  irreps $1 0 0 0$, $0 1 0 0$, $0 0 1 0$, $0 0 0 1$ are respectively $52$, $1274$, $273$, $26$. 
 Therefore $1 0 0 0$  denotes the hw of  the adjoint representation.
 The irrep with hw $0001$ is sometimes called ``basic'' because it has the smallest (classical) dimension.
The list of simple objects obtained at level $6$  is given below. It is an ordered set. The order is as follows:  If the weights $\lambda$ and $\mu$ appear respectively at levels $i$ and $j$, with $i < j$,  then $\lambda < \mu$.
This is a partial order which becomes a total order if we decide, when $i=j$, to sort the weights lexicographically according to their components, starting with the first (canonical order). For instance 
$0  0  0  3<0 1 0 1$ because the first appears at level 3 and the other at level 4, but  $0 1 0 1 < 1 0 1 0$ where both appear at level 4. A semicolon in the list separates the weights appearing at increasing levels.
The list of irreps reads:

$
 0  0  0  0 ;
   \\
 0  0  0  1 ;
   \\
 0  0  0  2 ,
 0  0  1  0 ,
 1  0  0  0 ;
   \\
 0  0  0  3 ,
 0  0  1  1 ,
 0  1  0  0 ,
 1  0  0  1 ;
   \\
 0  0  0  4 ,
 0  0  1  2 ,
 0  0  2  0 ,
 0  1  0  1 ,
 1  0  0  2 ,
 1  0  1  0 ,
 2  0  0  0 ;
  \\
 0  0  0  5 ,
 0  0  1  3 ,
 0  0  2  1 ,
 0  1  0  2 ,
 0  1  1  0 ,
 1  0  0  3 ,
 1  0  1  1 ,
 1  1  0  0 ,
 2  0  0  1 ;
  \\
 0  0  0  6 ,
 0  0  1  4 ,
 0  0  2  2 ,
 0  0  3  0 ,
 0  1  0  3 ,
 0  1  1  1 ,
 0  2  0  0 ,
 1  0  0  4 ,
 1  0  1  2 ,
 1  0  2  0 ,
 1  1  0  1 ,
 2  0  0  2 ,
 2  0  1  0 ,
 3  0  0  0.$

We call ``fundamental simple objects'' (at a chosen level $k$) those simple objects whose highest weights correspond to fundamental weights lying in the level-k Weyl alcove.
The four fundamental representations of $F_4$ do exist at level $3$, and therefore also at level $6$.
Notice  that these fundamental simple objects at level 6, namely those with labels 1000, 0100, 0010, 0001,  appear respectively at positions 5, 8, 4, 2 in the previous list.
 
   \subsubsection{Quantum dimensions}
     \label{quantumdimensionsofA}
We set $q = \exp (i \pi / \kappa)$ where $\kappa$ (the altitude) is $\kappa= g + k$, with $g = 9$ (dual Coxeter number of $F_4$) and $k = 6$ (the chosen level). So $\kappa = 15$ and $q$ is a 30th root of unity.
The quantum dimensions (or categorical dimensions) of simple objects can be obtained from the quantum version of the Weyl dimension formula. 
For an irrep $\lambda$ with components $(\lambda _1, \lambda _2, \lambda _3, \lambda _4)$, the q-dimension is $\text{qdim}(\lambda)=N/D$ where \\
{ \footnotesize
  $ N =
{}[\lambda _1+1] {}[\lambda _2+1] {}[\lambda _1+\lambda _2+2] {}[\frac{\lambda _3}{2}+\frac{1}{2}] {}[\lambda _2+\frac{\lambda
   _3}{2}+\frac{3}{2}] {}[\lambda _1+\lambda _2+\frac{\lambda _3}{2}+\frac{5}{2}] {}[\lambda _2+\lambda _3+2] {}[\lambda _1+\lambda _2+\lambda _3+3]
   {}[\lambda _1+2 \lambda _2+\lambda _3+4] {}[\frac{\lambda _4}{2}+\frac{1}{2}] {}[\frac{\lambda _3}{2}+\frac{\lambda _4}{2}+1] {}[\lambda _2+\frac{\lambda
   _3}{2}+\frac{\lambda _4}{2}+2] {}[\lambda _1+\lambda _2+\frac{\lambda _3}{2}+\frac{\lambda _4}{2}+3] {}[\lambda _2+\lambda _3+\frac{\lambda _4}{2}+\frac{5}{2}]
   {}[\lambda _1+\lambda _2+\lambda _3+\frac{\lambda _4}{2}+\frac{7}{2}] {}[\lambda _1+2 \lambda _2+\lambda _3+\frac{\lambda _4}{2}+\frac{9}{2}] {}[\lambda _1+2 \lambda
   _2+\frac{3 \lambda _3}{2}+\frac{\lambda _4}{2}+5] {}[\lambda _2+\lambda _3+\lambda _4+3] {}[\lambda _1+\lambda _2+\lambda _3+\lambda _4+4] {}[\lambda _1+2
   \lambda _2+\lambda _3+\lambda _4+5] {}[\lambda _1+2 \lambda _2+\frac{3 \lambda _3}{2}+\lambda _4+\frac{11}{2}] {}[\lambda _1+2 \lambda _2+2 \lambda _3+\lambda
   _4+6] {}[\lambda _1+3 \lambda _2+2 \lambda _3+\lambda _4+7] {}[2 \lambda _1+3 \lambda _2+2 \lambda _3+\lambda _4+8],$ and 
 $ D =   {{}[\frac{1}{2}]^2 [1]^3
   {}[\frac{3}{2}] {}[2]^3 {}[\frac{5}{2}]^2 {}[3]^3 {}[\frac{7}{2}] {}[4]^2 {}[\frac{9}{2}] {}[5]^2 {}[\frac{11}{2}] {}[6] {}[7]
   {}[8]},
 $}
where the square bracket $[n]$ stands for the q-number $\frac{q^n-q^{-n}}{q-q^{-1}}$. Since $F_4$ is not simply laced, one may notice that, with our definition of $q$, the argument $n$ of $[n]$ is not always an integer.

The fundamental irreps, 1000, 0100, 0010, 0001, have q-dimensions  
{\small
\[ 
\frac{[6] [\frac{13}{2}] [10]}{[1] [\frac{5}{2}] [3]},
\frac{[\frac{7}{2}] [4] [\frac{13}{2}] [7] [\frac{15}{2}] [10] [11]}{[1] [\frac{3}{2}] [2]^2 [\frac{5}{2}] [5] [\frac{11}{2}]},
\frac{[\frac{7}{2}] [\frac{13}{2}] [9] [10]}{[\frac{1}{2}] [2] [\frac{5}{2}] [3]},
\frac{[\frac{3}{2}] [4] [\frac{13}{2}] [9]}{[\frac{1}{2}] [2] [3] [\frac{9}{2}]},
\]
}
 or, equivalently,
 {\footnotesize
\[
\left\{
\frac{9}{2}+\frac{3 \sqrt{5}}{2}+\sqrt{3 \left(5+2 \sqrt{5}\right)},
\frac{37}{2}+\frac{17 \sqrt{5}}{2}+\sqrt{\frac{3}{2} \left(445+199 \sqrt{5}\right)},
\frac{21}{2}+\frac{9 \sqrt{5}}{2}+\sqrt{3 \left(85+38 \sqrt{5}\right)},
3+\sqrt{5}+\sqrt{3 \left(5+2 \sqrt{5}\right)}
\right\}
\]
}   
   
 \subsubsection{Global dimension}
      \label{globaldimensionofA}
   
The global dimension of  $\mathcal A(F_4,6)$, i.e. the sum of squares of categorical dimensions of its simple objects, is 
$$\vert {\mathcal A} \vert =  {dim_{\mathcal A}}=900 \left(501+224 \sqrt{5}+68 \sqrt{6 \left(5-\sqrt{5}\right)} +110 \sqrt{6 \left(5+\sqrt{5}\right)}\right)$$
It is the largest root of the polynomial $900 (1-84 x+806 x^2-2004 x^3+x^4)$.
Numerically, it is close to  $1.80324 \times10^6$.
It can also be written trigonometrically as 
{ \footnotesize
$$-\frac{35230013}{2}+\frac{268767 \sqrt{5}}{2}-36021446 \sin \left(\frac{\pi }{30}\right)+36393380 \sin \left(\frac{7 \pi }{30}\right)-35601380 \cos \left(\frac{\pi }{15}\right)+36511046 \cos
   \left(\frac{2 \pi }{15}\right)$$
   }

   The global dimension of ${\mathcal A}$ can also be expressed in terms of the first matrix entry of the generator $S$ of the modular group.
   One has 
   $$\vert {\mathcal A} \vert  = \frac{1}{S_{1,1}^2}$$
  The fact that this equality  holds provides a check of the previous determination of the quantum dimensions of simple objects.

    \subsubsection{Conformal weights}
 Remember that the conformal weight $h(\lambda)$ of an irrep $\lambda$ (which can be  specified by the highest weight of its classical counterpart) is given by $C(\lambda)/(2(k+g))$
 where $C(\lambda)$ is the Casimir eigenvalue of $\lambda$, and where the denominator of $h(\lambda)$ is twice the altitude.
 Calling $\rho$ the Weyl vector and $\langle \,, \, \rangle$ the canonical inner product, the Casimir eigenvalue is 
  given by   $C(\lambda) = \langle \lambda + \rho, \lambda + \rho \rangle -  \langle \rho,\rho  \rangle$.
     The factor $2$ in the denominator of  $h(\lambda)$ stands for the normalization factor of long roots.
  
 Using the order already introduced in the set of irreps, here is the list of conformal weights $h(\lambda)$:\\
{\small      
$(0, 2/5, 13/15, 4/5, 3/5, 7/5, 13/10, 6/5, 16/15, 2, 28/15, 9/5, 26/15, 8/5, 23/15, 4/3, 8/3, 5/2,
 12/5, 7/3, 34/15, \\ 11/5, 21/10, 2, 28/15, 17/5, 16/5, 46/15, 3, 3, 29/10, 14/5, 43/15, 41/15, 8/3, 13/5, 37/15, 12/5, 11/5)$.}

Only those irreps whose conformal weights are equal to $0$ modulo the integers can appear in the modular block of the identity, for any chosen modular invariant.
 More generally, irreps belonging to a specific modular block have the same conformal weights modulo the integers. These are necessary, but not sufficient, conditions. 
 For this reason, and because they will be used later, we give the list of reduced conformal weights $\overline h(\lambda)$, which is the same list as before, but with conformal weights taken modulo the integers, i.e. modulo $1$.\\
{\small      
$(0, 2/5, 13/15, 4/5, 3/5, 2/5, 3/10, 1/5, 1/15, 0, 13/15, 4/5, 11/15, 3/5, 8/15, 1/3, 2/3, 1/2, 2/5, 1/3, 4/15, 1/5, \\
 1/10, 0, 13/15, 2/5, 1/5, 1/15, 0, 0, 9/10, 4/5, 13/15, 11/15, 2/3, 3/5, 7/15, 2/5, 1/5)$.}

 \subsubsection{Central charges of $(F_4,6)$ and of $(F_4 \times A_2, (6, 2))$}
Remember that for a Lie group $G$ at level $k$, the central charge is given by $c = \dim(G) \times k /(k + g)$, where $g$ is the dual Coxeter number of $G$.
So, for $(F_4,6)$ one obtains  $(52 \times 6)/(6 + 9 )= 104/5$.
Notice that for $(A_2, 2)$, one gets $(8 \times 2)/(2 +3 )=16/5$, so that for the semisimple, but not simple, Lie group $F_4 \times A_2$ at level $(6, 2)$, one obtains a central charge equal to $104/5+16/5 = 24$.

 \subsection{The fusion ring}
 
 The category $\mathcal A(F_4,6)$ is modular. The generators $S$ and $T$ of the group $SL(2,\ZZ)$ can be obtained from the integrable representations by using the Kac-Peterson formula.
 
   \subsubsection{Fusion matrices}
   
 In turn, the fusion ring matrices $N_p$, where $p$ denotes a simple object, are obtained from the above modular generators $S$ and $T$ by  the Verlinde formula.
These matrices $N_f$, for $f = 0001$, $0010$, $0100$, and $1000$, in this order (starting from $0001$, the basic one), are displayed below. Dots stand for vanishing entries. In the present case, all  the fusion matrices are symmetric.

  {\tiny
\setlength{\arraycolsep}{1.5pt}
\[
N_{0001}=
\left(
\begin{array}{ccccccccccccccccccccccccccccccccccccccc}
 . & 1 & . & . & . & . & . & . & . & . & . & . & . & . & . & . & . & . & . & . & . & . & . & . & . & . & . & . & . & . & . & . & . & . & . & . & . & . & . \\
 1 & 1 & 1 & 1 & 1 & . & . & . & . & . & . & . & . & . & . & . & . & . & . & . & . & . & . & . & . & . & . & . & . & . & . & . & . & . & . & . & . & . & . \\
 . & 1 & 1 & 1 & . & 1 & 1 & . & 1 & . & . & . & . & . & . & . & . & . & . & . & . & . & . & . & . & . & . & . & . & . & . & . & . & . & . & . & . & . & . \\
 . & 1 & 1 & 1 & 1 & . & 1 & 1 & 1 & . & . & . & . & . & . & . & . & . & . & . & . & . & . & . & . & . & . & . & . & . & . & . & . & . & . & . & . & . & . \\
 . & 1 & . & 1 & . & . & . & . & 1 & . & . & . & . & . & . & . & . & . & . & . & . & . & . & . & . & . & . & . & . & . & . & . & . & . & . & . & . & . & . \\
 . & . & 1 & . & . & 1 & 1 & . & . & 1 & 1 & . & . & 1 & . & . & . & . & . & . & . & . & . & . & . & . & . & . & . & . & . & . & . & . & . & . & . & . & . \\
 . & . & 1 & 1 & . & 1 & 2 & 1 & 1 & . & 1 & 1 & 1 & 1 & 1 & . & . & . & . & . & . & . & . & . & . & . & . & . & . & . & . & . & . & . & . & . & . & . & . \\
 . & . & . & 1 & . & . & 1 & . & 1 & . & . & . & 1 & . & 1 & . & . & . & . & . & . & . & . & . & . & . & . & . & . & . & . & . & . & . & . & . & . & . & . \\
 . & . & 1 & 1 & 1 & . & 1 & 1 & 1 & . & . & . & . & 1 & 1 & 1 & . & . & . & . & . & . & . & . & . & . & . & . & . & . & . & . & . & . & . & . & . & . & . \\
 . & . & . & . & . & 1 & . & . & . & 1 & 1 & . & . & . & . & . & 1 & 1 & . & . & . & 1 & . & . & . & . & . & . & . & . & . & . & . & . & . & . & . & . & . \\
 . & . & . & . & . & 1 & 1 & . & . & 1 & 2 & 1 & 1 & 1 & . & . & . & 1 & 1 & 1 & . & 1 & 1 & . & . & . & . & . & . & . & . & . & . & . & . & . & . & . & . \\
 . & . & . & . & . & . & 1 & . & . & . & 1 & 1 & 1 & . & 1 & . & . & . & 1 & . & 1 & . & 1 & . & . & . & . & . & . & . & . & . & . & . & . & . & . & . & . \\
 . & . & . & . & . & . & 1 & 1 & . & . & 1 & 1 & 1 & 1 & 1 & . & . & . & . & 1 & 1 & . & 1 & 1 & . & . & . & . & . & . & . & . & . & . & . & . & . & . & . \\
 . & . & . & . & . & 1 & 1 & . & 1 & . & 1 & . & 1 & 1 & 1 & . & . & . & . & . & . & 1 & 1 & . & 1 & . & . & . & . & . & . & . & . & . & . & . & . & . & . \\
 . & . & . & . & . & . & 1 & 1 & 1 & . & . & 1 & 1 & 1 & 1 & 1 & . & . & . & . & . & . & 1 & 1 & 1 & . & . & . & . & . & . & . & . & . & . & . & . & . & . \\
 . & . & . & . & . & . & . & . & 1 & . & . & . & . & . & 1 & . & . & . & . & . & . & . & . & . & 1 & . & . & . & . & . & . & . & . & . & . & . & . & . & . \\
 . & . & . & . & . & . & . & . & . & 1 & . & . & . & . & . & . & 1 & 1 & . & . & . & . & . & . & . & 1 & 1 & . & . & . & . & . & 1 & . & . & . & . & . & . \\
 . & . & . & . & . & . & . & . & . & 1 & 1 & . & . & . & . & . & 1 & 2 & 1 & 1 & . & 1 & . & . & . & . & 1 & 1 & . & 1 & . & . & 1 & 1 & . & . & . & . & . \\
 . & . & . & . & . & . & . & . & . & . & 1 & 1 & . & . & . & . & . & 1 & 2 & 1 & 1 & . & 1 & . & . & . & . & 1 & 1 & . & 1 & . & . & 1 & 1 & . & . & . & . \\
 . & . & . & . & . & . & . & . & . & . & 1 & . & 1 & . & . & . & . & 1 & 1 & 1 & 1 & 1 & 1 & . & . & . & . & . & . & 1 & 1 & . & . & 1 & . & 1 & . & . & . \\
 . & . & . & . & . & . & . & . & . & . & . & 1 & 1 & . & . & . & . & . & 1 & 1 & 1 & . & 1 & 1 & . & . & . & . & . & . & 1 & 1 & . & . & 1 & 1 & . & . & . \\
 . & . & . & . & . & . & . & . & . & 1 & 1 & . & . & 1 & . & . & . & 1 & . & 1 & . & 1 & 1 & . & . & . & . & . & . & . & . & . & 1 & 1 & . & . & 1 & . & . \\
 . & . & . & . & . & . & . & . & . & . & 1 & 1 & 1 & 1 & 1 & . & . & . & 1 & 1 & 1 & 1 & 2 & 1 & 1 & . & . & . & . & . & . & . & . & 1 & 1 & 1 & 1 & 1 & . \\
 . & . & . & . & . & . & . & . & . & . & . & . & 1 & . & 1 & . & . & . & . & . & 1 & . & 1 & . & 1 & . & . & . & . & . & . & . & . & . & . & 1 & . & 1 & . \\
 . & . & . & . & . & . & . & . & . & . & . & . & . & 1 & 1 & 1 & . & . & . & . & . & . & 1 & 1 & 1 & . & . & . & . & . & . & . & . & . & . & . & 1 & 1 & 1 \\
 . & . & . & . & . & . & . & . & . & . & . & . & . & . & . & . & 1 & . & . & . & . & . & . & . & . & 1 & 1 & . & . & . & . & . & . & . & . & . & . & . & . \\
 . & . & . & . & . & . & . & . & . & . & . & . & . & . & . & . & 1 & 1 & . & . & . & . & . & . & . & 1 & 2 & 1 & . & 1 & . & . & 1 & . & . & . & . & . & . \\
 . & . & . & . & . & . & . & . & . & . & . & . & . & . & . & . & . & 1 & 1 & . & . & . & . & . & . & . & 1 & 2 & 1 & 1 & 1 & . & . & 1 & . & . & . & . & . \\
 . & . & . & . & . & . & . & . & . & . & . & . & . & . & . & . & . & . & 1 & . & . & . & . & . & . & . & . & 1 & 1 & . & 1 & . & . & . & 1 & . & . & . & . \\
 . & . & . & . & . & . & . & . & . & . & . & . & . & . & . & . & . & 1 & . & 1 & . & . & . & . & . & . & 1 & 1 & . & 1 & 1 & . & 1 & 1 & . & . & . & . & . \\
 . & . & . & . & . & . & . & . & . & . & . & . & . & . & . & . & . & . & 1 & 1 & 1 & . & . & . & . & . & . & 1 & 1 & 1 & 2 & 1 & . & 1 & 1 & 1 & . & . & . \\
 . & . & . & . & . & . & . & . & . & . & . & . & . & . & . & . & . & . & . & . & 1 & . & . & . & . & . & . & . & . & . & 1 & . & . & . & . & 1 & . & . & . \\
 . & . & . & . & . & . & . & . & . & . & . & . & . & . & . & . & 1 & 1 & . & . & . & 1 & . & . & . & . & 1 & . & . & 1 & . & . & 1 & 1 & . & . & . & . & . \\
 . & . & . & . & . & . & . & . & . & . & . & . & . & . & . & . & . & 1 & 1 & 1 & . & 1 & 1 & . & . & . & . & 1 & . & 1 & 1 & . & 1 & 2 & 1 & 1 & 1 & . & . \\
 . & . & . & . & . & . & . & . & . & . & . & . & . & . & . & . & . & . & 1 & . & 1 & . & 1 & . & . & . & . & . & 1 & . & 1 & . & . & 1 & 1 & 1 & . & 1 & . \\
 . & . & . & . & . & . & . & . & . & . & . & . & . & . & . & . & . & . & . & 1 & 1 & . & 1 & 1 & . & . & . & . & . & . & 1 & 1 & . & 1 & 1 & 1 & 1 & 1 & . \\
 . & . & . & . & . & . & . & . & . & . & . & . & . & . & . & . & . & . & . & . & . & 1 & 1 & . & 1 & . & . & . & . & . & . & . & . & 1 & . & 1 & 1 & 1 & . \\
 . & . & . & . & . & . & . & . & . & . & . & . & . & . & . & . & . & . & . & . & . & . & 1 & 1 & 1 & . & . & . & . & . & . & . & . & . & 1 & 1 & 1 & 1 & 1 \\
 . & . & . & . & . & . & . & . & . & . & . & . & . & . & . & . & . & . & . & . & . & . & . & . & 1 & . & . & . & . & . & . & . & . & . & . & . & . & 1 & . \\
\end{array}
\right)
\]
}

 {\tiny
\setlength{\arraycolsep}{1.5pt}
\[
N_{0010}=
\left(
\begin{array}{ccccccccccccccccccccccccccccccccccccccc}
 . & . & . & 1 & . & . & . & . & . & . & . & . & . & . & . & . & . & . & . & . & . & . & . & . & . & . & . & . & . & . & . & . & . & . & . & . & . & . & . \\
 . & 1 & 1 & 1 & 1 & . & 1 & 1 & 1 & . & . & . & . & . & . & . & . & . & . & . & . & . & . & . & . & . & . & . & . & . & . & . & . & . & . & . & . & . & . \\
 . & 1 & 1 & 2 & 1 & 1 & 2 & 1 & 2 & . & 1 & . & 1 & 1 & 1 & . & . & . & . & . & . & . & . & . & . & . & . & . & . & . & . & . & . & . & . & . & . & . & . \\
 1 & 1 & 2 & 2 & 1 & 1 & 2 & 1 & 2 & . & . & 1 & 1 & 1 & 1 & 1 & . & . & . & . & . & . & . & . & . & . & . & . & . & . & . & . & . & . & . & . & . & . & . \\
 . & 1 & 1 & 1 & . & . & 1 & . & 1 & . & . & . & . & . & 1 & . & . & . & . & . & . & . & . & . & . & . & . & . & . & . & . & . & . & . & . & . & . & . & . \\
 . & . & 1 & 1 & . & 1 & 2 & 1 & 1 & 1 & 2 & 1 & 1 & 2 & 1 & . & . & 1 & . & 1 & . & 1 & 1 & . & . & . & . & . & . & . & . & . & . & . & . & . & . & . & . \\
 . & 1 & 2 & 2 & 1 & 2 & 5 & 2 & 3 & 1 & 3 & 2 & 3 & 3 & 3 & 1 & . & . & 1 & 1 & 1 & 1 & 2 & 1 & 1 & . & . & . & . & . & . & . & . & . & . & . & . & . & . \\
 . & 1 & 1 & 1 & . & 1 & 2 & 1 & 2 & . & 1 & 1 & 1 & 1 & 2 & . & . & . & . & . & 1 & . & 1 & . & 1 & . & . & . & . & . & . & . & . & . & . & . & . & . & . \\
 . & 1 & 2 & 2 & 1 & 1 & 3 & 2 & 2 & . & 1 & 1 & 1 & 2 & 2 & 1 & . & . & . & . & . & . & 1 & 1 & 1 & . & . & . & . & . & . & . & . & . & . & . & . & . & . \\
 . & . & . & . & . & 1 & 1 & . & . & 1 & 2 & . & 1 & 1 & . & . & 1 & 2 & 1 & 1 & . & 2 & 1 & . & . & . & 1 & . & . & 1 & . & . & 1 & 1 & . & . & . & . & . \\
 . & . & 1 & . & . & 2 & 3 & 1 & 1 & 2 & 5 & 3 & 3 & 3 & 2 & . & 1 & 3 & 3 & 3 & 2 & 3 & 4 & 1 & 1 & . & . & 1 & . & 1 & 1 & . & 1 & 2 & 1 & 1 & 1 & . & . \\
 . & . & . & 1 & . & 1 & 2 & 1 & 1 & . & 3 & 2 & 3 & 2 & 2 & . & . & 1 & 2 & 2 & 2 & 1 & 3 & 1 & 1 & . & . & . & 1 & . & 1 & . & . & 1 & 1 & 1 & . & 1 & . \\
 . & . & 1 & 1 & . & 1 & 3 & 1 & 1 & 1 & 3 & 3 & 3 & 3 & 3 & 1 & . & 1 & 2 & 2 & 2 & 1 & 4 & 2 & 1 & . & . & . & . & . & 1 & 1 & . & 1 & 1 & 1 & 1 & 1 & . \\
 . & . & 1 & 1 & . & 2 & 3 & 1 & 2 & 1 & 3 & 2 & 3 & 2 & 3 & 1 & . & 1 & 1 & 1 & 1 & 2 & 3 & 1 & 2 & . & . & . & . & . & . & . & . & 1 & . & 1 & 1 & 1 & . \\
 . & . & 1 & 1 & 1 & 1 & 3 & 2 & 2 & . & 2 & 2 & 3 & 3 & 3 & 1 & . & . & 1 & 1 & 1 & 1 & 3 & 2 & 2 & . & . & . & . & . & . & . & . & . & 1 & 1 & 1 & 1 & 1 \\
 . & . & . & 1 & . & . & 1 & . & 1 & . & . & . & 1 & 1 & 1 & . & . & . & . & . & . & . & 1 & . & 1 & . & . & . & . & . & . & . & . & . & . & . & . & 1 & . \\
 . & . & . & . & . & . & . & . & . & 1 & 1 & . & . & . & . & . & 1 & 2 & . & 1 & . & 1 & . & . & . & 1 & 2 & 1 & . & 1 & . & . & 2 & 1 & . & . & . & . & . \\
 . & . & . & . & . & 1 & . & . & . & 2 & 3 & 1 & 1 & 1 & . & . & 2 & 5 & 3 & 3 & 1 & 3 & 2 & . & . & 1 & 3 & 3 & 1 & 3 & 2 & . & 3 & 4 & 1 & 1 & 1 & . & . \\
 . & . & . & . & . & . & 1 & . & . & 1 & 3 & 2 & 2 & 1 & 1 & . & . & 3 & 5 & 4 & 3 & 2 & 4 & 1 & . & . & 1 & 3 & 2 & 2 & 4 & 1 & 1 & 4 & 3 & 2 & 1 & 1 & . \\
 . & . & . & . & . & 1 & 1 & . & . & 1 & 3 & 2 & 2 & 1 & 1 & . & 1 & 3 & 4 & 3 & 3 & 3 & 4 & 1 & 1 & . & 1 & 2 & 1 & 2 & 3 & 1 & 1 & 4 & 2 & 3 & 1 & 1 & . \\
 . & . & . & . & . & . & 1 & 1 & . & . & 2 & 2 & 2 & 1 & 1 & . & . & 1 & 3 & 3 & 3 & 1 & 4 & 2 & 1 & . & . & 1 & 1 & 1 & 3 & 1 & . & 2 & 3 & 3 & 1 & 1 & . \\
 . & . & . & . & . & 1 & 1 & . & . & 2 & 3 & 1 & 1 & 2 & 1 & . & 1 & 3 & 2 & 3 & 1 & 2 & 3 & 1 & 1 & . & 1 & 1 & . & 1 & 1 & . & 2 & 3 & 1 & 1 & 2 & 1 & . \\
 . & . & . & . & . & 1 & 2 & 1 & 1 & 1 & 4 & 3 & 4 & 3 & 3 & 1 & . & 2 & 4 & 4 & 4 & 3 & 7 & 3 & 3 & . & . & 1 & 1 & 1 & 2 & 1 & 1 & 4 & 3 & 4 & 3 & 3 & 1 \\
 . & . & . & . & . & . & 1 & . & 1 & . & 1 & 1 & 2 & 1 & 2 & . & . & . & 1 & 1 & 2 & 1 & 3 & 1 & 2 & . & . & . & . & . & 1 & . & . & 1 & 1 & 2 & 1 & 2 & . \\
 . & . & . & . & . & . & 1 & 1 & 1 & . & 1 & 1 & 1 & 2 & 2 & 1 & . & . & . & 1 & 1 & 1 & 3 & 2 & 2 & . & . & . & . & . & . & . & . & 1 & 1 & 1 & 2 & 2 & 1 \\
 . & . & . & . & . & . & . & . & . & . & . & . & . & . & . & . & 1 & 1 & . & . & . & . & . & . & . & . & 1 & . & . & 1 & . & . & 1 & . & . & . & . & . & . \\
 . & . & . & . & . & . & . & . & . & 1 & . & . & . & . & . & . & 2 & 3 & 1 & 1 & . & 1 & . & . & . & 1 & 3 & 2 & . & 2 & 1 & . & 2 & 2 & . & . & . & . & . \\
 . & . & . & . & . & . & . & . & . & . & 1 & . & . & . & . & . & 1 & 3 & 3 & 2 & 1 & 1 & 1 & . & . & . & 2 & 3 & 2 & 3 & 3 & . & 2 & 3 & 2 & 1 & . & . & . \\
 . & . & . & . & . & . & . & . & . & . & . & 1 & . & . & . & . & . & 1 & 2 & 1 & 1 & . & 1 & . & . & . & . & 2 & 1 & 1 & 2 & 1 & . & 2 & 1 & 1 & . & . & . \\
 . & . & . & . & . & . & . & . & . & 1 & 1 & . & . & . & . & . & 1 & 3 & 2 & 2 & 1 & 1 & 1 & . & . & 1 & 2 & 3 & 1 & 2 & 2 & 1 & 2 & 3 & 1 & 1 & 1 & . & . \\
 . & . & . & . & . & . & . & . & . & . & 1 & 1 & 1 & . & . & . & . & 2 & 4 & 3 & 3 & 1 & 2 & 1 & . & . & 1 & 3 & 2 & 2 & 5 & 1 & 1 & 4 & 3 & 3 & 1 & 1 & . \\
 . & . & . & . & . & . & . & . & . & . & . & . & 1 & . & . & . & . & . & 1 & 1 & 1 & . & 1 & . & . & . & . & . & 1 & 1 & 1 & 1 & . & 1 & 1 & 1 & . & 1 & . \\
 . & . & . & . & . & . & . & . & . & 1 & 1 & . & . & . & . & . & 2 & 3 & 1 & 1 & . & 2 & 1 & . & . & 1 & 2 & 2 & . & 2 & 1 & . & 1 & 2 & . & 1 & 1 & . & . \\
 . & . & . & . & . & . & . & . & . & 1 & 2 & 1 & 1 & 1 & . & . & 1 & 4 & 4 & 4 & 2 & 3 & 4 & 1 & 1 & . & 2 & 3 & 2 & 3 & 4 & 1 & 2 & 5 & 3 & 3 & 2 & 2 & . \\
 . & . & . & . & . & . & . & . & . & . & 1 & 1 & 1 & . & 1 & . & . & 1 & 3 & 2 & 3 & 1 & 3 & 1 & 1 & . & . & 2 & 1 & 1 & 3 & 1 & . & 3 & 2 & 3 & 2 & 1 & . \\
 . & . & . & . & . & . & . & . & . & . & 1 & 1 & 1 & 1 & 1 & . & . & 1 & 2 & 3 & 3 & 1 & 4 & 2 & 1 & . & . & 1 & 1 & 1 & 3 & 1 & 1 & 3 & 3 & 3 & 2 & 2 & 1 \\
 . & . & . & . & . & . & . & . & . & . & 1 & . & 1 & 1 & 1 & . & . & 1 & 1 & 1 & 1 & 2 & 3 & 1 & 2 & . & . & . & . & 1 & 1 & . & 1 & 2 & 2 & 2 & 1 & 2 & 1 \\
 . & . & . & . & . & . & . & . & . & . & . & 1 & 1 & 1 & 1 & 1 & . & . & 1 & 1 & 1 & 1 & 3 & 2 & 2 & . & . & . & . & . & 1 & 1 & . & 2 & 1 & 2 & 2 & 2 & . \\
 . & . & . & . & . & . & . & . & . & . & . & . & . & . & 1 & . & . & . & . & . & . & . & 1 & . & 1 & . & . & . & . & . & . & . & . & . & . & 1 & 1 & . & . \\
\end{array}
\right)
\]
}

 {\tiny
\setlength{\arraycolsep}{1.5pt}
\[
N_{0100}=
\left(
\begin{array}{ccccccccccccccccccccccccccccccccccccccc}
 . & . & . & . & . & . & . & 1 & . & . & . & . & . & . & . & . & . & . & . & . & . & . & . & . & . & . & . & . & . & . & . & . & . & . & . & . & . & . & . \\
 . & . & . & 1 & . & . & 1 & . & 1 & . & . & . & 1 & . & 1 & . & . & . & . & . & . & . & . & . & . & . & . & . & . & . & . & . & . & . & . & . & . & . & . \\
 . & . & 1 & 1 & 1 & . & 2 & 2 & 1 & . & 1 & 1 & 1 & 2 & 1 & 1 & . & . & . & 1 & . & . & 1 & 1 & . & . & . & . & . & . & . & . & . & . & . & . & . & . & . \\
 . & 1 & 1 & 1 & . & 1 & 2 & 1 & 2 & . & 1 & 1 & 1 & 1 & 2 & . & . & . & . & . & 1 & . & 1 & . & 1 & . & . & . & . & . & . & . & . & . & . & . & . & . & . \\
 . & . & 1 & . & 1 & . & 1 & 1 & . & . & . & 1 & . & 1 & . & 1 & . & . & . & . & . & . & . & 1 & . & . & . & . & . & . & . & . & . & . & . & . & . & . & . \\
 . & . & . & 1 & . & 1 & 2 & 1 & 2 & . & 2 & 1 & 3 & 1 & 2 & . & . & 1 & 1 & 1 & 1 & 2 & 2 & . & 1 & . & . & . & . & 1 & . & . & . & 1 & . & 1 & . & . & . \\
 . & 1 & 2 & 2 & 1 & 2 & 5 & 2 & 3 & 1 & 4 & 3 & 4 & 4 & 4 & 1 & . & 1 & 2 & 2 & 2 & 2 & 4 & 2 & 2 & . & . & . & . & . & 1 & . & . & 1 & 1 & 1 & 1 & 1 & . \\
 1 & . & 2 & 1 & 1 & 1 & 2 & 2 & . & 1 & 1 & 2 & 1 & 3 & 1 & 2 & . & . & 1 & 1 & . & . & 2 & 2 & . & . & . & . & . & . & . & 1 & . & . & 1 & . & 1 & . & 1 \\
 . & 1 & 1 & 2 & . & 2 & 3 & . & 3 & . & 2 & 1 & 3 & 1 & 3 & . & . & . & 1 & . & 1 & 1 & 2 & . & 2 & . & . & . & . & . & . & . & . & . & . & 1 & . & 1 & . \\
 . & . & . & . & . & . & 1 & 1 & . & 1 & 2 & 1 & 1 & 2 & 1 & . & . & 2 & 1 & 3 & 1 & 1 & 2 & 1 & . & . & 1 & 1 & . & 1 & 1 & . & 2 & 2 & 1 & . & 1 & . & . \\
 . & . & 1 & 1 & . & 2 & 4 & 1 & 2 & 2 & 6 & 4 & 5 & 4 & 4 & 1 & 1 & 4 & 5 & 5 & 4 & 4 & 7 & 2 & 2 & . & 1 & 2 & 1 & 2 & 3 & 1 & 2 & 5 & 2 & 3 & 2 & 2 & . \\
 . & . & 1 & 1 & 1 & 1 & 3 & 2 & 1 & 1 & 4 & 3 & 3 & 4 & 2 & 1 & . & 2 & 3 & 4 & 2 & 2 & 5 & 3 & 1 & . & . & 1 & 1 & 1 & 2 & 1 & 1 & 2 & 3 & 2 & 2 & 1 & 1 \\
 . & 1 & 1 & 1 & . & 3 & 4 & 1 & 3 & 1 & 5 & 3 & 5 & 2 & 4 & . & 1 & 2 & 4 & 2 & 4 & 4 & 6 & 1 & 4 & . & . & 1 & 1 & 1 & 2 & . & . & 3 & 2 & 4 & 1 & 2 & . \\
 . & . & 2 & 1 & 1 & 1 & 4 & 3 & 1 & 2 & 4 & 4 & 2 & 5 & 3 & 2 & . & 2 & 2 & 4 & 2 & 1 & 5 & 4 & 1 & . & . & 1 & . & . & 1 & 1 & 1 & 2 & 2 & 1 & 3 & 1 & 1 \\
 . & 1 & 1 & 2 & . & 2 & 4 & 1 & 3 & 1 & 4 & 2 & 4 & 3 & 4 & . & . & 1 & 2 & 2 & 3 & 2 & 5 & 1 & 3 & . & . & . & 1 & . & 1 & . & . & 2 & 1 & 2 & 1 & 3 & . \\
 . & . & 1 & . & 1 & . & 1 & 2 & . & . & 1 & 1 & . & 2 & . & 1 & . & . & . & 1 & . & . & 1 & 2 & . & . & . & . & . & . & . & . & . & . & 1 & . & 1 & . & 1 \\
 . & . & . & . & . & . & . & . & . & . & 1 & . & 1 & . & . & . & 1 & 2 & 1 & 1 & . & 2 & 1 & . & . & . & 1 & 1 & . & 2 & 1 & . & 1 & 2 & . & 1 & . & . & . \\
 . & . & . & . & . & 1 & 1 & . & . & 2 & 4 & 2 & 2 & 2 & 1 & . & 2 & 6 & 5 & 5 & 3 & 4 & 5 & 1 & 1 & 1 & 3 & 4 & 2 & 4 & 4 & 1 & 3 & 6 & 3 & 3 & 2 & 1 & . \\
 . & . & . & . & . & 1 & 2 & 1 & 1 & 1 & 5 & 3 & 4 & 2 & 2 & . & 1 & 5 & 7 & 6 & 5 & 5 & 7 & 2 & 2 & . & 2 & 4 & 2 & 4 & 6 & 1 & 2 & 7 & 4 & 5 & 2 & 2 & . \\
 . & . & 1 & . & . & 1 & 2 & 1 & . & 3 & 5 & 4 & 2 & 4 & 2 & 1 & 1 & 5 & 6 & 7 & 4 & 2 & 7 & 4 & 1 & 1 & 2 & 4 & 2 & 2 & 5 & 2 & 3 & 6 & 5 & 3 & 4 & 2 & 1 \\
 . & . & . & 1 & . & 1 & 2 & . & 1 & 1 & 4 & 2 & 4 & 2 & 3 & . & . & 3 & 5 & 4 & 5 & 3 & 6 & 1 & 2 & . & 1 & 2 & 2 & 2 & 4 & 1 & 1 & 5 & 3 & 4 & 2 & 3 & . \\
 . & . & . & . & . & 2 & 2 & . & 1 & 1 & 4 & 2 & 4 & 1 & 2 & . & 2 & 4 & 5 & 2 & 3 & 5 & 5 & 1 & 3 & . & 2 & 2 & 1 & 3 & 3 & . & 1 & 4 & 2 & 4 & 1 & 2 & . \\
 . & . & 1 & 1 & . & 2 & 4 & 2 & 2 & 2 & 7 & 5 & 6 & 5 & 5 & 1 & 1 & 5 & 7 & 7 & 6 & 5 & 11 & 4 & 4 & . & 1 & 3 & 2 & 3 & 5 & 2 & 2 & 7 & 5 & 6 & 4 & 4 & 1 \\
 . & . & 1 & . & 1 & . & 2 & 2 & . & 1 & 2 & 3 & 1 & 4 & 1 & 2 & . & 1 & 2 & 4 & 1 & 1 & 4 & 5 & . & . & . & 1 & . & . & 2 & 1 & 1 & 2 & 3 & 1 & 3 & 1 & 1 \\
 . & . & . & 1 & . & 1 & 2 & . & 2 & . & 2 & 1 & 4 & 1 & 3 & . & . & 1 & 2 & 1 & 2 & 3 & 4 & . & 3 & . & . & . & . & 1 & 1 & . & . & 2 & 1 & 3 & 1 & 2 & . \\
 . & . & . & . & . & . & . & . & . & . & . & . & . & . & . & . & . & 1 & . & 1 & . & . & . & . & . & . & . & . & . & 1 & . & . & 1 & 1 & . & . & . & . & . \\
 . & . & . & . & . & . & . & . & . & 1 & 1 & . & . & . & . & . & 1 & 3 & 2 & 2 & 1 & 2 & 1 & . & . & . & 2 & 2 & . & 2 & 2 & . & 2 & 3 & 1 & 1 & 1 & . & . \\
 . & . & . & . & . & . & . & . & . & 1 & 2 & 1 & 1 & 1 & . & . & 1 & 4 & 4 & 4 & 2 & 2 & 3 & 1 & . & . & 2 & 3 & 2 & 3 & 4 & 1 & 3 & 5 & 3 & 2 & 1 & 1 & . \\
 . & . & . & . & . & . & . & . & . & . & 1 & 1 & 1 & . & 1 & . & . & 2 & 2 & 2 & 2 & 1 & 2 & . & . & . & . & 2 & 1 & 2 & 2 & 1 & 1 & 3 & 1 & 2 & 1 & . & . \\
 . & . & . & . & . & 1 & . & . & . & 1 & 2 & 1 & 1 & . & . & . & 2 & 4 & 4 & 2 & 2 & 3 & 3 & . & 1 & 1 & 2 & 3 & 2 & 3 & 3 & 1 & 1 & 4 & 2 & 3 & 1 & 1 & . \\
 . & . & . & . & . & . & 1 & . & . & 1 & 3 & 2 & 2 & 1 & 1 & . & 1 & 4 & 6 & 5 & 4 & 3 & 5 & 2 & 1 & . & 2 & 4 & 2 & 3 & 6 & 1 & 2 & 6 & 4 & 4 & 2 & 2 & . \\
 . & . & . & . & . & . & . & 1 & . & . & 1 & 1 & . & 1 & . & . & . & 1 & 1 & 2 & 1 & . & 2 & 1 & . & . & . & 1 & 1 & 1 & 1 & 1 & 1 & 1 & 2 & 1 & 1 & . & 1 \\
 . & . & . & . & . & . & . & . & . & 2 & 2 & 1 & . & 1 & . & . & 1 & 3 & 2 & 3 & 1 & 1 & 2 & 1 & . & 1 & 2 & 3 & 1 & 1 & 2 & 1 & 2 & 3 & 1 & 1 & 2 & 1 & . \\
 . & . & . & . & . & 1 & 1 & . & . & 2 & 5 & 2 & 3 & 2 & 2 & . & 2 & 6 & 7 & 6 & 5 & 4 & 7 & 2 & 2 & 1 & 3 & 5 & 3 & 4 & 6 & 1 & 3 & 7 & 5 & 5 & 3 & 3 & 1 \\
 . & . & . & . & . & . & 1 & 1 & . & 1 & 2 & 3 & 2 & 2 & 1 & 1 & . & 3 & 4 & 5 & 3 & 2 & 5 & 3 & 1 & . & 1 & 3 & 1 & 2 & 4 & 2 & 1 & 5 & 3 & 3 & 3 & 1 & . \\
 . & . & . & . & . & 1 & 1 & . & 1 & . & 3 & 2 & 4 & 1 & 2 & . & 1 & 3 & 5 & 3 & 4 & 4 & 6 & 1 & 3 & . & 1 & 2 & 2 & 3 & 4 & 1 & 1 & 5 & 3 & 5 & 1 & 3 & . \\
 . & . & . & . & . & . & 1 & 1 & . & 1 & 2 & 2 & 1 & 3 & 1 & 1 & . & 2 & 2 & 4 & 2 & 1 & 4 & 3 & 1 & . & 1 & 1 & 1 & 1 & 2 & 1 & 2 & 3 & 3 & 1 & 2 & 2 & 1 \\
 . & . & . & . & . & . & 1 & . & 1 & . & 2 & 1 & 2 & 1 & 3 & . & . & 1 & 2 & 2 & 3 & 2 & 4 & 1 & 2 & . & . & 1 & . & 1 & 2 & . & 1 & 3 & 1 & 3 & 2 & 1 & . \\
 . & . & . & . & . & . & . & 1 & . & . & . & 1 & . & 1 & . & 1 & . & . & . & 1 & . & . & 1 & 1 & . & . & . & . & . & . & . & 1 & . & 1 & . & . & 1 & . & . \\
\end{array}
\right)
\]
}

 {\tiny
\setlength{\arraycolsep}{1.5pt}
\[
N_{1000}=
\left(
\begin{array}{ccccccccccccccccccccccccccccccccccccccc}
 . & . & . & . & 1 & . & . & . & . & . & . & . & . & . & . & . & . & . & . & . & . & . & . & . & . & . & . & . & . & . & . & . & . & . & . & . & . & . & . \\
 . & 1 & . & 1 & . & . & . & . & 1 & . & . & . & . & . & . & . & . & . & . & . & . & . & . & . & . & . & . & . & . & . & . & . & . & . & . & . & . & . & . \\
 . & . & 1 & 1 & 1 & . & 1 & 1 & . & . & . & . & . & 1 & . & . & . & . & . & . & . & . & . & . & . & . & . & . & . & . & . & . & . & . & . & . & . & . & . \\
 . & 1 & 1 & 1 & . & . & 1 & . & 1 & . & . & . & . & . & 1 & . & . & . & . & . & . & . & . & . & . & . & . & . & . & . & . & . & . & . & . & . & . & . & . \\
 1 & . & 1 & . & 1 & . & . & 1 & . & . & . & . & . & . & . & 1 & . & . & . & . & . & . & . & . & . & . & . & . & . & . & . & . & . & . & . & . & . & . & . \\
 . & . & . & . & . & 1 & 1 & . & 1 & . & 1 & . & 1 & . & . & . & . & . & . & . & . & 1 & . & . & . & . & . & . & . & . & . & . & . & . & . & . & . & . & . \\
 . & . & 1 & 1 & . & 1 & 2 & 1 & 1 & . & 1 & 1 & 1 & 1 & 1 & . & . & . & . & . & . & . & 1 & . & . & . & . & . & . & . & . & . & . & . & . & . & . & . & . \\
 . & . & 1 & . & 1 & . & 1 & 1 & . & . & . & 1 & . & 1 & . & 1 & . & . & . & . & . & . & . & 1 & . & . & . & . & . & . & . & . & . & . & . & . & . & . & . \\
 . & 1 & . & 1 & . & 1 & 1 & . & 2 & . & . & . & 1 & . & 1 & . & . & . & . & . & . & . & . & . & 1 & . & . & . & . & . & . & . & . & . & . & . & . & . & . \\
 . & . & . & . & . & . & . & . & . & 1 & 1 & . & . & 1 & . & . & . & 1 & . & 1 & . & . & . & . & . & . & . & . & . & . & . & . & 1 & . & . & . & . & . & . \\
 . & . & . & . & . & 1 & 1 & . & . & 1 & 2 & 1 & 1 & 1 & 1 & . & . & 1 & 1 & 1 & 1 & 1 & 1 & . & . & . & . & . & . & . & . & . & . & 1 & . & . & . & . & . \\
 . & . & . & . & . & . & 1 & 1 & . & . & 1 & 1 & 1 & 1 & . & . & . & . & 1 & 1 & . & . & 1 & 1 & . & . & . & . & . & . & . & . & . & . & 1 & . & . & . & . \\
 . & . & . & . & . & 1 & 1 & . & 1 & . & 1 & 1 & 2 & . & 1 & . & . & . & 1 & . & 1 & 1 & 1 & . & 1 & . & . & . & . & . & . & . & . & . & . & 1 & . & . & . \\
 . & . & 1 & . & . & . & 1 & 1 & . & 1 & 1 & 1 & . & 2 & 1 & 1 & . & . & . & 1 & . & . & 1 & 1 & . & . & . & . & . & . & . & . & . & . & . & . & 1 & . & . \\
 . & . & . & 1 & . & . & 1 & . & 1 & . & 1 & . & 1 & 1 & 2 & . & . & . & . & . & 1 & . & 1 & . & 1 & . & . & . & . & . & . & . & . & . & . & . & . & 1 & . \\
 . & . & . & . & 1 & . & . & 1 & . & . & . & . & . & 1 & . & 1 & . & . & . & . & . & . & . & 1 & . & . & . & . & . & . & . & . & . & . & . & . & . & . & 1 \\
 . & . & . & . & . & . & . & . & . & . & . & . & . & . & . & . & 1 & 1 & . & . & . & 1 & . & . & . & . & 1 & . & . & 1 & . & . & . & . & . & . & . & . & . \\
 . & . & . & . & . & . & . & . & . & 1 & 1 & . & . & . & . & . & 1 & 2 & 1 & 1 & . & 1 & 1 & . & . & . & 1 & 1 & . & 1 & 1 & . & 1 & 1 & . & . & . & . & . \\
 . & . & . & . & . & . & . & . & . & . & 1 & 1 & 1 & . & . & . & . & 1 & 2 & 1 & 1 & 1 & 1 & . & . & . & . & 1 & 1 & 1 & 1 & . & . & 1 & 1 & 1 & . & . & . \\
 . & . & . & . & . & . & . & . & . & 1 & 1 & 1 & . & 1 & . & . & . & 1 & 1 & 2 & 1 & . & 1 & 1 & . & . & . & 1 & . & . & 1 & 1 & 1 & 1 & 1 & . & 1 & . & . \\
 . & . & . & . & . & . & . & . & . & . & 1 & . & 1 & . & 1 & . & . & . & 1 & 1 & 2 & . & 1 & . & . & . & . & . & 1 & . & 1 & . & . & 1 & . & 1 & . & 1 & . \\
 . & . & . & . & . & 1 & . & . & . & . & 1 & . & 1 & . & . & . & 1 & 1 & 1 & . & . & 2 & 1 & . & 1 & . & . & . & . & 1 & . & . & . & 1 & . & 1 & . & . & . \\
 . & . & . & . & . & . & 1 & . & . & . & 1 & 1 & 1 & 1 & 1 & . & . & 1 & 1 & 1 & 1 & 1 & 3 & 1 & 1 & . & . & . & . & . & 1 & . & . & 1 & 1 & 1 & 1 & 1 & . \\
 . & . & . & . & . & . & . & 1 & . & . & . & 1 & . & 1 & . & 1 & . & . & . & 1 & . & . & 1 & 2 & . & . & . & . & . & . & . & 1 & . & . & 1 & . & 1 & . & 1 \\
 . & . & . & . & . & . & . & . & 1 & . & . & . & 1 & . & 1 & . & . & . & . & . & . & 1 & 1 & . & 2 & . & . & . & . & . & . & . & . & . & . & 1 & . & 1 & . \\
 . & . & . & . & . & . & . & . & . & . & . & . & . & . & . & . & . & . & . & . & . & . & . & . & . & . & 1 & . & . & . & . & . & 1 & . & . & . & . & . & . \\
 . & . & . & . & . & . & . & . & . & . & . & . & . & . & . & . & 1 & 1 & . & . & . & . & . & . & . & 1 & 1 & 1 & . & 1 & . & . & 1 & 1 & . & . & . & . & . \\
 . & . & . & . & . & . & . & . & . & . & . & . & . & . & . & . & . & 1 & 1 & 1 & . & . & . & . & . & . & 1 & 1 & 1 & 1 & 1 & . & 1 & 1 & 1 & . & . & . & . \\
 . & . & . & . & . & . & . & . & . & . & . & . & . & . & . & . & . & . & 1 & . & 1 & . & . & . & . & . & . & 1 & . & . & 1 & . & . & 1 & . & . & . & . & . \\
 . & . & . & . & . & . & . & . & . & . & . & . & . & . & . & . & 1 & 1 & 1 & . & . & 1 & . & . & . & . & 1 & 1 & . & 1 & 1 & . & . & 1 & . & 1 & . & . & . \\
 . & . & . & . & . & . & . & . & . & . & . & . & . & . & . & . & . & 1 & 1 & 1 & 1 & . & 1 & . & . & . & . & 1 & 1 & 1 & 2 & 1 & . & 1 & 1 & 1 & . & . & . \\
 . & . & . & . & . & . & . & . & . & . & . & . & . & . & . & . & . & . & . & 1 & . & . & . & 1 & . & . & . & . & . & . & 1 & . & . & . & 1 & . & . & . & . \\
 . & . & . & . & . & . & . & . & . & 1 & . & . & . & . & . & . & . & 1 & . & 1 & . & . & . & . & . & 1 & 1 & 1 & . & . & . & . & 1 & 1 & . & . & 1 & . & . \\
 . & . & . & . & . & . & . & . & . & . & 1 & . & . & . & . & . & . & 1 & 1 & 1 & 1 & 1 & 1 & . & . & . & 1 & 1 & 1 & 1 & 1 & . & 1 & 2 & 1 & 1 & 1 & 1 & . \\
 . & . & . & . & . & . & . & . & . & . & . & 1 & . & . & . & . & . & . & 1 & 1 & . & . & 1 & 1 & . & . & . & 1 & . & . & 1 & 1 & . & 1 & 1 & 1 & 1 & . & . \\
 . & . & . & . & . & . & . & . & . & . & . & . & 1 & . & . & . & . & . & 1 & . & 1 & 1 & 1 & . & 1 & . & . & . & . & 1 & 1 & . & . & 1 & 1 & 2 & . & 1 & . \\
 . & . & . & . & . & . & . & . & . & . & . & . & . & 1 & . & . & . & . & . & 1 & . & . & 1 & 1 & . & . & . & . & . & . & . & . & 1 & 1 & 1 & . & 1 & 1 & 1 \\
 . & . & . & . & . & . & . & . & . & . & . & . & . & . & 1 & . & . & . & . & . & 1 & . & 1 & . & 1 & . & . & . & . & . & . & . & . & 1 & . & 1 & 1 & 1 & . \\
 . & . & . & . & . & . & . & . & . & . & . & . & . & . & . & 1 & . & . & . & . & . & . & . & 1 & . & . & . & . & . & . & . & . & . & . & . & . & 1 & . & . \\
\end{array}
\right)
\]
}

\begin{figure}[t]
\centerline{\includegraphics[width=13.8cm]{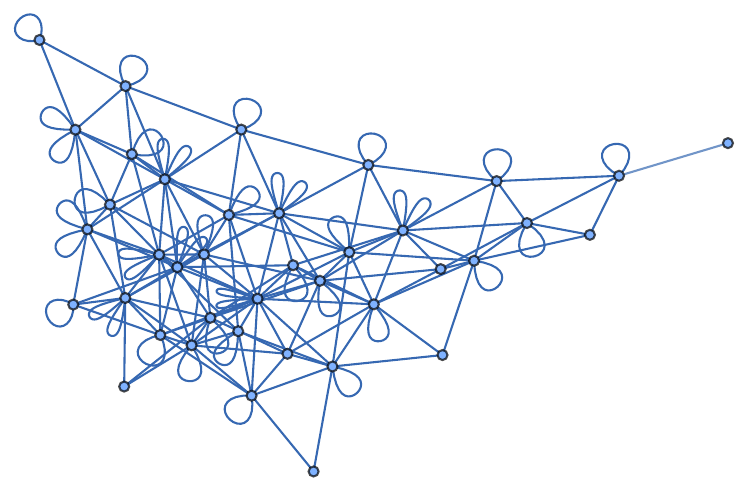}}
\caption{Fusion graph of the simple object labelled by the fundamental irrep $0001$ of $F_4$ at level $6$.}\label{FusionGraph0001}\vspace{-3mm}
\end{figure}

Each fusion matrix can be regarded as the adjacency matrix of a graph, actually a multigraph since multiple edges and self-loops do occur.
In particular, the fundamental fusion matrices $N_f$, with $f$ as above, associated to fundamental simple objects, give rise to four ``fundamental graphs". For instance, the fundamental $0001$, whose associated irrep has the smallest classical dimension ($26$), defines a graph (the adjacency graph of its fusion matrix) that is given by Fig.\ref{FusionGraph0001}.
Its rightmost vertex stands for the trivial $0000$ whose neighbor, which stands for $0001$, is connected to five vertices (including itself), as is clear from the second line of the fusion matrix $N_{0001}$.
  
\subsubsection{Fusion $\ZZ$-polynomials}
\label{FusionZpolynomials}
Finite-dimensional representations of a semisimple Lie algebra decompose into a direct sum of irreducible representations (irreps).
This is so, in particular, for tensor products of representations. 
In turn, arbitrary irreps can be expressed in terms of the fundamental ones as polynomials with integer coefficients.
This property generalizes to the context of fusion categories defined by affine Lie algebras at finite levels (the number of simple objects is finite).
In our case, namely $\mathcal A(F_4,6)$, the fusion ring is a free $Z$-module on the $39$ equivalence classes of simple objects, and it is isomorphic, as a $Z$-algebra, to a quotient of $Z[x,y,z,t]$, where $x,y,z,t$ are ``fundamental objects'', i.e., isomorphism classes of simple objects corresponding to integrable fundamental weights. In plain terms, all fusion matrices can be expressed polynomially in terms of the fundamental fusion matrices. These polynomials, which have integer coefficients, are called ``fusion polynomials''.
For definiteness, we call $x, y, z, t$ the variables respectively associated with the (commuting) fundamental generators $1000$, $0100$, $0010$,  $0001$ and call $P(p)$ the fusion polynomial associated with the simple object $p$.

The decomposition of the tensor product of the irrep $p$ with the fundamental $f$ is given by the lines of fusion matrices $N_f$ (or by their columns since fusion matrices are symmetric).
Selecting for instance $N_{0001}$, which seems to lead to simpler calculations, one can successively multiply each irrep by $0001$ (or equivalently, its associated fusion polynomial by $z$) and decompose.
Straightforward elimination techniques lead then to the list of polynomials  $P(p)$. Example: From line 5 of the fusion matrix $N_{0001}$ we see that non-zero matrix elements (they are equal to $1$) occur in columns 2, 4, and 9. 
Using the ordered list of simple objects, this means that $N_{0001} \times N_{1000} = N_{0001} + N_{0010} + N_{1001}$, equivalently $t \times x = t + z + P(1001)$, hence $P(1,0,0,1) = t(x-1)-z$.

 Warning: the order chosen at the end of the present section  to display these polynomials (left column, top to bottom, then right column), differs from the order that we have chosen for the list of simple objects.
 It only reflects the rather arbitrary sequence of different steps that we have followed to determine them.
 
One should remember that classical fusion polynomials usually differ from those obtained at finite levels,  for two reasons: \\ 
1)  Irreps that enter a classical tensor product decomposition may not exist at a chosen finite level. \\
2) Even if those irreps exist, some quantum intertwiners (for a given triple) may vanish. In other words, even if $a$, $b$, $c$ exist,  the quantum multiplicity of $a \otimes b \rightarrow c$ may be smaller than in the classical case.\\
Here is one example.  The classical (i.e. usual) tensor product of irreps of $F_4$ with highest weights $\{0,0,1,0\}$ and $\{0,0,0,6\}$ decomposes as a direct sum of irreps characterized by the following highest weight components (the last entry of each term is its multiplicity):\\
{\footnotesize 
$
\{ \{0,0,0,5\} , 1 \},
\{ \{0,0,0,6\} , 1 \},
\{ \{0,0,0,7\} , 1 \},
\{ \{0,0,1,3\} , 1 \},
\{ \{0,0,1,4\} , 2 \},
\{ \{0,0,1,5\} , 2 \},
\{ \{0,0,1,6\} , 1 \},
\{ \{0,0,2,3\} , 1 \},
$
$
\{ \{0,1,0,3\} , 1 \},
\{ \{0,1,0,4\} , 1 \},
\{ \{0,1,0,5\} , 1 \},
\{ \{1,0,0,4\} , 1 \},
\{ \{1,0,0,5\} , 2 \},
\{ \{1,0,0,6\} , 1 \},
\{ \{1,0,1,3\} , 1 \},
\{ \{1,0,1,4\} , 1 \}
$
}
However, at level 6 the irreps 
 {\footnotesize $\{0, 0, 0, 7\}, \{0, 0, 1, 5\}, \{0, 0, 1, 6\}, \{0, 0, 2, 3\}, \{0, 1, 0, 4\}, \{0, 1, 0, 5\}, \{1, 0, 0, 5\}, \{1, 0, 0, 6\},$ $\{1, 0, 1, 3\}, \{1, 0, 1, 4\}$}
  do not exist, so that we are left with the candidates  \\
  {\footnotesize $\{\{0,0,0,5\},1\},\{\{0,0,0,6\},1\},\{\{0,0,1,3\},1\}$, $\{\{0,0,1,4\},2\},\{\{0,1,0,3\},1\},\{\{1,0,0,4\},1\}$} in the putative decomposition. 
 Moreover, the multiplicities also differ: the correct decomposition at level $6$ is indeed
  {\footnotesize $\{0,0,1,0\}\otimes\{0,0,0,6\} = \{\{0,0,0,5\},1\}+\{\{0,0,1,3\},1\}+\{\{0,0,1,4\},1\}+\{\{0,1,0,3\},1\}+\{\{1,0,0,4\},1\}$}. 
  This can be checked either by using the fusion polynomials listed below, or by an explicit matrix multiplication of the fusion matrices  $N_{0010}$ and  $N_{0006}$.
  
Here is the list of fusion polynomials for $(F_4,6)$.

{\tiny
\newcommand{\cont}{\\[-0.2ex] &\hspace{-3.2em}}

\noindent
\begin{minipage}[t]{0.48\textwidth}
\begin{align*}
 P(0,0,0,0) &=1,\\
P(1,0,0,0) &= x,\\
P(0,1,0,0) &= y,\\
P(0,0,1,0) &= z,\\
P(0,0,0,1) &= t,\\
P(1,0,0,1) &= t(x-1)-z,\\
P(0,0,0,2) &= t^2-t-x-z-1,\\
P(0,0,1,1) &= -t^2+(-x+z+1)t-y+z+1,\\
P(2,0,0,0) &= -t^2+t+x^2-y+z,\\
P(1,0,1,0) &= zx+x+y-tz,\\
P(0,1,0,1) &= t^2+yt-(x+1)(z+1),\\
P(1,0,0,2) &= -x^2+\left(t^2-z-1\right)x-(t+1)z,\\
P(0,0,0,3) &= t^3-t^2-(x+2z+1)t+x+y,\\
P(0,0,1,2) &= -t^3+(-x+z+1)t^2+(x-y+z+1)t+x^2-z^2-y,\\
P(2,0,0,1) &= -t^3+t^2+\left(x^2-x-y+2z+1\right)t-y+z-x(z+1),\\
P(0,0,2,0) &= -t^3-(x-1)t^2+(x-y+2z+1)t+(z+1)(x+z),\\
P(1,1,0,0) &= t^3+(y-2(z+1))t+(x+1)y-z(z+1),\\
P(0,0,0,4) &= t^4-t^3-(x+3z+1)t^2+(2x+2y+z)t+z^2+x+y+xz+z,\\
P(1,0,1,1) &= (t-x)\left(t^2+(x-z-1)t+y-z-1\right),\\
P(0,1,0,2) &= yt^2+\left(x^2-(z+2)x-2y+z\right)t
\cont +z(z+1)-y(x+z+1),\\
P(2,0,0,2) &= -t^4+t^3+\left(x^2+x-y+3z+1\right)t^2-x(z+2)t-x^3
\cont +x(y-3z-1)+(y-z-1)z-x^2(z+1),\\
P(2,0,1,0) &= (z+1)x^2+\left(-t^2-zt+y+z+1\right)x
\cont -y(t+z-1)+z(t+z),\\
P(0,1,1,0) &= t^3+(x-1)t^2-(zx+x-y+z+1)t+x(y-z-1)+yz,\\
P(0,0,2,1) &= -t^4-(x-1)t^3+(2x-y+2z+1)t^2
\cont +\left(x^2+zx+z^2-y-1\right)t-x^2+y-yz-z-x(z+1),\\
P(1,0,2,0) &= t^4-t^3-\left(x^2+2z+2\right)t^2
\cont +\left(-z^2-z+2y+x(-y+z+1)+1\right)t
\cont +2yz+z+x^2(z+2)+x\left((z+1)^2+y\right),\\
P(3,0,0,0) &= -t^3+(1-2x)t^2+(x-y+3z+1)t+x^3+z^2-2xy-y+2xz+z,\\
P(0,2,0,0) &= (x-1)t^3+(y+z+1)t^2+(-2y+x(y-2z-2)+z+1)t
\cont +y^2-z^2-y-yz-z+x\left(y-(z+1)^2\right),\\
P(1,0,0,3) &= xt^3-zt^2-x(x+2z+1)t+(x+1)y+z(x+z+1),\\
P(1,1,0,1) &= t^4-t^3+(y-z-2)t^2+\left(-z^2-z+x(y+z+1)+1\right)t
\cont -z^2-xy-y-x^2(z+1)+1.
\end{align*}
\end{minipage}
\hfill
\begin{minipage}[t]{0.48\textwidth}
\begin{align*}
P(0,1,1,1) &= -t^4+(x-1)xt^3+\left(x^2+(y+z+1)x-y+2z+2\right)t^2
\cont +\left((y-2z-3)x^2+(-y+z+2)x+(z+1)^2-y\right)t
\cont -(y+1)z+x^2\left(-z^2-3z+y-3\right)
 +x\left(y^2-(z+2)y-2z^2-3z-1\right),\\
P(1,0,1,2) &= t^4-(z+2)t^3-\left(x^2-zx+x-y+z\right)t^2
\cont +\left((z+1)^2+x(-y+2z+2)\right)t
 +x^3-y(z+1)-x\left(z^2+2y\right),\\
P(0,0,2,2) &= xt^4+\left(-x^2+z+2\right)t^3
 -\left(x^3+(y+4z+4)x+y+3z+1\right)t^2
\cont -\left((2y-3z-3)x^2+\left(z^2+4z-3y+1\right)x+z^2-3y+3z+2\right)t
\cont -y^2+2z^2+2y+2yz+2z+x^3(z+2)+x^2\left(2z^2+6z+5\right)
\cont +x\left(-y^2+3(z+1)y+3z^2+6z+1\right)-1,\\
P(0,0,1,3) &= -\left(t^2-t-x-2z\right)\left(t^2+(x-z-1)t+y-z-1\right),\\
P(0,1,0,3) &= t^4+\left(-x^2+2x+y-1\right)t^3
\cont +\left(x^2-(y+3z+6)x-2y-2z-1\right)t^2
\cont +\left((-y+2z+2)x^2-(3z+2)x+z(-y+z+1)\right)t
\cont -x^3+2z^2+2y+yz+3z
\cont +x^2\left((z+2)^2-y\right)+x\left(-y^2+(z+4)y+3z^2+5z+4\right),\\
P(0,0,3,0) &= -t^5+\left(x^2-2x-1\right)t^4
 +\left(x^2+(y+z+1)x-2y+2z+4\right)t^3
\cont +\left(x^3+(y-2(z+2))x^2+(-y+5z+8)x+z^2+y+4z+2\right)t^2
\cont +\left(\left(2y-(z+2)^2\right)x^2+\left(y^2-(z+3)y-z^2+z-1\right)x-2(z+1)^2+y\right)t
\cont -3z^2-y-3z-x^3(z+2)-x^2\left(z^2+4z+y+3\right)
\cont -x\left(z^2+3z+2y(z+1)+1\right),\\
P(1,0,0,4) &= \left(-t^2+z+1\right)x^2+\left(t^4-(3z+2)t^2+(2y+z)t+(z+1)^2+y\right)x
\cont +y-yz+z\left(-t^3+2zt+t+z+1\right),\\
P(0,0,0,5) &= t^5-t^4-(x+4z+2)t^3+(2x+3y+2z+1)t^2
\cont +\left(3z^2+2xz+4z+x-y+1\right)t
 -x^2+z^2-2y-2yz-x(2y+z+1),\\
P(0,0,1,4) &= -t^5+(-x+z+1)t^4+\left(x^2-y+3z+1\right)t^3
\cont +\left(x^2+(y+4z+3)x-z(3z+2)\right)t^2
\cont +\left((y-2z-3)x^2+(-y+z+2)x+3(y-z-1)z\right)t
\cont +z^3-z^2-y+2yz-2z+x^2\left(-z^2-5z+y-4\right)
\cont +x\left(y^2-(z+3)y-2(z+1)^2\right),\\
P(0,0,0,6) &= t^6-t^5-(x+5z+2)t^4+\left(-x^2+4x+4y+3z+1\right)t^3
\cont +\left(6z^2+5z-3y+x(-y+z-3)+1\right)t^2
\cont +\left((-y+2z+1)x^2-(3y+6z+5)x+z^2-4y-5yz\right)t
\cont -z^3+y-yz+2z+x^2\left((z+2)^2-y\right)
+x\left(-y^2+(z+3)y+z^2+3z+2\right).
\end{align*}
\end{minipage}
}

Apart from their academic interest, the explicit knowledge  of these polynomials in four variables will be used at a later stage to provide a computational tool, or a checking device, because 
the relations that they describe hold not only in the fusion ring, but also in any module over this (commutative) ring.
They also provide a reasonably short bookkeeping device to  encode all the fusion matrices in terms of the fundamental ones  $x,y,z,t$.
Needless to say, all these fusion matrices have non-negative integer entries, as they should.

 \subsubsection{About fusion $\QQ$-polynomials}
  In the classical situation, there are infinitely many fusion polynomials and those associated with the fundamental representations are algebraically independent.
At finite levels, however, there are non-trivial relations between powers of fundamental irreps coming ultimately from the fact that we have only a finite number of simple objects.
 In turn, such relations usually allow one to reduce the number of generators.  
In the case of $(F_4, 6)$ for instance, it is relatively easy (see section \ref{FusionQpolynomialsThreeVars}) to express the fusion matrix of $0100$ as a polynomial in the fusion matrices of $1000$, $0010$, $0001$,  but the coefficients of this polynomial are only rational: they are not integers.
Since the already discussed fusion polynomials are $\ZZ$-polynomials in four variables, one can introduce $\QQ$-polynomials that express all simple objects polynomially in terms of only three fundamental ones (rather than four). The rational coefficients that appear in these expressions  usually involve  rather large numerators and denominators.

More generally, one may look for a minimal set of generators, and there is no reason to think that such generators are necessarily to be found among the simple irreps.
As we shall see in the next section, one can obtain a generating  family of $\QQ$-polynomials that depend only on a single variable.

 \subsubsection{Fusion $\QQ$-polynomials (one-variable polynomials)}
   \label{FusionQpolynomials}

We show that all the simples can be written polynomially (with rational coefficients) in terms of a single generator, which we can take as $u = 2 z + t$, where $z=0010$ and $t=0001$.
Existence of such a generator comes from the fact that, over the field $\QQ$ of rational numbers,  the fusion algebra $A_{\QQ} = A \otimes \QQ$  is a finite-dimensional semisimple algebra $\QQ^{39}$.
A primitive element $u$ is such that $A_{\QQ} = \QQ[u]$.
Over the complex numbers, the existence of a primitive element is a consequence of the standard Lagrange interpolation method: given $b=(b_1, b_2, .. b_{39})$ one can find a polynomial $P_b$ such that $P_b(u) = b$. It is known that the Lagrange method also works over the field of rational numbers.
Such an element $u$ should have a minimal polynomial of degree $39$. By inspection, none of the simples $\lambda$ qualify because their degrees (we sort the $\lambda$'s as in \ref{simpleobjects}) are 
{\scriptsize
$$\{1, 37, 34, 27, 28, 37, 15, 36, 30, 18, 33, 21, 28, 23, 28, 20, 29, 14, 37, 24, 22, 26, 15, 27, 30, 35, 36, 34, 26, 23, 16, 19, 28, 26, 25, 28, 23, 35, 34\}$$}
However, just by considering sums of two simples, for instance $N(0001) + 2 N(\lambda)$, the problem is solved. Indeed, the degrees of the minimal polynomials of this family are 
{\scriptsize
$$\{37, 37, 38, 39, 38, 38, 38, 38, 38, 37, 37, 39, 37, 38, 36, 36, 36, 37, 38, 37, 36, 37, 38, 37, 38, 37, 38, 38, 37, 37, 37, 38, 37, 36, 37, 38, 38, 38, 38\}$$}
Thus, we already have two candidates in this list, namely $N(0001) + 2 N(0010)$ and $N(0001) + 2 N(0020)$. For definiteness we stick to the first and call it $u$.

The matrices  $1, u, u^2, u^3, …, u^{38}$ form a basis of the regular representation of the fusion algebra, and the matrices $N(\lambda)$, with irreps $\lambda$ labelled from $1$ to $39$ form another basis.
In order to obtain explicit decompositions, the remaining problem is to find a matrix $K$ that relates the two bases. This is simply solved by evaluating $u$ on constant vectors.

Using this method we can express any simple $\lambda$ as a polynomial in $u$ of degree (at most) $38$. We write $\lambda=P_\lambda(u)$.
Actually, it is enough to use only  the  fundamental irreps since all the other irreps can be expressed in terms of the fundamental ones by using the already obtained $\ZZ$-fusion polynomials.
As before, we call $x,y,z,t$ the fundamental irreps $1000, 0100, 0010, 0001$ but the same symbols can sometimes denote their associated fusion matrices.

In the present framework, we do not necessarily need those $\QQ$-polynomials to obtain arbitrary fusion matrices since we already know all of them, either by brute force calculation using some variant of the Verlinde formula, or because they can be expressed in terms of the fundamental ones by making use of the $\ZZ$-polynomials of the previous section. The $\QQ$-polynomials can be used to check those calculations, but their main interest comes from the fact that they are valid not only in the fusion algebra but also in any representation of the latter, in particular in the modules that we shall consider later,   as well as in their associated ring of quantum symmetries. 

As a matter of fact, it turns out that, in those modules (nimreps), which are obtained as homomorphic images of the Verlinde algebra, the element $t = F(0001)$ is itself primitive, and the above $\QQ$-polynomials (of degree $38$ in $u$) can be replaced by polynomials (in $t$) of smaller degree, namely $11$. This statement does not hold in the fusion algebra itself where polynomials of degree $38$, in $u$, not in $t$, are required.
We shall come back to this discussion later.

The $\QQ$-polynomials in $u$ are unfortunately a bit large to be displayed on paper: their (rational) coefficients involve huge numerators and denominators.
The common denominator of the coefficients of $P_x(u)$, for instance, is an integer with 152 digits. These polynomials can be obtained from the author.

 \subsubsection{Another family of fusion $\QQ$-polynomials (in terms of three fundamental simples)}
    \label{FusionQpolynomialsThreeVars}
In a sense the present subsection is useless because we describe here a ``poor man's method'': We have already obtained $\QQ$-polynomials in  the single variable $u=2z+t$, 
so that exhibiting a new polynomial  in three variables $x,z,t$, which can only be used to return $y$ does not seem particularly useful. 
Nevertheless, apart from the fact that it was the first expression obtained by the author, 
the method is somewhat simpler, and the resulting expression, still huge, can be printed (see below). Finally, one can again use it as a check, either in the fusion algebra itself, or in its modules.

Tensor powers of simple objects (irreps) can be decomposed as linear combinations of  simple objects with non-negative coefficients. 
A simple method to obtain an explicit  decomposition ---let us say of $p^{\otimes s}$---  is to compute the $(s-1)$-th power matrix of the fusion matrix $N_p$. The line of the latter labelled by  $q$ (use the ordered list of irreps) will give the decomposition of $p^{s-1}\otimes q$ into simple objects. In particular it will give the decomposition of $p^{\otimes s}$ if we select the line labelled by $q=p$. For definiteness, let us now take $p=y$.
Then we have a decomposition of all $y^s$, let us say for $s \leq 39$, into simple objects. In turn, those simples can be written as $\ZZ$-polynomials in  $x,y,z,t$.
One can check that this collection of non-trivial polynomial equations can be explicitly solved for $y$. The result is the one given below.
Such a discussion is presented for Lie algebras of type $A$ in the textbook \cite{YellowBook}.

  In order to obtain the $\QQ$-polynomials for all simple objects, it is enough to use the one that expresses $y$ in terms of $x,z,t$ (we make a slight abuse of notation, identifying the objects, or their fusion matrices, with the associated fusion polynomials).  The expression of $y$ is given below.
  The interested reader can replace $x$, $z$, and $t$, in the next expression by the already given fusion matrices of $1000$, $0010$, $0001$, and check that the result is indeed equal to the fusion matrix of $0100$,
  associated with the polynomial called $y$.

{\scriptsize
\[
\newcommand{\Den}{12431163852040700049841800000}
\newcommand{\cont}{\\[-0.2ex] &\hspace*{-7.0em}}
y=\frac{1}{\Den}  \cdot
\begin{aligned}[t]
&
\cont
\cont
\bigg( -1018515905277473446768052154450\,t
-23951545339779825577015862696175\,t^2
\cont
-84970311024359199629345160022350\,t^3
+622532964512539810966726700856975\,t^4
\cont
+2663831010327851280459777959488555\,t^5
-6677548257166532022952925714544595\,t^6
\cont
-28558698530022879319311803195464055\,t^7
+40874206872637699604075349697814955\,t^8
\cont
+151691350176916386129943389049389392\,t^9
-166319001756803011535681224275806289\,t^{10}
\cont
-455047831335782438972706956411841942\,t^{11}
+465451762934185167186192034654113261\,t^{12}
\cont
+817086291230177675631429329743364889\,t^{13}
-880592029281989400399514134265407275\,t^{14}
\cont
-889868350234975960325828952077340289\,t^{15}
+1109991701562692019010720035055928512\,t^{16}
\cont
+559029480642995966926014197833596666\,t^{17}
-925386458071813277145558085392022546\,t^{18}
\cont
-149826542427600555855029734310199837\,t^{19}
+503435704935905031693623370141344520\,t^{20}
\cont
-41388950669796282647368403913709056\,t^{21}
-171922270673069927993519440797969705\,t^{22}
\cont
+48908959826771211737531586967677655\,t^{23}
+32896804842711391608161277817037552\,t^{24}
\cont
-17357566337773450830487375710233339\,t^{25}
-2020066582153203323797173009463752\,t^{26}
\cont
+2957394682284595994834545927102001\,t^{27}
-403171389623085470048510846773941\,t^{28}
\cont
-199455785987899309846970074366299\,t^{29}
+75170826325857915175265765923449\,t^{30}
\cont
-4054606272691878842071670402904\,t^{31}
-2730084186401628516717839434601\,t^{32}
\cont
+735294615879202907243357336246\,t^{33}
-85411727423162430280534630208\,t^{34}
\cont
+4952955915888494128523292525\,t^{35}
-116685533464302882034156729\,t^{36}
\cont
+\Den(2+x+z) \bigg)
\end{aligned}
\]
}
Notice that the last coefficient cancels the denominator, so that the above expression $y$ is the sum of a first degree polynomial in $x$ and $z$ with integer coefficients, and a polynomial in $t$, of degree $36$, with rational coefficients. This expression is easy to solve for $x$, or for $z$, if one prefers to decompose some chosen fusion matrix on another set of generators ---other than the three generators $x,z,t$.

 \section{Modular considerations}
  \label{ModularConsiderations}
 
  \subsection{An exceptional invariant of $F_4 \times A_2$ at level $(6,2)$}
  \label{sec:ZofF4A2}
  
  As we saw previously, the  Lie group $F_4 \times A_2$ at level $(6, 2)$ has central charge $24$.
  
 We already obtained the conformal weights of the integrable irreps of $F_4$ at level $6$.
Those with reduced conformal weight $\overline h$ equal to $0$ are  $0000, 0004, 1100, 0030, 0103$.
 Those with $\overline h=2/5$ are  $0001, 0003,  0021, 0006, 2010$.
 Those with  $\overline h=11/15$ are  $0101, 1012$.
 Those with $\overline h=1/3$ are $0102, 2000$.
 
At level $2$, $A_2$ has $6$ integrable irreps, with associated highest weights $00$, $11$, $10$, $01$, $20$, $02$, and with respective conformal weights equal to $0, 3/5, 4/15, 4/15, 2/3, 2/3$. 

 Now consider the four blocks:
  \begin{eqnarray*}
  Z_{F_4 \times A_2}(0,0) &=& ((0000 + 0004 + 1100+ 0030 ) \otimes (00)) \bigotimes ((0000 + 0004 + 1100+ 0030 ) \otimes (00)) \\
   Z_{F_4 \times A_2}(2/5,3/5) &=& ((0003 +  0021 + 0006 + 2010)  \otimes (11))  \bigotimes  ((0003 +  0021 + 0006 + 2010)  \otimes (11)) \\
    Z_{F_4 \times A_2}(11/15,4/15) &=&  ((0101 + 1012)  \otimes (10 + 01)) \bigotimes  ( (0101 + 1012)  \otimes (10 + 01))\\
     Z_{F_4 \times A_2}(1/3,2/3) &=& ((2000 + 0102 ) \otimes  (02 + 20)) \bigotimes ((2000 + 0102 ) \otimes  (02 + 20))
     \end{eqnarray*}
  and their sum, which defines the bilinear form 
   $$Z_{F_4 \times A_2} =  Z_{F_4 \times A_2}(0,0) +  Z_{F_4 \times A_2}(2/5,3/5) +  Z_{F_4 \times A_2}(11/15,4/15) +  Z_{F_4 \times A_2}(1/3,2/3)$$
       
  The conformal weight, for each of the previous blocks, modulo integers, is $0$.
  Indeed, for $F_4$ and $A_2$, we obtain for the different blocks the pairs  $(0,0), ( 2/5 , 3/5 ), (11/15 , 4/15 ),  (1/3, 2/3)$, 
  whose sums, modulo integers, vanish.
  
    Using the matrix representatives of the modular generators $S, T$ for $(F_4, 6)$ and $(A_2, 2)$,
     one can build those of $F_4 \times A_2$  and check that they commute with the matrix of the previous bilinear form. We obtain therefore a modular invariant for $F_4 \times A_2$ at level $(6,2)$. 

  Notice that the irrep $0103$ has an integer-valued conformal weight ($\overline h=0$) but does not appear in the block $Z_{F_4 \times A_2}(0,0)$, and  that  $0001$ has $\overline h=2/5$ but does not appear in the block $Z_{F_4 \times A_2}(2/5,3/5)$.
  The resulting invariant is obviously not the diagonal one (so it is not of $A$ type), and it is not a simple-current invariant either (so it is not of $D$ type); it is therefore exceptional ($E$ type).

 \subsection{The  exotic exceptional modular invariant of $F_4$ at level $6$}
  \label{modularZ}
 The exceptional modular invariant of $F_4$ at level $6$ can be obtained from the exceptional modular invariant of $F_4 \times A_2$  at level $(6,2)$ described previously in sec.~\ref{sec:ZofF4A2}  by reduction (or contraction) with respect to $A_2$.
 It  is again a sum of four blocks but two of them appear with multiplicity $2$:
 \begin{align*}
  Z_{F_4} =&(0000 + 0004 + 0030 + 1100) \otimes (0000 + 0004 + 0030 + 1100) + \\
 &(0003 + 0006 + 0021 + 2010) \otimes  (0003 + 0006 + 0021 + 2010)  + \\
 & 2  \; (0101 + 1012) \otimes   (0101 + 1012)  + 
  2 \;  (0102 + 2000) \otimes (0102 + 2000)    
 \end{align*}
 
 It will just be called $\mathcal{Z}$ in the sequel.
 The four irreps $0000, 0004, 1100, 0030$ appear in positions $1,10, 24,29$ in the list of simple objects of $(F_4,6)$ given in section \ref{simpleobjects}.
 The four irreps $0003,  0021, 0006, 2010$ appear in positions $6, 19, 26, 38$, the two irreps $0101, 1012$ in positions  $13, 34$ and the two irreps $2000, 0102$ in positions $16,20$.
 As a sparse array of dimension $39 \times 39$, the associated modular invariant matrix $\mathcal{Z}$ can therefore be defined as follows: 
 {\footnotesize
\begin{align*}
   &\{1, 1\}\mapsto 1, \{1, 10\}\mapsto 1, \{1, 24\}\mapsto 1, \{1, 29\}\mapsto 1,
   \{10, 1\}\mapsto 1, \{10, 10\}\mapsto 1, \{10, 24\}\mapsto 1, \{10, 29\}\mapsto 1, \\
   &\{24, 1\}\mapsto 1, \{24, 10\}\mapsto 1, \{24, 24\}\mapsto 1, \{24, 29\}\mapsto 1,
   \{29, 1\}\mapsto 1, \{29, 10\}\mapsto 1, \{29, 24\}\mapsto 1, \{29, 29\}\mapsto 1, \\
  & \{6, 6\}\mapsto 1, \{6, 26\}\mapsto 1, \{6, 19\}\mapsto 1, \{6, 38\}\mapsto 1,
   \{26, 6\}\mapsto 1, \{26, 26\}\mapsto 1, \{26, 19\}\mapsto 1, \{26, 38\}\mapsto 1, \\
  & \{19, 6\}\mapsto 1, \{19, 26\}\mapsto 1, \{19, 19\}\mapsto 1, \{19, 38\}\mapsto 1,
   \{38, 6\}\mapsto 1, \{38, 26\}\mapsto 1, \{38, 19\}\mapsto 1, \{38, 38\}\mapsto 1, \\
  & \{13, 13\}\mapsto 2, \{13, 34\}\mapsto 2,
   \{34, 13\}\mapsto 2, \{34, 34\}\mapsto 2,
   \{20, 20\}\mapsto 2, \{20, 16\}\mapsto 2,
   \{16, 20\}\mapsto 2, \{16, 16\}\mapsto 2
\end{align*}
}
 or, explicitly, as follows (the dots stand for zeroes).
{\tiny
\setlength{\arraycolsep}{1.5pt}
\[
\mathcal{Z}=
 \left(
\begin{array}{ccccccccccccccccccccccccccccccccccccccc}
 1 & . & . & . & . & . & . & . & . & 1 & . & . & . & . & . & . & . & . & . & . & . & . & . & 1 & . & . & . & . & 1 & . & . & . & . & . & . & . & . & . & . \\
 . & . & . & . & . & . & . & . & . & . & . & . & . & . & . & . & . & . & . & . & . & . & . & . & . & . & . & . & . & . & . & . & . & . & . & . & . & . & . \\
 . & . & . & . & . & . & . & . & . & . & . & . & . & . & . & . & . & . & . & . & . & . & . & . & . & . & . & . & . & . & . & . & . & . & . & . & . & . & . \\
 . & . & . & . & . & . & . & . & . & . & . & . & . & . & . & . & . & . & . & . & . & . & . & . & . & . & . & . & . & . & . & . & . & . & . & . & . & . & . \\
 . & . & . & . & . & . & . & . & . & . & . & . & . & . & . & . & . & . & . & . & . & . & . & . & . & . & . & . & . & . & . & . & . & . & . & . & . & . & . \\
 . & . & . & . & . & 1 & . & . & . & . & . & . & . & . & . & . & . & . & 1 & . & . & . & . & . & . & 1 & . & . & . & . & . & . & . & . & . & . & . & 1 & . \\
 . & . & . & . & . & . & . & . & . & . & . & . & . & . & . & . & . & . & . & . & . & . & . & . & . & . & . & . & . & . & . & . & . & . & . & . & . & . & . \\
 . & . & . & . & . & . & . & . & . & . & . & . & . & . & . & . & . & . & . & . & . & . & . & . & . & . & . & . & . & . & . & . & . & . & . & . & . & . & . \\
 . & . & . & . & . & . & . & . & . & . & . & . & . & . & . & . & . & . & . & . & . & . & . & . & . & . & . & . & . & . & . & . & . & . & . & . & . & . & . \\
 1 & . & . & . & . & . & . & . & . & 1 & . & . & . & . & . & . & . & . & . & . & . & . & . & 1 & . & . & . & . & 1 & . & . & . & . & . & . & . & . & . & . \\
 . & . & . & . & . & . & . & . & . & . & . & . & . & . & . & . & . & . & . & . & . & . & . & . & . & . & . & . & . & . & . & . & . & . & . & . & . & . & . \\
 . & . & . & . & . & . & . & . & . & . & . & . & . & . & . & . & . & . & . & . & . & . & . & . & . & . & . & . & . & . & . & . & . & . & . & . & . & . & . \\
 . & . & . & . & . & . & . & . & . & . & . & . & 2 & . & . & . & . & . & . & . & . & . & . & . & . & . & . & . & . & . & . & . & . & 2 & . & . & . & . & . \\
 . & . & . & . & . & . & . & . & . & . & . & . & . & . & . & . & . & . & . & . & . & . & . & . & . & . & . & . & . & . & . & . & . & . & . & . & . & . & . \\
 . & . & . & . & . & . & . & . & . & . & . & . & . & . & . & . & . & . & . & . & . & . & . & . & . & . & . & . & . & . & . & . & . & . & . & . & . & . & . \\
 . & . & . & . & . & . & . & . & . & . & . & . & . & . & . & 2 & . & . & . & 2 & . & . & . & . & . & . & . & . & . & . & . & . & . & . & . & . & . & . & . \\
 . & . & . & . & . & . & . & . & . & . & . & . & . & . & . & . & . & . & . & . & . & . & . & . & . & . & . & . & . & . & . & . & . & . & . & . & . & . & . \\
 . & . & . & . & . & . & . & . & . & . & . & . & . & . & . & . & . & . & . & . & . & . & . & . & . & . & . & . & . & . & . & . & . & . & . & . & . & . & . \\
 . & . & . & . & . & 1 & . & . & . & . & . & . & . & . & . & . & . & . & 1 & . & . & . & . & . & . & 1 & . & . & . & . & . & . & . & . & . & . & . & 1 & . \\
 . & . & . & . & . & . & . & . & . & . & . & . & . & . & . & 2 & . & . & . & 2 & . & . & . & . & . & . & . & . & . & . & . & . & . & . & . & . & . & . & . \\
 . & . & . & . & . & . & . & . & . & . & . & . & . & . & . & . & . & . & . & . & . & . & . & . & . & . & . & . & . & . & . & . & . & . & . & . & . & . & . \\
 . & . & . & . & . & . & . & . & . & . & . & . & . & . & . & . & . & . & . & . & . & . & . & . & . & . & . & . & . & . & . & . & . & . & . & . & . & . & . \\
 . & . & . & . & . & . & . & . & . & . & . & . & . & . & . & . & . & . & . & . & . & . & . & . & . & . & . & . & . & . & . & . & . & . & . & . & . & . & . \\
 1 & . & . & . & . & . & . & . & . & 1 & . & . & . & . & . & . & . & . & . & . & . & . & . & 1 & . & . & . & . & 1 & . & . & . & . & . & . & . & . & . & . \\
 . & . & . & . & . & . & . & . & . & . & . & . & . & . & . & . & . & . & . & . & . & . & . & . & . & . & . & . & . & . & . & . & . & . & . & . & . & . & . \\
 . & . & . & . & . & 1 & . & . & . & . & . & . & . & . & . & . & . & . & 1 & . & . & . & . & . & . & 1 & . & . & . & . & . & . & . & . & . & . & . & 1 & . \\
 . & . & . & . & . & . & . & . & . & . & . & . & . & . & . & . & . & . & . & . & . & . & . & . & . & . & . & . & . & . & . & . & . & . & . & . & . & . & . \\
 . & . & . & . & . & . & . & . & . & . & . & . & . & . & . & . & . & . & . & . & . & . & . & . & . & . & . & . & . & . & . & . & . & . & . & . & . & . & . \\
 1 & . & . & . & . & . & . & . & . & 1 & . & . & . & . & . & . & . & . & . & . & . & . & . & 1 & . & . & . & . & 1 & . & . & . & . & . & . & . & . & . & . \\
 . & . & . & . & . & . & . & . & . & . & . & . & . & . & . & . & . & . & . & . & . & . & . & . & . & . & . & . & . & . & . & . & . & . & . & . & . & . & . \\
 . & . & . & . & . & . & . & . & . & . & . & . & . & . & . & . & . & . & . & . & . & . & . & . & . & . & . & . & . & . & . & . & . & . & . & . & . & . & . \\
 . & . & . & . & . & . & . & . & . & . & . & . & . & . & . & . & . & . & . & . & . & . & . & . & . & . & . & . & . & . & . & . & . & . & . & . & . & . & . \\
 . & . & . & . & . & . & . & . & . & . & . & . & . & . & . & . & . & . & . & . & . & . & . & . & . & . & . & . & . & . & . & . & . & . & . & . & . & . & . \\
 . & . & . & . & . & . & . & . & . & . & . & . & 2 & . & . & . & . & . & . & . & . & . & . & . & . & . & . & . & . & . & . & . & . & 2 & . & . & . & . & . \\
 . & . & . & . & . & . & . & . & . & . & . & . & . & . & . & . & . & . & . & . & . & . & . & . & . & . & . & . & . & . & . & . & . & . & . & . & . & . & . \\
 . & . & . & . & . & . & . & . & . & . & . & . & . & . & . & . & . & . & . & . & . & . & . & . & . & . & . & . & . & . & . & . & . & . & . & . & . & . & . \\
 . & . & . & . & . & . & . & . & . & . & . & . & . & . & . & . & . & . & . & . & . & . & . & . & . & . & . & . & . & . & . & . & . & . & . & . & . & . & . \\
 . & . & . & . & . & 1 & . & . & . & . & . & . & . & . & . & . & . & . & 1 & . & . & . & . & . & . & 1 & . & . & . & . & . & . & . & . & . & . & . & 1 & . \\
 . & . & . & . & . & . & . & . & . & . & . & . & . & . & . & . & . & . & . & . & . & . & . & . & . & . & . & . & . & . & . & . & . & . & . & . & . & . & . \\
\end{array}
\right)
\]
}
 One can check that $\mathcal{Z}$ indeed commutes with the generators $S$ and $T$ of the modular group.

Notice that $\mathcal{Z}$ contains two $2 \times 2$ blocks with coefficients equal to $2$ (they originate from the $A_2$ factors $10+01$ and $20+02$). All the other non-trivial entries are equal to $1$.

The $12=4+4+2+2$ irreps of $(F_4, 6)$ that appear in the previous modular invariant constitute a strict subset of the set of (integrable)  irreps of $(F_4, 6)$, which has cardinality $39$. 
By analogy with the exceptional modular invariants of $SU(2)$, these particular irreps are called (generalized) exponents.
With the chosen ordering of irreps,  the exponent list is  $\{1, 6, 10, 13, 16, 19, 20, 24, 26, 29, 34, 38\}$.

 The above modular invariant was first obtained by Schellekens and Yankielowicz in 1989, \cite{Schellekens-Yankielowicz}.
 Warning: The CERN preprint (1989) contains a typo: the factors of $2$ that should appear in the invariant are incorrectly typed ``1'' in the preprint; they  appear correctly  in the published version (Nucl. Phys. B).
 
 \paragraph{A generalized trigonometric identity.}
 Using the $q$-dimensions of the  irreps of $(F_4, 6)$, and the value of $\vert{\mathcal A}\vert$, given  in sections   \ref{quantumdimensionsofA} and  \ref{globaldimensionofA}, one can check that the following trigonometric identity (whose origin can be traced back to the unitarity of the modular matrix $S$) holds:
  \[
\sum_{m,n}   \text{qdim}(m) \,  \mathcal{Z}_{m,n} \, \text{qdim}(n) = \sum_n  \text{qdim}(n)^2 \,   = \vert {\mathcal A}_k \vert 
\]

  \subsection{Comments on conformal embeddings, VOAs, and Niemeier lattices}
    
The  considerations of sec.~\ref{sec:ZofF4A2} could suggest the existence of a conformal embedding into an ambient Lie algebra whose affine counterpart, at level $1$, would have $c=24$ and would decompose as previously under branching to $F_4 \times A_2$;
those simple objects would have integer conformal weights since they should give $0$ modulo $1$ under this branching.
However, such an ambient Lie algebra $\mathfrak g$ (that one may suppose simple)  cannot exist. This has been known for a long time.
One reason is that the list of conformal embeddings is known, and the case under study just does not show up in the list. 
For completeness, let us explain explicitly  why this is so.

If $\mathfrak g$ is simply laced,  the central charge at level $1$ coincides with the rank of $\mathfrak g$.
Otherwise the rank is multiplied by a corrective factor $x = (\gamma + 1)/(g+1)$, where $g$ is the dual Coxeter number of $\mathfrak g$ and $\gamma$ its Coxeter number.
\\
If $\mathfrak g$ is simple and simply laced, the only two candidates are $A_{24}$ and $D_{24}$, but
neither the first nor the second possesses a Lie subalgebra isomorphic to $F_4 \times A_2$ (in the first case, for instance, it is easy to check that $F_4 \times A_2$ does not possess any faithful representation in dimension $25$).
 If $\mathfrak g$ is simple but not simply laced, there are no candidates at all. Indeed, let us look at the list and calculate the value of $c$: \\
$G_2$ : $(g, \gamma)= (4,6)$. $x= 7/5$. So $c = 14 \times 7/5 = 98/5$.\\
$F_4$ :$(g, \gamma)= (9,12)$. $x = 13/10$. So $c = 52 \times 13/10=338/5$.\\
$B_r$ :  $(g, \gamma)=(2r-1,2r)$. $x = (2r+1)/2r$. So  $c = (2r+1)^2/2$.\\
$C_r$ :$(g, \gamma)= (r+1,2r)$. $x = (2r+1)/(r+2)$. So $c =  r(2r+1)^2/(r+2)$.\\
The equations $(2 r + 1)^2/2 = 24$ and $r (2 r + 1)^2/(r + 2) = 24$ have no integer solutions.

Therefore, the  modular invariant of $F_4 \times A_2$ displayed in sec.~\ref{sec:ZofF4A2}  cannot be obtained from a conformal embedding of Lie algebras.
The associated conformal field theory appears in Schellekens' list of holomorphic CFTs \cite{Schellekens.Meromorphic}  --- it is the case numbered 14.
It was shown (see \cite{LamAndLin} and references therein),
that all 70 Lie algebras in Schellekens’ list can be realized as the weight-one Lie algebras of some holomorphic vertex operator algebras of central charge 24.
In a sense, the notion of conformal embedding (or extension) still holds, provided one uses a vertex operator algebra (VOA) rather than affine Lie algebras.
 In a few words, this construction, in the present case, goes as follows:
 One starts from a Niemeier lattice in dimension $24$ with root lattice $A_8^3$.
The VOA of this lattice is a module for the affine Lie algebra $A_{8,1} \oplus A_{8,1} \oplus A_{8,1}$ (which has central charge $3 \times 8 = 24$). The grade-1 subspace of this lattice VOA is the finite-dimensional Lie algebra $A_8 \oplus A_8 \oplus A_8$. One then considers an outer  \(\mathbb{Z}_{3}\)  automorphism of the VOA that cyclically permutes those three ${A}_{8,1}$ modules, and  builds a fixed-point sub-VOA which is now a module for the affine Lie algebra ${A}_{8,3}$. 
The grade-1 component of this new orbifold VOA is no longer $A_{8}^{3}$ but the smaller finite-dimensional Lie algebra $A_8 \oplus A_2 \oplus A_2$.
A subsequent \(\mathbb{Z}_{2}\) folding  acts on the \(A_{8}\) component to turn it into \(F_{4}\). 
Because it folds the diagram in half, the level of the associated affine Lie algebra doubles from 3 to 6, giving \({F}_{4,6}\). Simultaneously, it twists and merges the two \(A_{2,1}\) components into a single \({A}_{2,2}\). The resulting orbifold VOA is a module for  \({F}_{4,6} \oplus {A}_{2,2}\). \\
Finally, as explained in the previous section, the exotic modular invariant of $(F_4,6)$ that is our starting point in the present paper ---the one discovered in 1989 by Schellekens and Yankielowicz---  is obtained by reduction, or contraction, with respect to $A_2$ from the one obtained for $F_4 \times A_2$.

The present section was included to put our problem in perspective, but, as already mentioned, only the expression of $\mathcal{Z}$, for $F_4$ at level $6$, the one given in sec.~\ref{modularZ}, is needed for our purposes.

 \section{Modular splitting equation}
 \label{modularsplittingeq}

The modular splitting equation is a consequence of the module compatibility equation expressing $\mathcal O$ as an $A_k$-bimodule:
$ (m n) x (p q)=(m (n x p) q)$. It is obtained by writing this equation in terms
of fusion matrices $N_{m}$ and toric matrices $W_{x,y}$, then by choosing $x = 1$. 
Recall that, by definition, $m \, n = \sum_p (N_m)_{np} \, p$, and $ m \, x \, n = \sum_y (W_{x,y})_{mn} \, y$.

Explicitly,
for each of the $39^2$ possible pairs $(m, n)$ this equation reads\footnote{In the present case the fusion matrices are symmetric, and the transpose symbol plays no role.}:
\begin{equation}
N_{m}\,\mathcal{Z}\,N_{n}^{tr}=\sum_{x=1}^{r_O} (W_{1,x})_{mn} \, W_{x,1}
\label{Kmse}
\end{equation}
where
$r_O=Tr(\mathcal{Z}\mathcal{Z}^{\dagger})$ is the dimension of the
quantum symmetry algebra ($r_O=64$).

\paragraph{Determination of the toric matrices (with one twist)}
\label{sec:eqsofmodsplit}
 The matrices $W_{x,1}$, of size $r_{\mathcal A} \times r_{\mathcal A}$, are obtained by solving the above modular splitting equation. One can use the following algorithm,
 that we only sketch here because it was already discussed several times in the literature \cite{RCGSsu4,CTR,IsasiSchieber}.
For each choice of the pair $({m , n})$, \ie for each pair of fusion matrices $(N_{m}, N_{n})$, we define the
matrix ${\mathcal K} _{mn} = N_{m}\,\mathcal{Z}\,N_{n}^{tr}$.
Notice that $\mathcal K$ itself can be regarded as an array with $r_{\mathcal A}^4$ components.
The modular splitting equation can be viewed as the linear expansion of the matrix ${\mathcal K}
_{mn}$ on the set of toric matrices $W_{x,1}$, where the coefficients
of this expansion are the non-negative integers $(W_{1,x})_{mn}$.
Usually, the family of toric matrices is not free: the $r_O$ toric
matrices $W_{x,1}$ are not linearly independent and span, like the matrices
${\mathcal K}_{mn}$, a vector space of dimension $r_W < r_O$. 
In the present case $r_W=40$ and $r_O=64$.

For each $(m,n)$, let us define the quantity $\vert \vert {\mathcal K} _{m n} \vert \vert^2 = \sum_x |(W_{1,x})_{mn}|^2$, We call this quantity, somewhat abusively, the ``squared norm'' of ${\mathcal K} _{m n}$ because it is the sum of the squares of the components of ${\mathcal K} _{m n}$ with respect to the family of toric matrices ${W_{x,1}}$, where $x \in O$.
From the modular splitting equation we see that, for each pair $(m,n)$, the squared norm is equal to the diagonal matrix entry of ${\mathcal K} _{m n}$. Indeed, 
$\vert \vert {\mathcal K} _{m n} \vert
\vert^2 = ({\mathcal K} _{m n})_{\overline{m} \, \overline{n}} = \sum_{p,
q} (N_{m})_{\overline{m} p}\, ({\mathcal{Z}})_{pq} \,(N_{n}^{tr})_{q
\overline{n}}$. In our case, duality is trivial, so $\overline{m}=m$ and $\overline{n}=n$.

The toric matrices $W_{x,1}$ are then obtained by considering matrices ${\mathcal K}_{mn}$
of increasing norms $1,2,3,\ldots$ For example those ${\mathcal K}_{mn}$
 of norm $1$
immediately define toric matrices (since the sum appearing on the right-hand side of the modular splitting equation involves only one term), in particular, one recovers $W_{1,1} =
\mathcal Z$. Then we analyse those ${\mathcal K}_{mn}$ of norm $2$, taking into account the fact that we have already determined several toric matrices in a previous step. 
Then we analyse those of norm $3$ and so on. The process
ultimately stops since the rank $r_W$ is finite. 
 The algorithm is not straightforward because ambiguities arise in the decomposition of ${\mathcal K}_{mn}$ stemming from the fact that the family of toric matrices is not free.
 In particular, those cases where $W_{x,1}=W_{y,1}$ although $x \neq y$ (non-trivial multiplicity) require special attention.

Let us give some details.  $\mathcal K$ has $2078205$ non-zero components. For each $p=(m,n)$, $\mathcal K(p)$ is a vector with $r_{\mathcal A}^2$ components.
 Call ${\mathcal P}(X)$ the $r_{\mathcal A}\times r_{\mathcal A}$ matrix obtained by partitioning the array $X$, with $r_{\mathcal A}^2$ components, into $r_{\mathcal A}$ rows of $r_{\mathcal A}$ elements.
 The possible squared norms are the $\mathcal K(p, p)$. There are $486$ possible norms.
 
 There are $12$ entries $p$ such that $\mathcal K(p, p)=1$.\\
  This happens when $p \in \{1, 2, 4, 5, 40, 41, 43, 44, 118, 119, 157, 158\}$. 
 We obtain in this way $12$ toric matrices ${\mathcal P}(\mathcal K(p))$; in particular ${\mathcal{Z}}={\mathcal P}(\mathcal K(1))$. The rank of this family of matrices is also $12$.
 
There are $14$ entries $p$ such that $\mathcal K(p, p)=2$. \\
This happens when $p \in \{3, 26, 39, 42, 65, 78, 79, 80, 182, 976, 977, 980, 1483, 1484\}$.
 For each such $p$, $\mathcal K(p)$ should be a sum of two toric matrices, so one checks whether
 $\mathcal K(p)$ is independent of the already obtained ones, or not; if it is, one often obtains by subtraction of the known toric matrices (and using the non-negativity constraints), a new toric matrix that we add to the list.
  Otherwise we keep it for later study. After this second step, we obtain a rank-$20$ family of $20$ toric matrices.  

Then one looks at those entries $p$ for which $\mathcal K(p, p)=3$, analyzing the $\mathcal K(p)$ (now a sum of three) along the same lines as above. And so on.
It turns out that it is enough to study those cases where $\mathcal K(p, p)=1, 2, 3, \ldots, 9$, because, after the 9th step, we obtain a family of $40$ independent toric matrices, and we know in advance that the rank should be $r_W=40$. Therefore all the (distinct) toric matrices are then determined.

The final task is to determine the multiplicities, and this can be done by reanalysing the previous constraints \ie the norms ranging from 1 to 9, using the $40$ toric matrices just determined, as well as the constraints obtained from those norms larger than $9$ --- all these constraints, even when not needed, can be used as a check on the resulting family (with multiplicity).

The final result is as follows.
There are forty distinct toric matrices $W_{x,1}$. Eight have multiplicity $3$. Eight have multiplicity $2$. Twenty-four have multiplicity $1$.
 The total is $24+8\times 3+ 8\times 2 =64$, as it should.
 
 Our numbering (from $1$ to $40$) of the distinct toric matrices reflects the relative arbitrariness of the determination process and has no intrinsic significance. Since this ordering will also be used later in this section, we record the multiplicities as follows: 
 members $21, 22, 23, 29, 36, 37, 38, 40$ of the list have multiplicity $3$, and members numbered $24, 25, 26, 27, 31, 33, 34, 39$ have multiplicity $2$. The others have multiplicity $1$.
 Later in this section we shall revert to labels $x$ running from $1$ to $64$, taking the multiplicities of the $W$ matrices into account (\ie we may have $W_x = W_y$ even though $x \neq y$).

 \paragraph{About conjugations and labels}
\label{sec:conjugations}
A general discussion of conjugations can be found in \cite{RCGSsu4}. 
In the present situation, complex conjugation (the ``bar operation'') is trivial: fusion matrices $N_m$ are symmetric.
This leads to the equality $W_{1,x}=W_{x,1}$; we may therefore denote both by $W_x$ without ambiguity. In the present situation we also have $W_{x,y}=W_{y,x}$ for the toric matrices with two twists.
However the toric matrices themselves, $W_{x,y}$ or $W_x$, are not necessarily symmetric. More generally one could define chiral transposition (the operation $c$): $W_{x^c y^c}= W_{x,y}^{tr}$, with, in particular, $W_{x^c}= W_{x}^{tr}$, and the chiral adjoint (the $\dagger$ operation), obtained as the composition of the complex conjugate and the chiral transposition, but this is not needed here, since complex conjugacy is trivial.

One finds $12$ distinct symmetric toric matrices (but $24$ when multiplicities are taken into account).
Using the same labelling as before, these matrices are numbered $1, 6, 14, 18, 23, 25, 26, 33, 37, 38, 39, 40$.
Among them, four (namely 23, 37, 38, 40) have multiplicity $3$, four (namely 25, 26, 33, 39) have multiplicity $2$, and four (namely 1, 6, 14, 18) have multiplicity $1$.

The previous labelling (from $1$ to $40$) that we used to denote distinct toric matrices is not convenient for labelling the simple objects of $\mathcal O$ because the latter should have labels $x$ running from $1$ to $64$.
So, we decide to take multiplicity into account and shift the labels accordingly. For instance, the toric matrix that was in position $21$ of the list, and that was the first to appear with non-trivial multiplicity (namely $3$), will correspond to the $x$-labels $21, 22, 23$. The next toric matrix, which was in position $22$, and also had multiplicity $3$ will then correspond to the $x$-labels $24, 25, 26$. \etc

Using these new labels, from $1$ to $64$, the symmetric toric matrices are those corresponding to the following $x$ labels:
 $1$, $6$, $14$, $18$, $(27, 28 , 29)$, $(32 , 33)$,
$(34, 35)$,
$(46, 47)$,
$(54, 55, 56)$,
$(57, 58, 59)$,
$(60, 61)$,
$(62, 63, 64)$.
The labels appearing within the same parentheses denote simple objects of $\mathcal O$ that define equal toric matrices $W_{x}$ (multiplicities).

The ambichiral toric matrices are symmetric. As we shall see in the next section, they correspond to labels $J = \{1, 18, 32, 33, 60, 61 \}$. Among these six toric matrices, only four are distinct.

Introducing the characters $\chi$, we can write $Z_{x,y}= \overline \chi W_{x,y} \chi$.
The quantities $Z_{x,y}$ are interpreted as twisted partition functions for a boundary conformal field theory in the presence of defect lines of type $x$ and $y$.
We also set $Z_{x}= Z_{x,1}$. The usual partition function is the modular invariant $Z_{1}$.

 \paragraph{Explicit results}
 \label{bilinears}
 
It would not be reasonable to print the $64$ toric matrices $W_x$ of size $39\times 39$ (or even the $40$ matrices that are distinct). The interested reader can download them as Mathematica files from the author's web page. 
Since the matrices $W$ define bilinear forms, one could restrict our attention to those that are symmetric and print the associated quadratic forms (the matrices $W$ being rather sparse, the latter expressions are reasonably short), but to save space we shall only display those associated with the ambichiral toric matrices:

Define
\begin{align*}
U &=\chi (0,0,0,3)+\chi (0,0,0,6)+\chi (0,0,2,1)+\chi (2,0,1,0) \\
V &= \chi (0,0,0,0)+\chi (0,0,0,4)+\chi (0,0,3,0)+\chi (1,1,0,0) \\
A &=\chi (0,1,0,1)+\chi (1,0,1,2) \\
B &= \chi (0,1,0,2)+\chi (2,0,0,0)
 \end{align*}
Then,
 \begin{align*}
 Z_{x=1}=& \overline V V + \overline U U + 2 (\overline A A + \overline B B), \text{ (this is the modular invariant $\mathcal{Z} = Z_{F_4}$ again)}\\
 Z_{x=18}=& \overline U U + ( \overline U V + \overline V U) + 2 ( \overline A A + \overline A B + \overline B A),\\
 Z_{x=32}= Z_{x=33}=& \overline U A + \overline A U + \overline A A + \overline V B + \overline B V, \\
 Z_{x=60}= Z_{x=61}=&  \overline U (A + B) + ( \overline A + \overline B) U + ( \overline V + \overline B) A + \overline A (V + B).
 \end{align*}
Equivalently, 
 \[ 
 Z_{x=1}= \begin{pmatrix}
\overline U & \overline V & \overline A &\overline B 
 \end{pmatrix}.
 \begin{pmatrix}
 1 & 0 & 0 & 0 \\
 0 & 1 & 0 & 0 \\
 0 & 0 & 2 & 0 \\
 0 & 0 & 0 & 2
 \end{pmatrix}.
 \begin{pmatrix}
 U \\ V \\ A \\ B 
 \end{pmatrix},
\]

  \[ 
 Z_{x=18}= \begin{pmatrix}
\overline U & \overline V & \overline A &\overline B 
 \end{pmatrix} .
 \begin{pmatrix}
 1 & 1 & 0 & 0 \\
 1 & 0 & 0 & 0 \\
 0 & 0 & 2 & 2 \\
 0 & 0 & 2 & 0
 \end{pmatrix}.
 \begin{pmatrix}
 U \\ V \\ A \\ B 
 \end{pmatrix},
\]

  \[ 
  Z_{x=32}= Z_{x=33}=
 \begin{pmatrix}
\overline U & \overline V & \overline A &\overline B 
 \end{pmatrix}
 . \begin{pmatrix}
 0 & 0 & 1 & 0 \\
 0 & 0 & 0 & 1 \\
 1 & 0 & 1 & 0 \\
 0 & 1 & 0 & 0
 \end{pmatrix}.
 \begin{pmatrix}
 U \\ V \\ A \\ B 
 \end{pmatrix},
\]

  \[ 
  Z_{x=60}= Z_{x=61}=
 \begin{pmatrix}
\overline U & \overline V & \overline A &\overline B 
 \end{pmatrix} .
 \begin{pmatrix}
 0 & 0 & 1 & 1 \\
 0 & 0 & 1 & 0 \\
 1 & 1 & 0 & 1 \\
 1 & 0 & 1 & 0
 \end{pmatrix}.
 \begin{pmatrix}
 U \\ V \\ A \\ B 
 \end{pmatrix}.
\]

 \section{Left and right intertwiners}
     \label{Intertwiners}
    \subsection{Intertwiners}
    
    For each triplet $(m,n,x)$, where $m,n$ denote irreps of $A$ and $x$ an irrep of $O$, we define the matrices $$K^x_{mn} = N_m W_{x,1} N_n$$
    The intertwining equations (one matrix equation for each triplet) read $$K^x_{mn} =  \sum_{y=1}^{r_O} (W_{x,y})_{mn} \, W_{y,1}$$
    The modular splitting equation studied in the previous section was obtained as the special case $x=1$ (the left-hand side was known from the a priori knowledge of the modular invariant matrix).
    Now the situation is similar since we have determined the left-hand side (the matrices $W_{x,1}$). 
    Solving these equations for each triplet $(m,n,x)$ would allow one to determine all the toric matrices with two twists $W_{x,y}$, or, equivalently,  the entire collection of intertwiners $V_{m,n}$, 
    since  $(V_{m,n})_{x,y} =(W_{x,y})_{m,n}$. 
     It is technically simpler to restrict one's attention to the cases where $x$ is arbitrary, $m$ is chosen to be a fundamental irrep $f$, and $n$ is trivial ($n = 1$), or conversely.
    Indeed,  the determination of the coefficients $({W_{x,y}})_{f,1}$ and $({W_{x,y}})_{1,f}$ of the fundamental intertwiners $V_f^L=V_{f,1}$ and $V_f^R=V_{1,f}$ allows one to obtain the fundamental left and right Ocneanu graphs, and to identify the quantum subgroup $M$ (and the quantum module $M^\prime$, if any), together with their collection of annular matrices (nimreps). 
    
    In practice, we calculate $K^x_{f,1}$, for the four fundamental irreps $f$, using our knowledge of the fusion matrices and of the toric matrices with a single twist.
    Then, we identify the result with the right-hand side of the intertwining equations, written  $\sum_{y=1}^{r_O} (V_{f}^L)_{x,y}  \, W_{y}$, where  $V_{f}^L$ is a matrix $64 \times 64$ with indeterminate coefficients.
    In other words, for each $f$ (hence four families of indeterminate coefficients), we solve the equations  $N_f W_{x} = \sum_{y} (V_{f}^L)_{x,y}  \, W_{y}$.
    We then do the same for $K^x_{1,f}$ and solve for $V_{f}^R$.
    
    Equivalently, we may introduce a  $64 \times (39 \times 39)$ array $K_0$, a table of $64$ one-dimensional arrays of length $39^2$, where the row $K_0(x) = Vec(W_x)$ is obtained by flattening the matrix $W_x $, i.e. by  stacking the rows of $W_x = (W_x)_{m,n}$ side by side,  
   and rewrite the previous left and right intertwining equations ($T$ standing for transpose) as follows:  $$(N_f \otimes N_1) . K_0^T  = K_0^T . V_{f}^L$$ and $$(N_1 \otimes N_f) . K_0^T  = K_0^T . V_{f}^R$$

    Actually, those intertwining equations do not determine all the indeterminate coefficients, and we have to rely on other types of constraints in order to solve for $V_{f}^L$ and $V_{f}^R$.    
  
    \subsection{About $V_{n}^L$ and $V_{n}^R$: technical details}

    Most parameters occurring in the matrix elements of the left and right intertwiners can be fixed from the following considerations: 
    \begin{itemize} 
     \item  The intertwiners should solve the previous  intertwining equations.
  \item  The matrices $V_{n}^L$ are related to the matrices $V_{n}^R$ by a permutation (which is independent of the choice of $n$).
 \item    The coefficients of matrices $V_{n}^L$ (and of $V_{n}^R$) should be non-negative integers.
\item    The matrices $V_{n}^L$  and $V_{n}^R$ should be symmetric (not a generic feature, but one that should hold in the present case).
    \end{itemize}
    The matrices $V_{n}^L$  (or $V_{n}^R$) can be obtained from the fundamental ones $V_{f}^L$ using fusion polynomials of $F_4$ at level $6$.
  It is therefore enough to restrict our attention to the fundamental ones, $V_{f}^L$ and  $V_{f}^R$. 
Taking into account the previous constraints, we still have many undetermined parameters, even when we restrict our attention to a single fundamental irrep, for instance $t=0001$.
 Imposing that  matrices $V_{n}^L$ (or $V_{n}^R$), for different irreps $n$, commute, is also a constraint that should be obeyed, but it does not help much, at this stage, because it introduces too many quadratic relations between unknown parameters.

We shall see later (in sections \ref{annularmatricesM} and \ref{annularmatricesMprime}) how these remaining parameters are determined. In the present section, however, we use the results obtained there and continue our general description of $V^L_f$ and $V^R_f$. Let us only mention here that determining the remaining parameters is relatively easy for $M$ but much more involved for $M^\prime$.

 \subsection{General features of $V^L_f$ and $V^R_f$}
    
  For each $f_i$ (i=1,2,3,4), the adjacency graph of  $V^L_f$ has four connected components forming two pairs of isomorphic components. 
  
The above fact can be recognised after solving only the  ``first round of constraints'' discussed previously,  i.e. even before determining all the parameters entering the left and right intertwiners $V^L_f$ and $V^R_f$.
 
    The graph of $V^L_f$ (for each $f$) contains two isomorphic copies of a ``quantum subgroup'' $M$, called $M_1^L$, $M_2^L$, and two isomorphic copies of another module $M^\prime$, called $M_1^{\prime L}$, $M_2^{\prime L}$.
    Each connected component contains $16$ vertices.
    Using the labeling described previously, the four components are as follows:
    
$M_1^{L} : \quad	 {1, 5, 9, 11, 16, 18, 20, 32, 33, 36, 37, 43, 44, 45, 60, 61}$

$M_2^{L} : \quad	{2, 6, 10, 12, 13, 14, 15, 30, 31, 34, 35, 42, 46, 47, 48, 49}$

$M_1^{\prime L} : \quad {3, 7, 17, 19, 24, 25, 26, 51, 52, 53, 54, 55, 56, 57, 58, 59}$

$M_2^{\prime L}  : \quad {4, 8, 21, 22, 23, 27, 28, 29, 38, 39, 40, 41, 50, 62, 63, 64}$
    
The graph of $V^R_f$ (for each $f$) displays the same features.
 
$M_1^{R} : \quad {1, 2, 3, 4, 13, 17, 18, 30, 31, 32, 33, 48, 49, 50, 60, 61}$

$M_2^{R} : \quad{5, 6, 7, 8, 14, 16, 19, 34, 35, 36, 37, 38, 43, 44, 46, 47}$

$M_1^{\prime R} :\quad {9, 10, 15, 20, 21, 22, 23, 39, 40, 41, 54, 55, 56, 57, 58, 59}$

$M_2^{\prime R} : \quad  {11, 12, 24, 25, 26, 27, 28, 29, 42, 45, 51, 52, 53, 62, 63, 64}$

For each fundamental $f$, the intersection of the adjacency graphs of $V^L_f$ and $V^R_f$ is a subgraph with 6 vertices, which are the same for all the $f$'s, and correspond to the simple objects of the ambichiral (dyslectic) category $\mathcal J$. Recall that $\mathcal J$ is modular, although $\mathcal O$ is not.
With the same labeling as before, those vertices (called ambichiral points) are  $\{1, 18, 32, 33, 60, 61\}$.

For $f=1000$, $0100$, $0010$, $0001$, and therefore also for all irreps, one has
$$V_{f}^R = \pi_{LR}^{-1} \,.\, V_{f}^L \,.\, \pi_{LR}$$
where $\pi_{LR}$ is the permutation matrix defined by the permutation
$\varpi_{LR} =$ {(1, 5, 9, 11, 2, 6, 10, 12, 3, 7, 4, 8, 16, 14, 19, 13, 20, 18, 15, \
17, 24, 25, 26, 21, 22, 23, 27, 28, 29, 43, 44, 32, 33, 34, 35, 48, \
49, 42, 51, 52, 53, 38, 30, 31, 50, 46, 47, 36, 37, 45, 39, 40, 41, \
54, 55, 56, 57, 58, 59, 60, 61, 62, 63, 64)}.
  
With a slight abuse of terminology, one can say that the connected component of the graph that contains the unit of $O$ defines ``the graph'' of a quantum subgroup $M$. Strictly speaking, this statement is incorrect since one can define one left graph and one right graph, for each irrep $n$ ($39$ choices) of the fusion algebra of $F_4$ at level $6$. 
In a later section we shall display the restriction of $V^L_f$ (or  of $V^R_f$) to the connected component of the identity, i.e. the nimreps of the quantum subgroup $M$, for the four fundamental irreps $f$.
For definiteness we identify $M$ with $M_1^{L}$.
Then we shall do the same for the quantum module $M^\prime$ that, for definiteness, we identify with $M_1^{\prime L}$.

     \subsection{The Ocneanu graphs of $V^L_f$ and $V^R_f$}
     \label{OcneanuGraphs}
           
It would not be reasonable to print the $2 \times 39$ matrices of size $64 \times 64$, nor would it be reasonable to display the entire family of adjacency graphs (which are not connected). Even restricting our attention to the four fundamental irreps would be somewhat cumbersome and displaying their full adjacency graphs is not very helpful because many non-zero matrix elements are larger than $1$ (multiplicity), a feature that tends to clutter the figures.
 
\begin{figure}
\centerline{\scalebox{0.7}{\includegraphics{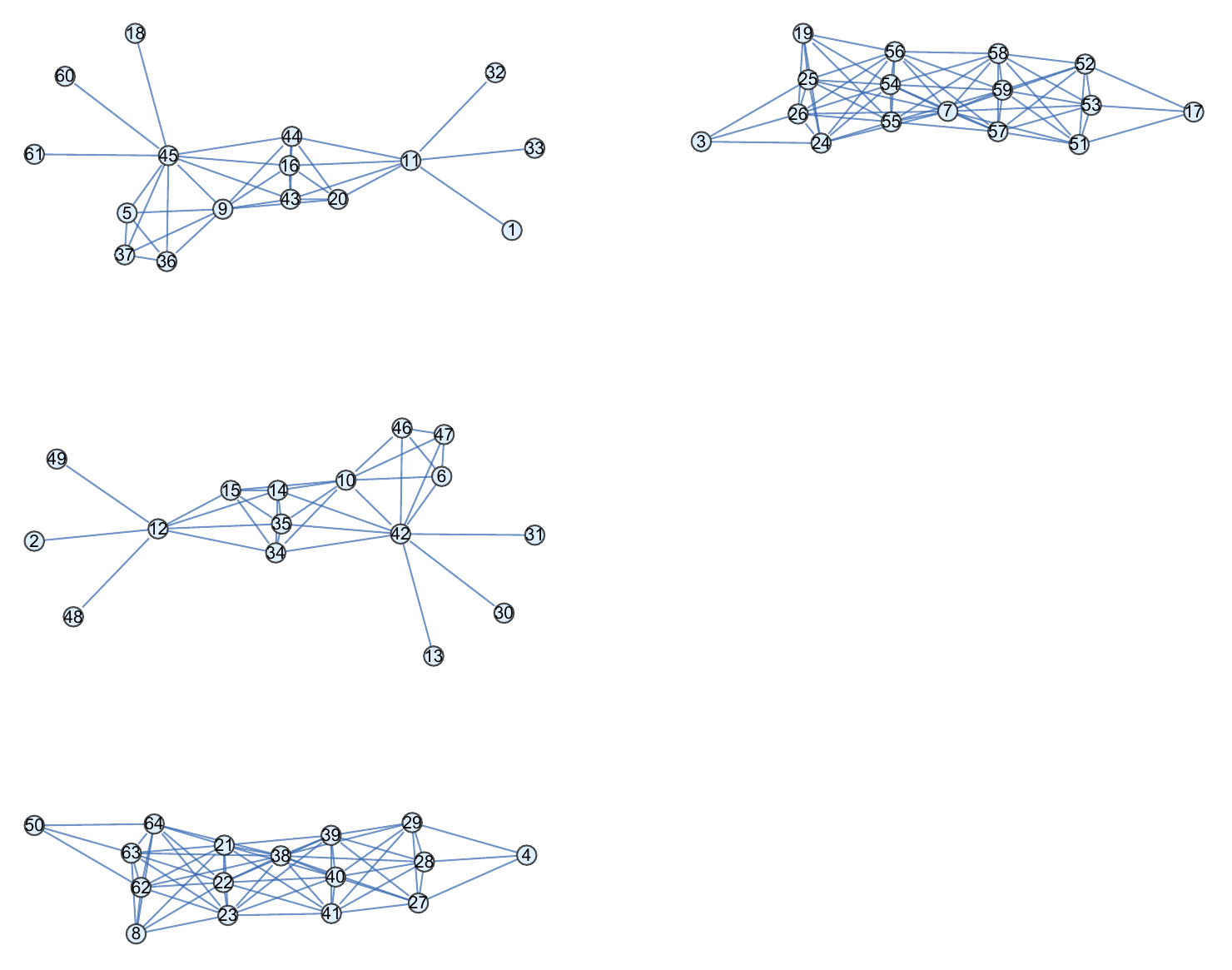}}} 
\caption{The left chiral graph of quantum symmetries for $1000$ (adjacency graph of $V^L_{1000}$). Connected components: two copies of the quantum subgroup $M$ and two copies of the quantum module $M^\prime$. Multiplicities and self-loops are not displayed. Ambichiral points: $\{1, 18, 32, 33, 60, 61\}$.}
\label{FL1000}
\end{figure}

\begin{figure}
\centerline{\scalebox{0.7}{\includegraphics{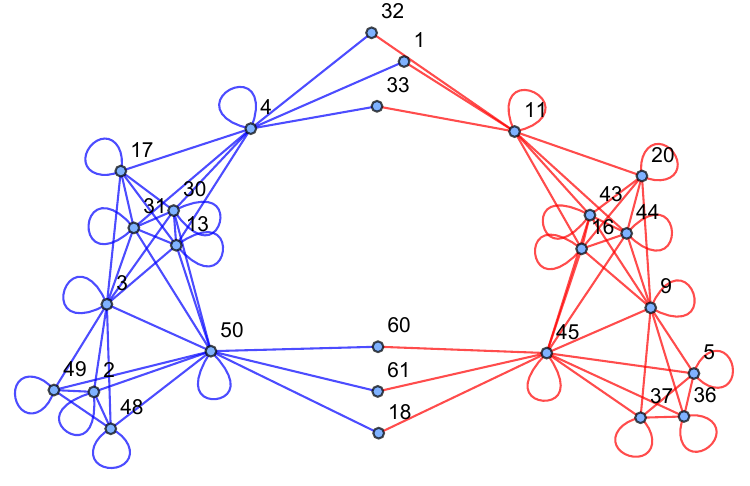}}} 
\caption{The left and right chiral graphs of quantum symmetries for $1000$ restricted to the  connected component $M$ (restrictions of the adjacency graphs of $V^L_{1000}$ and $V^R_{1000}$).  Left and right graphs intersect at the ambichiral points: $\{1, 18, 32, 33, 60, 61\}$. Multiplicities are not displayed.}
\label{FLR1000M}
\end{figure}

   We display the left chiral graph of quantum symmetries for the irrep $1000$ on Fig. \ref{FL1000}. Notice its four connected components forming two pairs of equivalent components and defining the graphs of $M$ and $M^\prime$ (for that irrep).
   The right chiral graph for the same irrep $1000$ looks the same. It is usually instructive to plot simultaneously the left and right chiral graphs, but in our case the picture would not be very clear, even if we discard self-loops and multiedges, because of the large number of vertices and edges. However, we plot simultaneously in Fig.  \ref{FLR1000M} the restriction of both $V^L_{1000}$ and $V^R_{1000}$ to their first connected component called $M$.
      Notice that left and right chiral graphs intersect at the ambichiral vertices.
   
\begin{figure}
\centerline{\scalebox{0.7}{\includegraphics{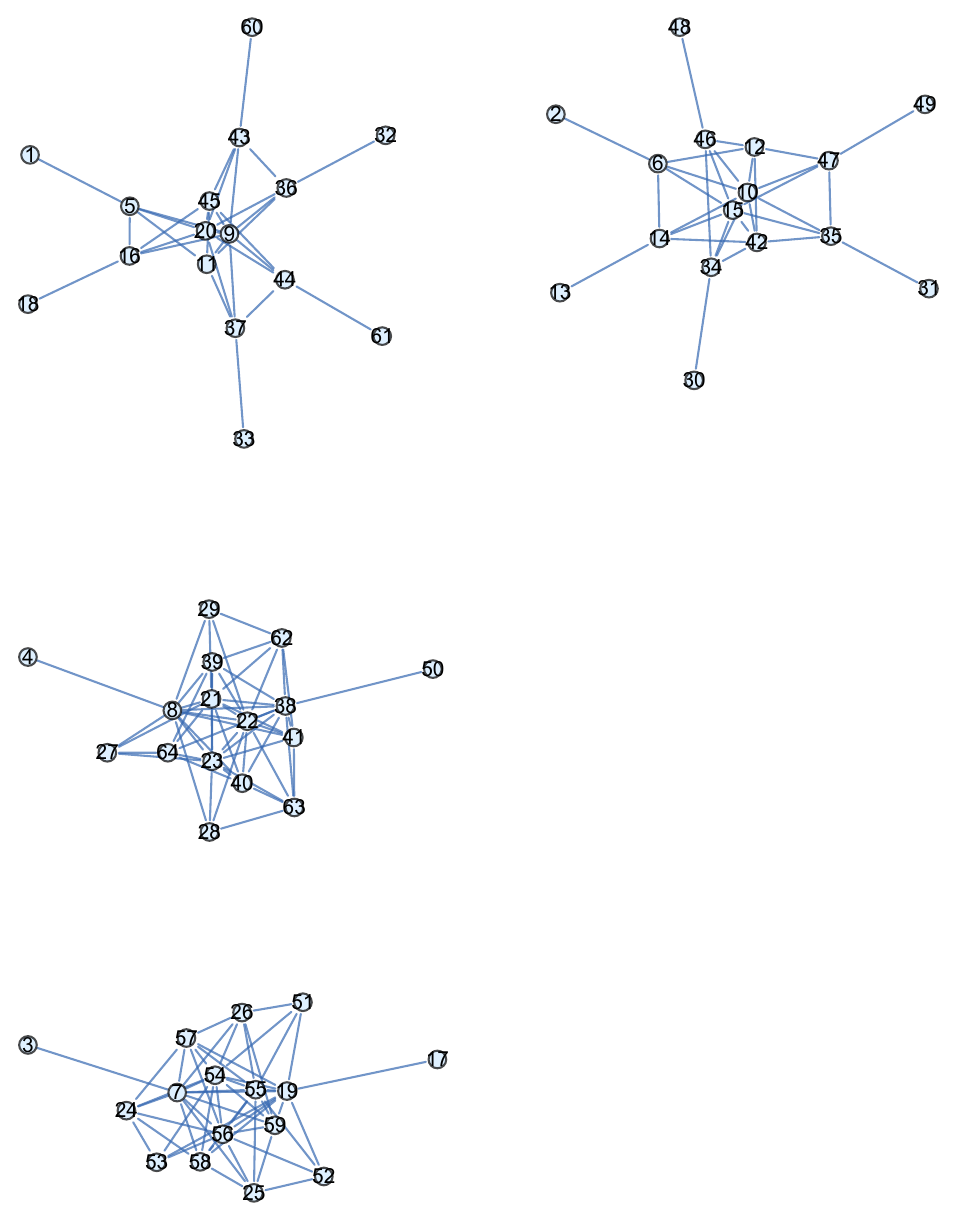}}} 
\caption{The left chiral graph of quantum symmetries for $0001$ (adjacency graph of $V^L_{0001}$). Connected components: two copies of the quantum subgroup $M$ and two copies of the quantum module $M^\prime$. Multiplicities and self-loops are not displayed. Ambichiral points: $\{1, 18, 32, 33, 60, 61\}$.}
\label{FL0001}
\end{figure}

\begin{figure}
\centerline{\scalebox{0.7}{\includegraphics{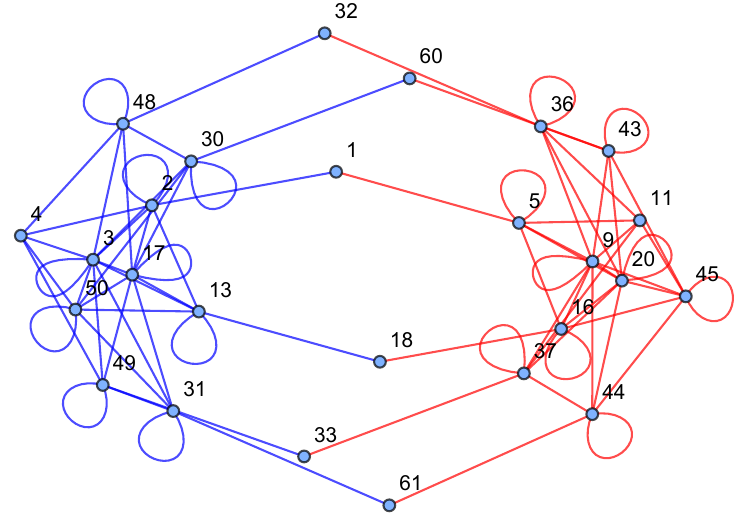}}} 
\caption{The left and right chiral graphs of quantum symmetries for $0001$ restricted to the  connected component $M$ (restrictions of the adjacency graphs of $V^L_{0001}$ and $V^R_{0001}$).  Left and right graphs intersect at the ambichiral points: $\{1, 18, 32, 33, 60, 61\}$. Multiplicities are not displayed.}
\label{FLR0001M}
\end{figure}

   In the same way, we display the left chiral graph of quantum symmetries for the irrep $0001$ in Fig. \ref{FL0001}, and plot simultaneously in Fig. \ref{FLR0001M} the restriction of both $V^L_{0001}$ and $V^R_{0001}$ to their first connected component.
   
   The adjacency graphs of the left and right intertwiners relative to the irreps $0100$ and $0010$ are not displayed in this section, but their restrictions to their first and third connected components, defining the quantum subgroup $M$ and the quantum module $M^\prime$ will be given in the next section.

   \section{The quantum subgroup $M$}
      \label{qModuleM}
 
      \subsection{The quantum subgroup $M$ as a module over the fusion ring}
         \label{annularmatricesM}
   
The following four matrices $F_{0001}$, $ F_{0010}$, $ F_{0100}$, $ F_{1000}$ are obtained as the restriction of the left intertwining matrices  $V^L_f$ to the connected component of the identity $M_1^{L}$. 
These four matrices are particular nimreps (or annular matrices) for the quantum subgroup $M$. 
More generally, the annular matrices (or nimreps) representing the fusion matrices $N_n$ are denoted $F_n$, with $n=1,2,\ldots 39$, or, with explicit Dynkin labels, as $F_{ijkl}$, with the labels sorted as in section \ref{simpleobjects}, for instance $F_5=F_{1000}$.

After solving the various constraints described at the beginning of the previous section, one finds that the matrix elements of $F_{0001}$, for instance, contain five undetermined parameters that have to obey a few inequalities in order for the matrix elements to be non-negative integers. It is easy to find a solution.  One can make similar comments for $F_{0010}$. 
The two solutions (choices of parameters) should be compatible, meaning that the annular matrices $F_\lambda$ obtained from them by using the $\ZZ$ or $\QQ$ polynomials discussed in a previous section should also give matrices with non-negative integers. Since the number of parameters is relatively small, it is easy to find a solution, and one can convince oneself that the existence of several solutions reflects the existence of trivial symmetries (permutations of particular  basis elements).
   
The adjacency graphs of the four ``fundamental'' annular matrices $F_f$ are given below, to the right. We have enumerated their vertices from $1$ to $16$, these labels corresponding to the natural order of the rows or columns of the chosen restriction. However,  the reader should be warned that these labels do not coincide with the labelling used for the $64$ vertices of $V^L_f$. 
For instance, the six ambichiral vertices, which are called $\{1, 18, 32, 33, 60, 61 \}$ in the graphs of $V^{L,R}_f$ appear as $(1,6,8,9,15,16)$ in the graphs of $F_f$ --- the latter labels denoting the positions of these special vertices in the connected component of the identity $M_1^{L}$. Notice that they are displayed ``around the others'' in the following figures.
\begin{center}
\begin{minipage}[t]{0.62\textwidth}
\vspace{0pt}
{\footnotesize
\[
\setlength{\arraycolsep}{2.5pt}
 F_{0001}=
 \left(
\begin{array}{cccccccccccccccc}
 . & 1 & . & . & . & . & . & . & . & . & . & . & . & . & . & . \\
 1 & 1 & 1 & 1 & 1 & . & 1 & . & . & . & . & . & . & . & . & . \\
 . & 1 & 4 & 2 & 2 & . & 3 & . & . & 1 & 1 & 2 & 2 & 2 & . & . \\
 . & 1 & 2 & . & . & . & . & . & . & 1 & 1 & . & . & 1 & . & . \\
 . & 1 & 2 & . & 2 & 1 & 1 & . & . & . & . & . & . & 1 & . & . \\
 . & . & . & . & 1 & . & . & . & . & . & . & . & . & . & . & . \\
 . & 1 & 3 & . & 1 & . & 1 & . & . & 1 & 1 & 1 & 1 & 2 & . & . \\
 . & . & . & . & . & . & . & . & . & 1 & . & . & . & . & . & . \\
 . & . & . & . & . & . & . & . & . & . & 1 & . & . & . & . & . \\
 . & . & 1 & 1 & . & . & 1 & 1 & . & 1 & . & 1 & . & . & . & . \\
 . & . & 1 & 1 & . & . & 1 & . & 1 & . & 1 & . & 1 & . & . & . \\
 . & . & 2 & . & . & . & 1 & . & . & 1 & . & 2 & . & 1 & 1 & . \\
 . & . & 2 & . & . & . & 1 & . & . & . & 1 & . & 2 & 1 & . & 1 \\
 . & . & 2 & 1 & 1 & . & 2 & . & . & . & . & 1 & 1 & 1 & . & . \\
 . & . & . & . & . & . & . & . & . & . & . & 1 & . & . & . & . \\
 . & . & . & . & . & . & . & . & . & . & . & . & 1 & . & . & . \\
\end{array}
\right)
\]
}
\end{minipage}
\hfill
\begin{minipage}[t]{0.32\textwidth}
\vspace{0pt}
\centering
\includegraphics[width=\textwidth]{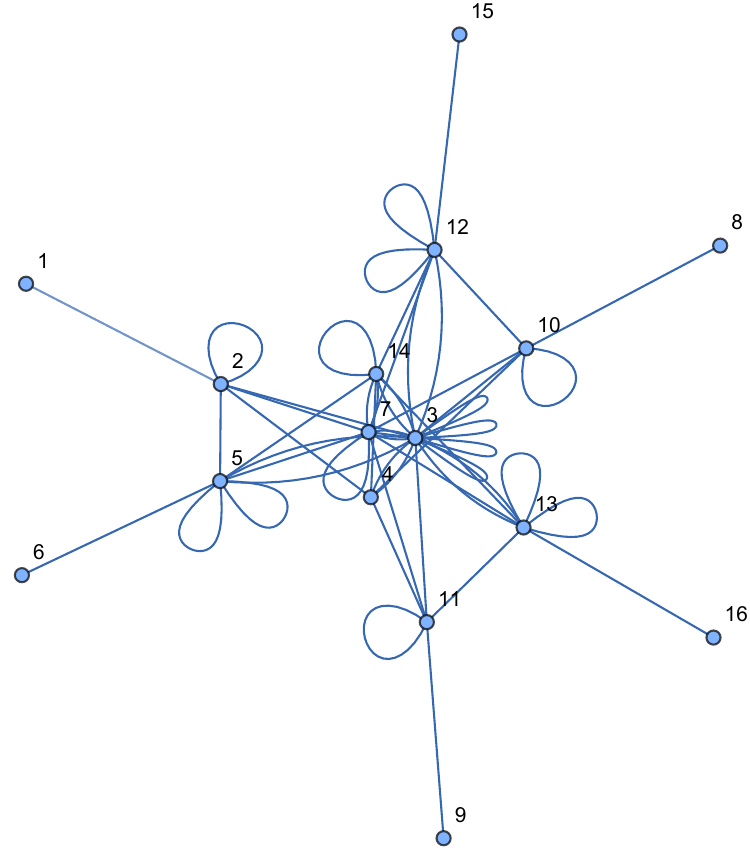}
\captionof{figure}{Adjacency graph of the matrix of the generator $0001$ of $(F_4,6)$ acting on the quantum subgroup $M$.}
\label{gp0001}
\end{minipage}
\end{center}
\begin{center}
\begin{minipage}[t]{0.62\textwidth}
\vspace{0pt}
{\footnotesize
\[
\setlength{\arraycolsep}{2.5pt}
 F_{0010}=
 \left(
\begin{array}{cccccccccccccccc}
 . & . & 1 & . & . & . & . & . & . & . & . & . & . & . & . & . \\
 . & 1 & 4 & 2 & 2 & . & 3 & . & . & 1 & 1 & 2 & 2 & 2 & . & . \\
 1 & 4 & 18 & 6 & 7 & 1 & 12 & 1 & 1 & 4 & 4 & 7 & 7 & 9 & 1 & 1 \\
 . & 2 & 6 & . & 2 & . & 3 & . & . & 2 & 2 & 2 & 2 & 3 & . & . \\
 . & 2 & 7 & 2 & 3 & . & 4 & . & . & 2 & 2 & 3 & 3 & 4 & . & . \\
 . & . & 1 & . & . & . & 1 & . & . & . & . & . & . & . & . & . \\
 . & 3 & 12 & 3 & 4 & 1 & 6 & . & . & 3 & 3 & 4 & 4 & 6 & 1 & 1 \\
 . & . & 1 & . & . & . & . & . & . & . & . & . & . & . & . & . \\
 . & . & 1 & . & . & . & . & . & . & . & . & . & . & . & . & . \\
 . & 1 & 4 & 2 & 2 & . & 3 & . & . & 1 & 1 & 2 & 2 & 2 & . & . \\
 . & 1 & 4 & 2 & 2 & . & 3 & . & . & 1 & 1 & 2 & 2 & 2 & . & . \\
 . & 2 & 7 & 2 & 3 & . & 4 & . & . & 2 & 2 & 3 & 3 & 4 & . & . \\
 . & 2 & 7 & 2 & 3 & . & 4 & . & . & 2 & 2 & 3 & 3 & 4 & . & . \\
 . & 2 & 9 & 3 & 4 & . & 6 & . & . & 2 & 2 & 4 & 4 & 3 & . & . \\
 . & . & 1 & . & . & . & 1 & . & . & . & . & . & . & . & . & . \\
 . & . & 1 & . & . & . & 1 & . & . & . & . & . & . & . & . & . \\
\end{array}
\right)
\]
}
\end{minipage}
\hfill
\begin{minipage}[t]{0.32\textwidth}
\vspace{0pt}
\centering
\includegraphics[width=\textwidth]{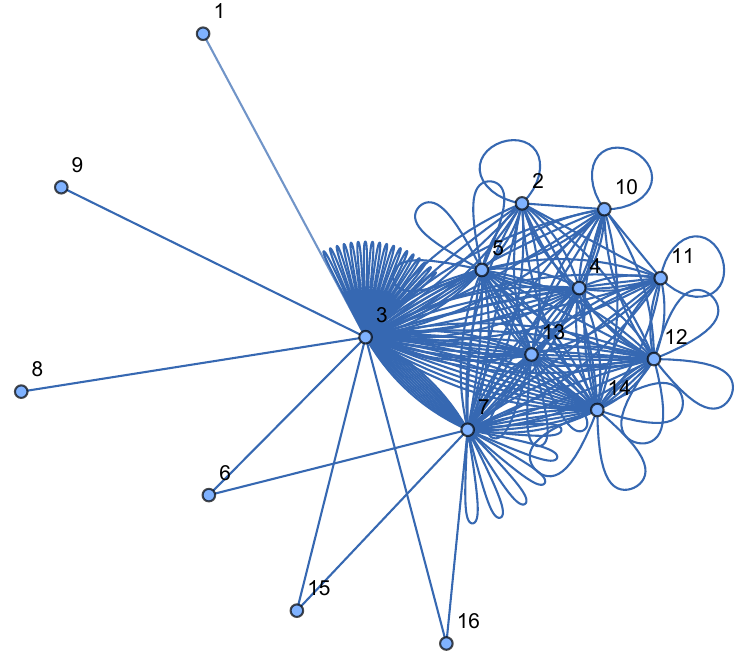}
\captionof{figure}{Adjacency graph of the matrix of the generator $0010$ of $(F_4,6)$ acting on the quantum subgroup $M$.}
\label{gp0010}
\end{minipage}
\end{center}
\begin{center}
\begin{minipage}[t]{0.62\textwidth}
\vspace{0pt}
{\footnotesize
\[
\setlength{\arraycolsep}{2.5pt}
 F_{0100}=
 \left(
\begin{array}{cccccccccccccccc}
 . & . & . & 1 & . & . & 1 & . & . & . & . & 1 & 1 & . & . & . \\
 . & 2 & 9 & . & 1 & . & 3 & . & . & 3 & 3 & 3 & 3 & 5 & 1 & 1 \\
 . & 9 & 32 & 7 & 12 & 1 & 17 & . & . & 9 & 9 & 12 & 12 & 17 & 1 & 1 \\
 1 & . & 7 & 8 & 5 & . & 10 & 1 & 1 & . & . & 5 & 5 & 2 & . & . \\
 . & 1 & 12 & 5 & 3 & . & 9 & 1 & 1 & 3 & 3 & 6 & 6 & 5 & 1 & 1 \\
 . & . & 1 & . & . & . & . & . & . & 1 & 1 & 1 & 1 & 1 & . & . \\
 1 & 3 & 17 & 10 & 9 & . & 15 & 1 & 1 & 3 & 3 & 9 & 9 & 7 & . & . \\
 . & . & . & 1 & 1 & . & 1 & . & . & . & . & . & 1 & . & . & . \\
 . & . & . & 1 & 1 & . & 1 & . & . & . & . & 1 & . & . & . & . \\
 . & 3 & 9 & . & 3 & 1 & 3 & . & . & 2 & 3 & 1 & 3 & 5 & . & 1 \\
 . & 3 & 9 & . & 3 & 1 & 3 & . & . & 3 & 2 & 3 & 1 & 5 & 1 & . \\
 1 & 3 & 12 & 5 & 6 & 1 & 9 & . & 1 & 1 & 3 & 3 & 6 & 5 & . & 1 \\
 1 & 3 & 12 & 5 & 6 & 1 & 9 & 1 & . & 3 & 1 & 6 & 3 & 5 & 1 & . \\
 . & 5 & 17 & 2 & 5 & 1 & 7 & . & . & 5 & 5 & 5 & 5 & 10 & 1 & 1 \\
 . & 1 & 1 & . & 1 & . & . & . & . & . & 1 & . & 1 & 1 & . & . \\
 . & 1 & 1 & . & 1 & . & . & . & . & 1 & . & 1 & . & 1 & . & . \\
\end{array}
\right)
\]
}
\end{minipage}
\hfill
\begin{minipage}[t]{0.32\textwidth}
\vspace{0pt}
\centering
\includegraphics[width=\textwidth]{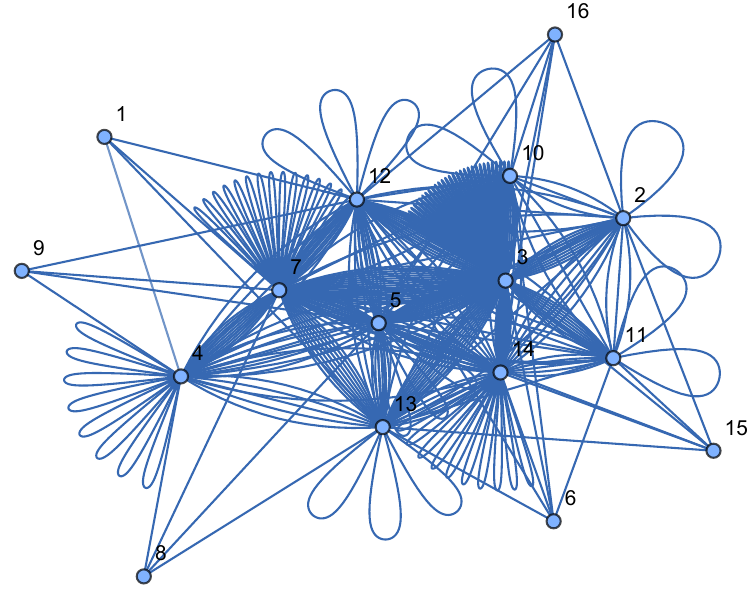}
\captionof{figure}{Adjacency graph of the matrix of the generator $0100$ of $(F_4,6)$ acting on the quantum subgroup $M$.}
\label{gp0100}
\end{minipage}
\end{center}

\begin{center}
\begin{minipage}[t]{0.62\textwidth}
\vspace{0pt}
{\footnotesize
\[
\setlength{\arraycolsep}{2.5pt}
 F_{1000}=
 \left(
\begin{array}{cccccccccccccccc}
 . & . & . & 1 & . & . & . & . & . & . & . & . & . & . & . & . \\
 . & 1 & 2 & . & . & . & . & . & . & 1 & 1 & . & . & 1 & . & . \\
 . & 2 & 6 & . & 2 & . & 3 & . & . & 2 & 2 & 2 & 2 & 3 & . & . \\
 1 & . & . & 3 & 1 & . & 3 & 1 & 1 & . & . & 1 & 1 & . & . & . \\
 . & . & 2 & 1 & 1 & . & 2 & . & . & . & . & 1 & 1 & 1 & . & . \\
 . & . & . & . & . & . & . & . & . & . & . & . & . & 1 & . & . \\
 . & . & 3 & 3 & 2 & . & 3 & . & . & . & . & 2 & 2 & . & . & . \\
 . & . & . & 1 & . & . & . & . & . & . & . & . & . & . & . & . \\
 . & . & . & 1 & . & . & . & . & . & . & . & . & . & . & . & . \\
 . & 1 & 2 & . & . & . & . & . & . & 1 & 1 & . & . & 1 & . & . \\
 . & 1 & 2 & . & . & . & . & . & . & 1 & 1 & . & . & 1 & . & . \\
 . & . & 2 & 1 & 1 & . & 2 & . & . & . & . & 1 & 1 & 1 & . & . \\
 . & . & 2 & 1 & 1 & . & 2 & . & . & . & . & 1 & 1 & 1 & . & . \\
 . & 1 & 3 & . & 1 & 1 & . & . & . & 1 & 1 & 1 & 1 & 3 & 1 & 1 \\
 . & . & . & . & . & . & . & . & . & . & . & . & . & 1 & . & . \\
 . & . & . & . & . & . & . & . & . & . & . & . & . & 1 & . & . \\
\end{array}
\right)
\]
}
\end{minipage}
\hfill
\begin{minipage}[t]{0.32\textwidth}
\vspace{0pt}
\centering
\includegraphics[width=\textwidth]{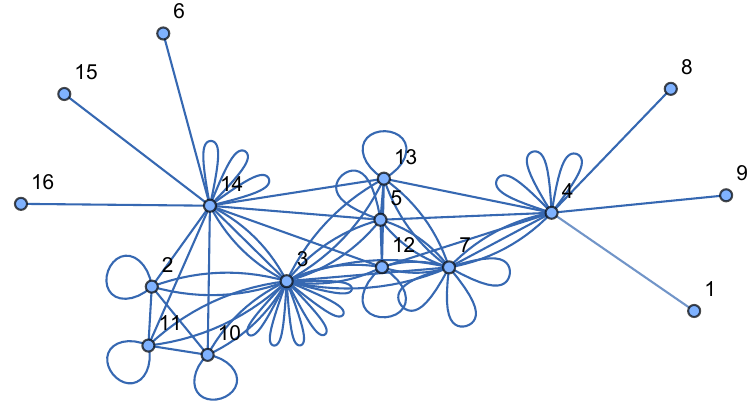}
\captionof{figure}{Adjacency graph of the matrix of the generator $1000$ of $(F_4,6)$ acting on the quantum subgroup $M$.}
\label{gp1000}
\end{minipage}
\end{center}

The above adjacency graphs contain multiple edges and self-loops (non-trivial multiplicities) that make the figures difficult to read.
Only the graph associated with the basic representation (the fundamental $0001$), and to some extent the graph of $1000$, are reasonably readable.

The following non-trivial checks can be performed: 
   \begin{itemize} 
  \item  The four matrices $F_f$ are symmetric.
 \item    The four matrices $F_f$ commute.
  \item   Their coefficients are non-negative integers. 
\item   The matrix elements of the $39$ matrices $F_n$ obtained from the fundamental ones $F_f$ by using the fusion $\ZZ$-polynomials of section \ref{FusionZpolynomials}  are non-negative integers.
{\scriptsize Example. $N_{1011}$ is obtained from the polynomial $P(1,0,1,1)=(t-x)(t^2+(x-z-1)t+y-z-1)$ by setting $x=N_{1000}$, $y=N_{0100}$, $z=N_{0010}$, $t=N_{0001}$. Its representative, the annular matrix $F_{1011}$, is obtained from the same polynomial by setting $x=F_{1000}$, $y=F_{0100}$, $z=F_{0010}$, $t=F_{0001}$.}
 \item  The matrix obtained from the three matrices $F_{0001}$, $F_{0010}$, $F_{1000}$,  by using the fusion $\QQ$-polynomial of section \ref{FusionQpolynomialsThreeVars}  indeed coincides with the $F_{0100}$ matrix.
  \item  All the matrices $F_f$ can be obtained from the $\QQ$-polynomials in $u$ (of degree 38) described in section \ref{FusionQpolynomials} by setting $u = 2 z + t$, i.e. replacing $u$ by $2 F_{0010} + F_{0001}$.
    \item  All the matrices $F_f$ can be obtained from the $\QQ$-polynomials in $t$ (of degree 11) given in the next section.
  \end{itemize}

    \subsection{The $\QQ$-polynomials  of degree $11$ in $t$} 
    \label{polynomialsof0001}
    In the quantum subgroup $M$, the element $t=F_{0001}$ is primitive (we recall that $N_{0001}$ is not primitive in the fusion ring). The generating polynomials have degree $11$ in $M$ --- this comes from the fact that the minimal polynomial $min(t)$ of $F_{0001}$ has degree $12$:
      {\scriptsize
    $$min(t) = (1 + 12 t + 14 t^2 - 12 t^3 + t^4) (1 + 3 t - t^2 - 3 t^3 + 
     t^4) (1 - 3 t - t^2 + 3 t^3 + t^4)$$}
The three fundamental matrices $F_{1000}$, $F_{0100}$, $F_{0010}$ can be obtained as: 
    {\scriptsize
    \begin{align*}   
r_x =& (1/660) (1 + 3 t - t^2 - 3 t^3 + t^4) (2016 - 5001 t - 14817 t^2 + 
     13308 t^3 + 14466 t^4 - 5170 t^5 - 2049 t^6 + 227 t^7)\\
r_y =& (1/
   660) (443 - 2469 t - 2190 t^2 + 20386 t^3 - 2940 t^4 - 33996 t^5 + 
     13520 t^6 + 16266 t^7 - 8610 t^8 - 896 t^9 + 813 t^{10} - 63 t^{11})\\
r_z =& -(1/660) (1 + 3 t - t^2 - 3 t^3 + t^4) (1668 - 2655 t - 4235 t^2 + 
     5004 t^3 + 3508 t^4 - 1540 t^5 - 513 t^6 + 59 t^7)
         \end{align*}}

\subsection{Essential matrices for $M$}

\noindent
\begin{minipage}[t]{0.43\textwidth}
\vspace{0pt}
\centering
{\footnotesize
\setlength{\arraycolsep}{1pt}
\renewcommand{\arraystretch}{0.7}
\[
Ess_1=\quad
\left(
\begin{array}{cccccccccccccccc}
 1 & . & . & . & . & . & . & . & . & . & . & . & . & . & . & . \\
 . & 1 & . & . & . & . & . & . & . & . & . & . & . & . & . & . \\
 . & . & . & . & 1 & . & 1 & . & . & . & . & . & . & . & . & . \\
 . & . & 1 & . & . & . & . & . & . & . & . & . & . & . & . & . \\
 . & . & . & 1 & . & . & . & . & . & . & . & . & . & . & . & . \\
 . & 1 & 1 & . & 1 & 1 & . & . & . & . & . & . & . & 1 & . & . \\
 . & . & 2 & . & 1 & . & 1 & . & . & . & . & 1 & 1 & 1 & . & . \\
 . & . & . & 1 & . & . & 1 & . & . & . & . & 1 & 1 & . & . & . \\
 . & . & 1 & . & . & . & . & . & . & 1 & 1 & . & . & 1 & . & . \\
 1 & 1 & 1 & 1 & 1 & . & 1 & . & . & . & . & . & . & . & . & . \\
 . & 1 & 3 & 1 & 2 & . & 2 & . & . & 1 & 1 & 1 & 1 & 2 & . & . \\
 . & . & 2 & 1 & 1 & . & 2 & . & . & . & . & 1 & 1 & 1 & . & . \\
 . & 1 & 3 & . & . & . & 1 & . & . & 1 & 1 & 1 & 1 & 2 & 1 & 1 \\
 . & . & 2 & 2 & 1 & . & 2 & . & . & . & . & 1 & 1 & . & . & . \\
 . & 1 & 2 & . & . & . & 1 & . & . & 1 & 1 & 1 & 1 & 1 & . & . \\
 . & . & . & 1 & . & . & 1 & 1 & 1 & . & . & . & . & . & . & . \\
 . & 1 & 1 & . & 1 & . & . & . & . & . & . & . & . & . & . & . \\
 . & 1 & 3 & 1 & 1 & . & 2 & . & . & 1 & 1 & 1 & 1 & 1 & . & . \\
 . & 1 & 4 & 1 & 2 & 1 & 2 & . & . & 1 & 1 & 1 & 1 & 2 & . & . \\
 . & . & 3 & 2 & 1 & . & 3 & 1 & 1 & 1 & 1 & 2 & 2 & 1 & . & . \\
 . & 1 & 3 & 1 & 1 & . & 1 & . & . & 1 & 1 & 1 & 1 & 2 & . & . \\
 . & 1 & 3 & . & . & . & 1 & . & . & 1 & 1 & 1 & 1 & 2 & . & . \\
 . & 1 & 5 & 1 & 2 & . & 3 & . & . & 1 & 1 & 2 & 2 & 2 & . & . \\
 1 & . & 1 & 2 & 1 & . & 2 & . & . & . & . & 1 & 1 & . & . & . \\
 . & 1 & 2 & . & 1 & . & . & . & . & 1 & 1 & . & . & 1 & . & . \\
 . & . & . & . & 1 & 1 & . & . & . & . & . & . & . & . & . & . \\
 . & . & 1 & . & 1 & . & 1 & . & . & . & . & . & . & 1 & . & . \\
 . & 1 & 2 & 1 & 1 & . & 1 & . & . & . & . & 1 & 1 & 1 & . & . \\
 1 & 1 & 1 & . & 1 & . & 1 & . & . & . & . & . & . & 1 & . & . \\
 . & . & 2 & . & . & . & 1 & . & . & 1 & 1 & 1 & 1 & 1 & . & . \\
 . & 1 & 3 & 1 & 1 & . & 2 & . & . & 1 & 1 & 1 & 1 & 1 & . & . \\
 . & . & 1 & 1 & . & . & 1 & . & . & . & . & . & . & . & . & . \\
 . & . & 1 & 1 & . & . & 1 & . & . & . & . & 1 & 1 & 1 & . & . \\
 . & . & 4 & 1 & 1 & . & 2 & . & . & 1 & 1 & 2 & 2 & 2 & 1 & 1 \\
 . & 1 & 2 & 1 & 2 & . & 2 & . & . & . & . & 1 & 1 & 1 & . & . \\
 . & 1 & 3 & . & 1 & . & 1 & . & . & 1 & 1 & 1 & 1 & 2 & . & . \\
 . & . & 1 & 1 & 1 & . & 2 & . & . & . & . & 1 & 1 & 1 & . & . \\
 . & 1 & 2 & . & 1 & 1 & 1 & . & . & . & . & . & . & 1 & . & . \\
 . & . & . & 1 & 1 & . & . & . & . & . & . & . & . & . & . & . \\
\end{array}
\right)
\]
}
\end{minipage}
\hfill
\begin{minipage}[t]{0.53\textwidth}
\vspace{0pt}
\small

The module structure of $M$ over the fusion algebra is described by the
square matrices $F_n$, with matrix elements $(F_n)_{a,b}$, where
$n=1,\ldots,39$ and $a,b=1,\ldots,16$.

Setting
\[
(Ess_a)_{n,b}=(F_n)_{a,b},
\]
it can also be described by the rectangular matrices $Ess_a$.
Among these rectangular matrices, sometimes called ``essential matrices'',
the first one plays a special role. It is sometimes called the
``main intertwiner'' because
\[
N_n Ess_1 = Ess_1 F_n,
\]
for all $n=1,\ldots,39$.

Since we know all the $F_n$'s, we can determine all the $Ess_a$'s.
In particular, the first one is displayed here.

As we shall see below, it determines the induction--restriction rules.

\end{minipage}

  \subsection{Induction--restriction rules: The associated quantum homogeneous space and the quantum vector bundle paradigm}
  
Classically, the dimension of the space ${\Gamma_a}$ of vector-valued functions defined on the quotient of a finite group $A$ by a subgroup $M$, and valued in a vector space $a$, is equal to $\vert {\Gamma_a} \vert =dim(a)\times \vert A/M \vert$.
If $a$ is a representation space for $M$,  the dimension $\vert {\Gamma_a} \vert$ can also be obtained by decomposing this space of functions  into a sum of irreducible representations of $A$. In particular, if $a$ is the trivial representation of $M$, one obtains in this way  a decomposition of the algebra ${\mathcal F}$ of complex-valued functions on $A/M$ and we have $\vert {\mathcal F} \vert=\vert {\Gamma_1} \vert = \vert A/M \vert$.
 From a geometrical point of view, $A$ is a principal $M$-bundle over $A/M$, and for each irrep $a$ of $M$ one builds an associated vector bundle over $A/M$ whose space of sections ${\Gamma_a}$ can be decomposed, using induction, as a sum $\Gamma_a = \sum _{n  \uparrow \Gamma_a} \,  n$, where $n$ runs over irreps of $A$.

Returning to the present situation, we define  $\Gamma_a = \sum _{n  \uparrow \Gamma_a} \,  n$, where the sum runs over those simple objects $n$ of  ${\mathcal A}$ that appear in the column $a$ of the essential matrix $Ess_{{1}}$.
We think of $\Gamma_a$ as the space of sections of a ``quantum vector bundle''  over a would-be homogeneous space  $\mathcal{A}/\mathcal{M}$ which is both discrete and quantum, and we decompose $\Gamma_a$  into a sum of simple objects of $\mathcal{A}$.
Taking $a={1}$, one obtains  ${\mathcal F}=\Gamma_{1} =\sum _{n  \uparrow \Gamma_{1}} n$, and ${\mathcal F}$ (the étale algebra) plays the role of an algebra of functions on a would-be homogeneous space.

 We have already calculated the quantum dimensions $\text{qdim}(n)$ of the simple objects of the fusion algebra $\mathcal A$ in  \ref{quantumdimensionsofA}, and the global dimension $\vert {\mathcal A}\vert$ of $\mathcal A$ in  \ref{globaldimensionofA}. 
 The (quantum) dimension of $\Gamma_a$ can be calculated as the sum of quantum dimensions $\sum _{n  \uparrow \Gamma_a} \,  \text{qdim}(n)$ It can also be expressed in terms of the quantum dimensions\footnote{Warning: For quantum dimensions, the notation $\text{qdim}(n)$ refers to simple objects of $\mathcal A$, whereas $q_a$ refers to simple objects of $\mathcal M$.} $q_a$ of the simple objects $a$ of ${\mathcal M}$ as $q_a \times \vert {\mathcal A}/{\mathcal M}  \vert = q_a \times \vert {\mathcal A}\vert / \vert {\mathcal M}  \vert $.
 
 In particular, taking $a=1$, one obtains $$\vert {\mathcal F} \vert = \sum _{n  \uparrow \Gamma_1} \,  \text{qdim}(n)$$ and $\vert {\mathcal M}  \vert = 
 \vert {\mathcal A}\vert / \vert {\mathcal F}  \vert$.
 
  \subsection{\'Etale algebra and dyslectic modules}
 \label{dyslecticdecomp}

  The ambichiral vertices of $O$, labelled $\{1, 18, 32, 33, 60, 61\}$, or the modular vertices of $M$, labelled $\{1,6,8,9,15,16\}$, determine, by induction,  modules whose decomposition into simple objects of $\mathcal A$ is given by the columns of $Ess_1$ with the same labels. In particular, we obtain the decomposition of the commutative algebra ${\mathcal F}$. Notice that the members of the two pairs of modular vertices $(8,9)$ and $(15,16)$ correspond to equivalent $A$-objects.
\begin{align*}
  \Gamma_{1} = {\mathcal F} &= \{0, 0, 0, 0\}    \oplus \{0, 0, 0, 4\}   \oplus \{1, 1, 0, 0\}   \oplus \{0, 0, 3, 0 \} \\
   \Gamma_{6} &=  \{0, 0, 0, 3\}    \oplus \{0, 0, 2, 1\}   \oplus \{0, 0, 0, 6\}   \oplus \{2, 0, 1, 0 \} \\
  \Gamma_{8} \simeq  \Gamma_{9} &=  \{2, 0, 0, 0\}    \oplus \{0, 1, 0, 2\}    \\
    \Gamma_{15} \simeq  \Gamma_{16} &= \{0, 1, 0, 1\}    \oplus \{1, 0, 1, 2\}  \\
  \end{align*}

 \subsection{Quantum dimensions and other numerical results}
 
\subsubsection{Values of $\vert{\mathcal F}\vert$ and $\vert {\Gamma_a} \vert$}

From $\vert \Gamma_a \vert = \sum _{n  \uparrow \Gamma_a} \,  \text{qdim}(n)$, we obtain: 

{\footnotesize
\[
\begin{array}{@{}l@{\quad}l@{}}
\vert{\mathcal F}\vert =
\vert {\Gamma_{1}} \vert=105+45 \sqrt{5}+4 \sqrt{15 \left(85+38 \sqrt{5}\right)}, &
\vert {\Gamma_{2}} \vert=1080+480 \sqrt{5}+\sqrt{30 \left(77305+34571 \sqrt{5}\right)},\\[3pt]
\vert {\Gamma_{3}} \vert=3 \left(1465+655 \sqrt{5}+\sqrt{30 \left(143285+64079  \sqrt{5}\right)}\right),&
\vert {\Gamma_{4}} \vert=   3 \left(450+200 \sqrt{5}+\sqrt{30 \left(13325+5959 \sqrt{5}\right)}\right),\\[3pt]
\vert {\Gamma_{5}} \vert=   1740+780 \sqrt{5}+\sqrt{30 \left(202385+90509 \sqrt{5}\right)},&
\vert {\Gamma_{6}} \vert=   165+75 \sqrt{5}+2   \sqrt{30 \left(445+199 \sqrt{5}\right)},\\[3pt]
\vert {\Gamma_{7}} \vert=   3 \left(905+405 \sqrt{5}+2 \sqrt{15 \left(27365+12238 \sqrt{5}\right)}\right),&
 \vert {\Gamma_{8}} \vert=  105+45 \sqrt{5}+4 \sqrt{15 \left(85+38 \sqrt{5}\right)},\\[3pt]
 \vert {\Gamma_{9}} \vert=  105+45  \sqrt{5}+4 \sqrt{15 \left(85+38 \sqrt{5}\right)},&
  \vert {\Gamma_{10}} \vert= 1080+480 \sqrt{5}+\sqrt{30 \left(77305+34571 \sqrt{5}\right)},\\[3pt]
\vert {\Gamma_{11}} \vert=   1080+480 \sqrt{5}+\sqrt{30 \left(77305+34571 \sqrt{5}\right)},&
\vert {\Gamma_{12}} \vert=   1740+780  \sqrt{5}+\sqrt{30 \left(202385+90509 \sqrt{5}\right)},\\[3pt]
\vert {\Gamma_{13}} \vert=   1740+780 \sqrt{5}+\sqrt{30 \left(202385+90509 \sqrt{5}\right)},&
 \vert {\Gamma_{14}} \vert=  3 \left(725+325 \sqrt{5}+\sqrt{30 \left(34885+15601   \sqrt{5}\right)}\right),\\[3pt]
 \vert {\Gamma_{15}} \vert=  165+75 \sqrt{5}+2 \sqrt{30 \left(445+199 \sqrt{5}\right)},&
\vert {\Gamma_{16}} \vert=   165+75 \sqrt{5}+2 \sqrt{30 \left(445+199 \sqrt{5}\right)}.
   \end{array}
   \]
   }

  \subsubsection{Value of $\vert{\mathcal M}\vert$ (from $\vert{\mathcal F}\vert$ and $\vert{\mathcal A}\vert$ )}
  
  From $\vert {\mathcal F} \vert = \vert {\Gamma_{1}} \vert$, from $\vert {\mathcal A} \vert$   in  \ref{globaldimensionofA},
   and $\vert {\mathcal M}  \vert =  \vert {\mathcal A}\vert / \vert {\mathcal F}  \vert$, we obtain: 
   $$\vert {\mathcal M}  \vert = 15 \left(75+33 \sqrt{5}+2 \sqrt{6 \left(445+199 \sqrt{5}\right)}\right)$$

Remark: In our case, the decomposition of ${\mathcal F}$ into simple objects of ${\mathcal A}$  (the first column of the matrix $Ess_1$) can also be read from the first modular block of the partition function, $V = \chi (0,0,0,0)+\chi (0,0,0,4)+\chi (0,0,3,0)+\chi (1,1,0,0)$, as in sections \ref{modularZ}  or  \ref{bilinears}.  

We also have 
$$ \vert {\mathcal A}/{\mathcal M}  \vert = 105+45 \sqrt{5}+4 \sqrt{15 \left(85+38 \sqrt{5}\right)}$$
 
\subsubsection{q-dimensions $q_a$ of simple objects $a$ of ${\mathcal M}$ (from $\vert {\Gamma_a} \vert$)}

From $q_a = \vert \Gamma_a \vert  /    \vert \Gamma_1 \vert =  \vert \Gamma_a \vert /   \vert {\mathcal A}/{\mathcal M}  \vert$, we get 
{\footnotesize
\[
\begin{array}{@{}l@{\quad}l@{}}
q_{1}= 1, &
q_{2}= 3+\sqrt{5}+\sqrt{3 \left(5+2 \sqrt{5}\right)} ,\\[3pt]
q_{3}= \frac{3}{2} \left(7+3 \sqrt{5}\right)+\sqrt{3 \left(85+38 \sqrt{5}\right)},&
q_{4}= \frac{3}{2} \left(3+\sqrt{5}\right)+\sqrt{3 \left(5+2 \sqrt{5}\right)} ,\\[3pt]
 q_{5}=4+2 \sqrt{5}+\sqrt{\frac{3}{2} \left(25+11 \sqrt{5}\right)} ,&
q_{6}= \frac{1}{2} \left(1+\sqrt{5}\right) ,\\[3pt]
q_{7}= 6+3 \sqrt{5}+\sqrt{\frac{3}{2} \left(65+29 \sqrt{5}\right)} ,&
 q_{8}=1 ,\\[3pt]
q_{9}= 1 ,&
q_{10}= 3+\sqrt{5}+\sqrt{3 \left(5+2 \sqrt{5}\right)} ,\\[3pt]
q_{11}= 3+\sqrt{5}+\sqrt{3 \left(5+2 \sqrt{5}\right)} ,&
q_{12}= 4+2 \sqrt{5}+\sqrt{\frac{3}{2} \left(25+11 \sqrt{5}\right)} ,\\[3pt]
q_{13}= 4+2 \sqrt{5}+\sqrt{\frac{3}{2} \left(25+11 \sqrt{5}\right)} ,&
q_{14}= 6+3 \sqrt{5}+\sqrt{\frac{3}{2} \left(25+11 \sqrt{5}\right)} ,\\[3pt]
q_{15}= \frac{1}{2} \left(1+\sqrt{5}\right) ,&
q_{16}= \frac{1}{2} \left(1+\sqrt{5}\right).
\end{array}
   \]
   }

\subsubsection{q-dimensions $q_a$ (from Perron-Frobenius)}

Let us consider the eigensystem associated, for instance, with the annular matrix $F_{0001}$.
Its largest eigenvalue is $3+\sqrt{5}+\sqrt{3 \left(5+2 \sqrt{5}\right)}$. One may check that the associated normalized eigenvector (we normalize its first component to $1$) has the same components $q_a$ as those given previously. 
One may follow the same method, starting, for instance, from the $1000$ fundamental irrep, or from any other, instead of $0001$, to calculate, or recover, the values of $q_a$. This provides a useful check of the previous calculations. 

\subsubsection{Recovering $\vert{\mathcal M}\vert$ from $q_a$ (summing squares)}
Now that we have determined the values of $q_a$, and checked these results using various methods, we may also recover the obtained value for the global dimension of $\mathcal M$ by summing squares. Indeed one should have the equality:
$$\vert{\mathcal M}  \vert= \sum_a (q_a)^2$$ 
with a summation over the simple objects $a$ of ${\mathcal M}$. 
This equality provides another check of the calculations.

\subsubsection{Recovering $\vert{\mathcal M}\vert$ from ambichiral data}

For a modular invariant and a quantum subgroup associated with a conformal embedding of an affine Lie algebra at level $k$ into one at level $1$, calling $\vert{\mathcal J}\vert$  the global dimension of the fusion algebra of the latter, and using the same notation as before, one would have $\vert{\mathcal A}\vert = \vert{\mathcal O}\vert =\vert {\mathcal M}\vert \times \vert{\mathcal M}\vert /\vert{\mathcal J}\vert$, so that $\vert{\mathcal M}\vert$ could be obtained as $\sqrt{ \vert{\mathcal A}\vert \times \vert{\mathcal J}\vert  }$. 
The number $\vert{\mathcal J}\vert$
 can also be obtained as the sum  $\sum_{\hat a} (q_{\hat a})^2$ with ${\hat a}$ belonging to the set of modular generators of ${\mathcal M}$, or, equivalently, to the set of simple objects of ${\mathcal O}$ that generate its ambichiral subalgebra.
 
In the situation studied in this paper  we do not have a conformal embedding into an affine Lie algebra but a conformal extension within a vertex operator algebra.
However, we observe that the above property still holds (we did not try to prove it on general grounds): 
The modular generators are those associated with vertices ${\hat a} \in \{1, 6, 8, 9, 15, 16\}$ of the adjacency graphs of the matrices $F_n$; so,
using the values of $q_a$ obtained previously, one finds $\vert{\mathcal J}\vert = \sum_{\hat a} (q_{\hat a})^2 = \frac{3}{2} \left(5+\sqrt{5}\right)$, and, using the known value of  $\vert{\mathcal A}\vert$ one can indeed check that $\vert{\mathcal M}\vert = \sqrt{ \vert{\mathcal A}\vert \times \vert{\mathcal J}\vert}$.

Remark. The previous equations relating $\vert{\mathcal A}\vert, \vert{\mathcal F}\vert, \vert{\mathcal M}\vert, \vert{\mathcal J}\vert$, imply $\vert{\mathcal A} \vert = \vert {\mathcal F} \vert^2 \, \vert {\mathcal J} \vert$, but the set ${\mathcal J}$ contains the unit (inherited from ${\mathcal O}$), so that $\vert {\mathcal J} \vert \geq 1$. Therefore $\vert {\mathcal A} \vert \geq \vert {\mathcal F} \vert^2$. 
Actually, the same inequality holds not only for $\mathcal F=\vert {\Gamma_{1}} \vert$ but for all the $\vert {\Gamma_{\hat a}} \vert$ when $\hat a$ is ambichiral.

  \subsection{The quantum subgroup $M$: a discussion of its self-fusion}
 \label{selffusionofM}
     In order to discuss the structure of $M$ as a fusion-based ring,   compatible with the action of the fusion ring $A$, we work with the same basis $a,b,\ldots =1,2,\ldots,16$ used to label the vertices of the graphs previously used to display the module structure of $M$ over $A$.  We shall call $g_a$ the self-fusion matrices (``graph matrices'') that encode the multiplication by these basis elements. The unit is $g_1$ and it is trivially identified with $F_1$. The annular matrices $F_n$ can be decomposed in terms of the self-fusion matrices $g_a$.
    Formally, these decompositions express, in a particular case, the compatibility  between the module action (matrices $F_n$), and the self-fusion (matrices $g_a$).
For all $39$ irreps $n$, one should have  
\begin{equation}
F_n = \sum_b (F_n)_{1,b} \, g_b
\label{compatibilityspecial}
\end{equation}
 where $b$ runs over the $16$ basis elements of $M$. In other words, 
the decomposition of $F_n$ in terms of the $g_a$'s, stemming from the fact that  $g_1 \times g_a = g_a$,  is given by the first row of the matrix $F_n$. 

The following array displays the first row of each $F_n$, where the labels $n$
of the $F_n$ are sorted and defined as in section \ref{simpleobjects}.

\medskip
\noindent
\begin{minipage}[t]{0.40\textwidth}
\vspace{0pt}
\centering
{\scriptsize
\[
\setlength{\arraycolsep}{2.2pt}
\renewcommand{\arraystretch}{0.78}
\left(
\begin{array}{cccccccccccccccc}
 1 & . & . & . & . & . & . & . & . & . & . & . & . & . & . & . \\
 . & 1 & . & . & . & . & . & . & . & . & . & . & . & . & . & . \\
 . & . & . & . & 1 & . & 1 & . & . & . & . & . & . & . & . & . \\
 . & . & 1 & . & . & . & . & . & . & . & . & . & . & . & . & . \\
 . & . & . & 1 & . & . & . & . & . & . & . & . & . & . & . & . \\
 . & 1 & 1 & . & 1 & 1 & . & . & . & . & . & . & . & 1 & . & . \\
 . & . & 2 & . & 1 & . & 1 & . & . & . & . & 1 & 1 & 1 & . & . \\
 . & . & . & 1 & . & . & 1 & . & . & . & . & 1 & 1 & . & . & . \\
 . & . & 1 & . & . & . & . & . & . & 1 & 1 & . & . & 1 & . & . \\
 1 & 1 & 1 & 1 & 1 & . & 1 & . & . & . & . & . & . & . & . & . \\
 . & 1 & 3 & 1 & 2 & . & 2 & . & . & 1 & 1 & 1 & 1 & 2 & . & . \\
 . & . & 2 & 1 & 1 & . & 2 & . & . & . & . & 1 & 1 & 1 & . & . \\
 . & 1 & 3 & . & . & . & 1 & . & . & 1 & 1 & 1 & 1 & 2 & 1 & 1 \\
 . & . & 2 & 2 & 1 & . & 2 & . & . & . & . & 1 & 1 & . & . & . \\
 . & 1 & 2 & . & . & . & 1 & . & . & 1 & 1 & 1 & 1 & 1 & . & . \\
 . & . & . & 1 & . & . & 1 & 1 & 1 & . & . & . & . & . & . & . \\
 . & 1 & 1 & . & 1 & . & . & . & . & . & . & . & . & . & . & . \\
 . & 1 & 3 & 1 & 1 & . & 2 & . & . & 1 & 1 & 1 & 1 & 1 & . & . \\
 . & 1 & 4 & 1 & 2 & 1 & 2 & . & . & 1 & 1 & 1 & 1 & 2 & . & . \\
 . & . & 3 & 2 & 1 & . & 3 & 1 & 1 & 1 & 1 & 2 & 2 & 1 & . & . \\
 . & 1 & 3 & 1 & 1 & . & 1 & . & . & 1 & 1 & 1 & 1 & 2 & . & . \\
 . & 1 & 3 & . & . & . & 1 & . & . & 1 & 1 & 1 & 1 & 2 & . & . \\
 . & 1 & 5 & 1 & 2 & . & 3 & . & . & 1 & 1 & 2 & 2 & 2 & . & . \\
 1 & . & 1 & 2 & 1 & . & 2 & . & . & . & . & 1 & 1 & . & . & . \\
 . & 1 & 2 & . & 1 & . & . & . & . & 1 & 1 & . & . & 1 & . & . \\
 . & . & . & . & 1 & 1 & . & . & . & . & . & . & . & . & . & . \\
 . & . & 1 & . & 1 & . & 1 & . & . & . & . & . & . & 1 & . & . \\
 . & 1 & 2 & 1 & 1 & . & 1 & . & . & . & . & 1 & 1 & 1 & . & . \\
 1 & 1 & 1 & . & 1 & . & 1 & . & . & . & . & . & . & 1 & . & . \\
 . & . & 2 & . & . & . & 1 & . & . & 1 & 1 & 1 & 1 & 1 & . & . \\
 . & 1 & 3 & 1 & 1 & . & 2 & . & . & 1 & 1 & 1 & 1 & 1 & . & . \\
 . & . & 1 & 1 & . & . & 1 & . & . & . & . & . & . & . & . & . \\
 . & . & 1 & 1 & . & . & 1 & . & . & . & . & 1 & 1 & 1 & . & . \\
 . & . & 4 & 1 & 1 & . & 2 & . & . & 1 & 1 & 2 & 2 & 2 & 1 & 1 \\
 . & 1 & 2 & 1 & 2 & . & 2 & . & . & . & . & 1 & 1 & 1 & . & . \\
 . & 1 & 3 & . & 1 & . & 1 & . & . & 1 & 1 & 1 & 1 & 2 & . & . \\
 . & . & 1 & 1 & 1 & . & 2 & . & . & . & . & 1 & 1 & 1 & . & . \\
 . & 1 & 2 & . & 1 & 1 & 1 & . & . & . & . & . & . & 1 & . & . \\
 . & . & . & 1 & 1 & . & . & . & . & . & . & . & . & . & . & . \\.
\end{array}
\right)
\]
}
\end{minipage}
\hfill
\begin{minipage}[t]{0.56\textwidth}
\vspace{0pt}
\small

We read off
\[
F_1=g_1,\quad
F_2=g_2,\quad
F_3=g_5+g_7,\quad
F_4=g_3,\quad
F_5=g_4,
\]
and
\[
F_6=g_2+g_3+g_5+g_6,\quad \etc
\]
For instance,
\[
F_8=F_{0100}=g_4+g_7+(g_{12}+g_{13}).
\]

Conversely, one can use these equations to express many matrices $g_a$
in terms of the annular matrices $F_n$. One finds
\[
\begin{aligned}
g_1 &= F_1, &
g_2 &= F_2, &
g_3 &= F_4, &
g_4 &= F_5,\\
g_5 &= F_{39}-g_4,\\
g_6 &= F_{26}-g_5,\\
g_7 &= F_{32}-g_3-g_4,\\
g_{14} &= F_{27}-g_3-g_5-g_7 .
\end{aligned}
\]

For the remaining ones, only the sums are determined:
\[
\begin{aligned}
g_{8+9} &=
F_{16}-g_4-g_7,\\
g_{10+11} &=
g_2g_4-g_2-2g_3-g_{14},\\
g_{12+13} &=
F_{28}-g_2-2g_3-g_4-g_5-g_7-g_{14},\\
g_{15+16} &=
F_{13}-g_2-3g_3-g_7-g_{10+11}\\
&\qquad{}-g_{12+13}-2g_{14}.
\end{aligned}
\]

The above procedure is not sufficient to disentangle $8$ from $9$,
$10$ from $11$, $12$ from $13$, and $15$ from $16$.
\end{minipage}

 Before discussing the ``couplings" $(8,9)$, $(10,11)$, $(12,13)$, and $(15,16)$, let us display the graph matrices $g_j$ obtained so far.
The matrix $g_1=F_1=F_{0000}$ is the identity, and we already displayed  $g_2=F_2=F_{0001}$, $g_3=F_4=F_{0010}$, and $g_4=F_5=F_{1000}$ in section \ref{annularmatricesM}. Matrices $g_5, g_6, g_7, g_{14}$ are as follows:
{\footnotesize
\setlength{\arraycolsep}{2.5pt}
\renewcommand{\arraystretch}{0.8}
\begin{align*}
g_5=&
\left(
\begin{array}{cccccccccccccccc}
 . & . & . & . & 1 & . & . & . & . & . & . & . & . & . & . & . \\
 . & 1 & 2 & . & 2 & 1 & 1 & . & . & . & . & . & . & 1 & . & . \\
 . & 2 & 7 & 2 & 3 & . & 4 & . & . & 2 & 2 & 3 & 3 & 4 & . & . \\
 . & . & 2 & 1 & 1 & . & 2 & . & . & . & . & 1 & 1 & 1 & . & . \\
 1 & 2 & 3 & 1 & 3 & 1 & 2 & . & . & . & . & . & . & 1 & . & . \\
 . & 1 & . & . & 1 & . & . & . & . & . & . & . & . & . & . & . \\
 . & 1 & 4 & 2 & 2 & . & 3 & . & . & 1 & 1 & 2 & 2 & 2 & . & . \\
 . & . & . & . & . & . & . & . & . & . & . & 1 & . & . & . & . \\
 . & . & . & . & . & . & . & . & . & . & . & . & 1 & . & . & . \\
 . & . & 2 & . & . & . & 1 & . & . & 1 & . & 2 & . & 1 & 1 & . \\
 . & . & 2 & . & . & . & 1 & . & . & . & 1 & . & 2 & 1 & . & 1 \\
 . & . & 3 & 1 & . & . & 2 & 1 & . & 2 & . & 3 & . & 1 & 1 & . \\
 . & . & 3 & 1 & . & . & 2 & . & 1 & . & 2 & . & 3 & 1 & . & 1 \\
 . & 1 & 4 & 1 & 1 & . & 2 & . & . & 1 & 1 & 1 & 1 & 2 & . & . \\
 . & . & . & . & . & . & . & . & . & 1 & . & 1 & . & . & . & . \\
 . & . & . & . & . & . & . & . & . & . & 1 & . & 1 & . & . & . \\
\end{array}
\right)
, \quad
g_6=&
\left(
\begin{array}{cccccccccccccccc}
 . & . & . & . & . & 1 & . & . & . & . & . & . & . & . & . & . \\
 . & . & . & . & 1 & . & . & . & . & . & . & . & . & . & . & . \\
 . & . & 1 & . & . & . & 1 & . & . & . & . & . & . & . & . & . \\
 . & . & . & . & . & . & . & . & . & . & . & . & . & 1 & . & . \\
 . & 1 & . & . & 1 & . & . & . & . & . & . & . & . & . & . & . \\
 1 & . & . & . & . & 1 & . & . & . & . & . & . & . & . & . & . \\
 . & . & 1 & . & . & . & . & . & . & . & . & . & . & . & . & . \\
 . & . & . & . & . & . & . & . & . & . & . & . & . & . & 1 & . \\
 . & . & . & . & . & . & . & . & . & . & . & . & . & . & . & 1 \\
 . & . & . & . & . & . & . & . & . & . & . & 1 & . & . & . & . \\
 . & . & . & . & . & . & . & . & . & . & . & . & 1 & . & . & . \\
 . & . & . & . & . & . & . & . & . & 1 & . & 1 & . & . & . & . \\
 . & . & . & . & . & . & . & . & . & . & 1 & . & 1 & . & . & . \\
 . & . & . & 1 & . & . & . & . & . & . & . & . & . & 1 & . & . \\
 . & . & . & . & . & . & . & 1 & . & . & . & . & . & . & 1 & . \\
 . & . & . & . & . & . & . & . & 1 & . & . & . & . & . & . & 1 \\
\end{array}
\right)
, \\
g_7=&
\left(
\begin{array}{cccccccccccccccc}
 . & . & . & . & . & . & 1 & . & . & . & . & . & . & . & . & . \\
 . & 1 & 3 & . & 1 & . & 1 & . & . & 1 & 1 & 1 & 1 & 2 & . & . \\
 . & 3 & 12 & 3 & 4 & 1 & 6 & . & . & 3 & 3 & 4 & 4 & 6 & 1 & 1 \\
 . & . & 3 & 3 & 2 & . & 3 & . & . & . & . & 2 & 2 & . & . & . \\
 . & 1 & 4 & 2 & 2 & . & 3 & . & . & 1 & 1 & 2 & 2 & 2 & . & . \\
 . & . & 1 & . & . & . & . & . & . & . & . & . & . & . & . & . \\
 1 & 1 & 6 & 3 & 3 & . & 6 & 1 & 1 & 1 & 1 & 3 & 3 & 3 & . & . \\
 . & . & . & . & . & . & 1 & . & . & . & . & . & . & . & . & . \\
 . & . & . & . & . & . & 1 & . & . & . & . & . & . & . & . & . \\
 . & 1 & 3 & . & 1 & . & 1 & . & . & 1 & 1 & 1 & 1 & 2 & . & . \\
 . & 1 & 3 & . & 1 & . & 1 & . & . & 1 & 1 & 1 & 1 & 2 & . & . \\
 . & 1 & 4 & 2 & 2 & . & 3 & . & . & 1 & 1 & 2 & 2 & 2 & . & . \\
 . & 1 & 4 & 2 & 2 & . & 3 & . & . & 1 & 1 & 2 & 2 & 2 & . & . \\
 . & 2 & 6 & . & 2 & . & 3 & . & . & 2 & 2 & 2 & 2 & 3 & . & . \\
 . & . & 1 & . & . & . & . & . & . & . & . & . & . & . & . & . \\
 . & . & 1 & . & . & . & . & . & . & . & . & . & . & . & . & . \\
\end{array}
\right)
, \quad
g_{14}=&
\left(
\begin{array}{cccccccccccccccc}
 . & . & . & . & . & . & . & . & . & . & . & . & . & 1 & . & . \\
 . & . & 2 & 1 & 1 & . & 2 & . & . & . & . & 1 & 1 & 1 & . & . \\
 . & 2 & 9 & 3 & 4 & . & 6 & . & . & 2 & 2 & 4 & 4 & 3 & . & . \\
 . & 1 & 3 & . & 1 & 1 & . & . & . & 1 & 1 & 1 & 1 & 3 & 1 & 1 \\
 . & 1 & 4 & 1 & 1 & . & 2 & . & . & 1 & 1 & 1 & 1 & 2 & . & . \\
 . & . & . & 1 & . & . & . & . & . & . & . & . & . & 1 & . & . \\
 . & 2 & 6 & . & 2 & . & 3 & . & . & 2 & 2 & 2 & 2 & 3 & . & . \\
 . & . & . & . & . & . & . & . & . & . & . & . & . & 1 & . & . \\
 . & . & . & . & . & . & . & . & . & . & . & . & . & 1 & . & . \\
 . & . & 2 & 1 & 1 & . & 2 & . & . & . & . & 1 & 1 & 1 & . & . \\
 . & . & 2 & 1 & 1 & . & 2 & . & . & . & . & 1 & 1 & 1 & . & . \\
 . & 1 & 4 & 1 & 1 & . & 2 & . & . & 1 & 1 & 1 & 1 & 2 & . & . \\
 . & 1 & 4 & 1 & 1 & . & 2 & . & . & 1 & 1 & 1 & 1 & 2 & . & . \\
 1 & 1 & 3 & 3 & 2 & 1 & 3 & 1 & 1 & 1 & 1 & 2 & 2 & 3 & 1 & 1 \\
 . & . & . & 1 & . & . & . & . & . & . & . & . & . & 1 & . & . \\
 . & . & . & 1 & . & . & . & . & . & . & . & . & . & 1 & . & . \\
\end{array}
\right)
\end{align*}
}
  Notice  that $g_8$ alone would determine the missing ones. 
 Indeed $g_9$ is then obtained from the known $g_{8+9}$;
then, from the $g_6$ matrix, and up to a commutativity issue that will be discussed later,  we read  $g_{15}=g_6 . g_8$, and $g_{16}=g_6 . g_9$;
 then,  from the $g_2$ matrix, we read $g_{10}=g_2 . g_8$ and $g_{11}=g_2 . g_9$. In the same way $g_{16}$ is obtained from the known $g_{15+16}$; and, from the $g_2$ matrix, we read $g_{12}=g_2 . g_{15}$ and $g_{13}=g_2 . g_{16}$.
 The resulting family of twelve matrices obtained so far,  $\{g_1, g_2, g_3, g_4, g_5, g_6, g_7,g_{8+9}, g_{10+11}, g_{12+13}, g_{14}, g_{15+16}\}$, is a commutative family of rank $12$.   

In order to decouple the four pairs $(8,9)$, $(10,11)$, $(12,13)$, and $(15,16)$ we can use several ingredients.
First of all we notice that the module category with $16$ simple objects that we consider is a category of modules over an \'etale algebra (hence a commutative algebra) but the Grothendieck ring of $M$ itself need not be commutative. The already determined $8$ matrices $g_a$, with $a \in \{1,2,3,4,5,6,7,14\}$ commute, but we should not a priori impose commutation with, or among, those that are still to be determined (the four pairs). We may use the following:  1. Non-negativity of the matrix elements. 2. The subring generated by the ambichiral elements $a \in \{1, 6, 8, 9, 15, 16\}$ should be commutative: this follows from the fact that these are the ``modular points'' (braiding  is defined).
3. The rigidity condition (compatibility with duality) which, in simple terms, reads as follows: if $g_a$ is the self-fusion matrix representing a simple object $a$, then its transpose $g_a^T$ implements duality and should correspond to a simple object $a^\vee$ of $M$. 4. The compatibility equations used in \ref{compatibilityspecial} are a particular case of a more general compatibility between modular structure and self-fusion that holds for all the simple objects of $M$, not only for the identity. Namely, for all simple objects $a$ of $M$ and simple objects $n$ of $(F_4,6)$, one should have 
\begin{equation}
F_n . \,  g_a = \sum_b (F_n)_{a,b} \, g_b.
\label{compatibilitygeneral}
\end{equation}

Using these constraints we find two classes of solutions which, up to relabelling, lead to two solutions only.

We could not find a convincing argument that would allow us to make a definitive choice, and it may be that the determination of the structure constants describing the multiplication cannot be done without a thorough study of the multiplicative structure of the ring $O$ itself (because the self-fusion of $M$ is supposed to be inherited from the $O$ multiplication), something that we did not carry out.
So we leave this problem for the future (for future work.). 
The two solutions and their determination are  discussed in the Appendix.

The first is a non-commutative solution, which looks nice but should maybe be discarded because we could not find a way for it to obey simultaneously the rigidity condition and the general compatibility conditions with the $A$ action.
The algebra generated by the $g_a$'s, for this solution, is isomorphic to the direct sum of eight copies of the complex numbers, and four copies of the algebra of $2\times2$ complex matrices, stemming from the existence of four pairs  $(8,9)$, $(10,11)$, $(12,13)$, and $(15,16)$, of objects that could not be decoupled at the end of the previous discussion.

The other is a commutative solution that enjoys all the required properties but looks a bit trivial.
It is a system of $16$ commuting matrices $g_a$ (with non-negative integer entries, of course). There are eight self-dual symmetric matrices with real eigenvalues, and eight other matrices, namely the members of the four pairs previously described, that have complex eigenvalues.

Notice that a braiding of $M$ would imply commutativity of self-fusion, but $M$ is not braided, a fact that does not prevent commutativity of the self-fusion. Our analysis selects the commutative solution (because it obeys all the constraints) but we cannot totally exclude the other, so that, as already written, we take this as an open question. 
Notice that the determination, or lack thereof, of the multiplicative structure of the ring $O$ and of the subring $M$ does not affect the results described in the main body of the present paper, which concern the explicit description of the action of the fusion ring on its various modules in terms of matrices with non-negative integer entries (nimreps). This was the main objective of the present work.

 \section{The quantum module $M^\prime$}  
    \label{qModuleMprime}

This section is similar to the previous one but there are two main differences between them. 
The first is that, after having solved the intertwining equations discussed in section \ref{Intertwiners},
 getting rid of the remaining parameters in the restriction of the left and right intertwiners $V^L_f$ and $V^R_f$ to $M^\prime$ (i.e. to $M_1^\prime$  or to its copy $M_2^\prime$) is much more difficult than for their restriction to $M$ (i.e. to $M_1$ or to its copy $M_2$). The determination of these parameters therefore deserves a separate discussion. Another difference from the previous case is that $M^\prime$ is ``only'' a quantum module. The natural self-fusion structure inherited from the multiplication of $O$ does not naturally restrict to $M^\prime$.
 Although $M$ and $M^\prime$ are both associated with the same modular invariant, and with the same algebra of quantum symmetries $O$, they are inequivalent as ``nimreps'', and they describe two distinct module categories. This also justifies a specific discussion.
 
      \subsection{Solving constraints for the $M^\prime$ parameters}
      \label{annularmatricesMprimeParameters}

The four matrices $F_{0001}^\prime$, $F_{0010}^\prime$, $ F_{0100}^\prime$, $F_{1000}^\prime$ describing the module structure of the quantum module $M^\prime$ over the fusion algebra can be obtained by restriction of the left or right intertwining matrices to the third, or fourth, connected component of the graphs $V^L_f$ or $V^R_f$.       
After solving the various constraints described in the intertwiners section, one finds that the matrix elements of $t=F_{0001}^\prime$, for instance, still contain 29 parameters, whose possible ranges of values, taking into account the non-negative integer constraints, leads approximately to $2 \times 10^9$ possibilities. Something similar occurs with $z=F_{0010}^\prime$.
Finding compatible pairs $t$ and $z$ to build an acceptable $u$ to be plugged in the rational polynomials of degree $38$ discussed in sec.~\ref{FusionQpolynomials} becomes an issue.

Since $M$ and $M^\prime$ are both associated with (and determined by) the same modular invariant, the zeros of the minimal polynomials of the matrices $F(\lambda)$, or $F^\prime(\lambda)$, should match the values $S(\lambda,\mu)/S(0,\mu)$ where $S$ is the modular $S$ matrix and where $\mu$ runs over exponents.
 In particular, the matrices $F_{0001}$ and $F_{0001}^\prime$ are equivalent over the field of complex numbers (even though they turn out to be inequivalent as matrices with non-negative integer entries  describing distinct nimreps). Since they have the same minimal polynomials, one can use for $M^\prime$ the same reduced rational fusion polynomials that we discussed in sec.~\ref{polynomialsof0001} in the case of $M$. These are of degree $11$ in the variable $t=0001$, which is therefore primitive for both $M$ and $M^\prime$.

This is a simplifying remark because, from the known matrix $F_{0001}^\prime$,  one can use these  polynomials to obtain all the $F_{f}^\prime$ where $f$ is fundamental, then, using the $\ZZ$-polynomials, 
all the $F_{m}^\prime$, for all simple $m$. 
From the corresponding known matrices for $M$, from the known permutations exchanging the two copies of $M$ and the two copies of $M^\prime$, 
one obtains the left intertwiners $V^L_{0010}$ and $V^L_{0001}$ of $O$, so that one can use the rational polynomials of degree 38 in $u=2z+t$, which is a primitive element for ${\mathcal A}$, 
 to obtain all the left intertwiners (and also the right ones by using the appropriate permutation).

However, determining a ``good'' $F_{0001}^\prime$, which is the starting point of the above procedure, is not immediately obvious because choosing a random solution among those satisfying the constraints on the $29$ remaining parameters of $t$ generally leads to unacceptable results for the annular matrices $F_{m}^\prime$ (they will have matrix elements that  are not non-negative integers), although it gives an acceptable matrix $F_{0001}^\prime$. A refinement of the above method, sketched below, is needed.

   The minimal polynomials of $F(\lambda)$, for $\lambda$ equal to $0001$, $0010$, $0100$, $1000$, calculated for $M$,  are respectively given by
    \begin{align*}
&  \left(z^4-12 z^3+14 z^2+12 z+1\right) \left(z^4-3 z^3-z^2+3 z+1\right)^2 \left(z^4+3 z^3-z^2-3 z+1\right) \\
&  z^8 \left(z^4-42 z^3-51 z^2+72 z-9\right) \left(z^4+3 z^3-6 z^2-18 z-9\right)\\
&   \left(z^4-74 z^3-4 z^2+56 z+16\right) \left(z^4-14 z^3-34 z^2+41 z+1\right) \left(z^4+4 z^3-4 z^2-z+1\right)^2\\
&  z^8 \left(z^4-18 z^3+69 z^2-72 z-9\right) \left(z^4-3 z^3-21 z^2-27 z-9\right)
        \end{align*}
 Notice that the first is symmetric, up to alternating signs.    
 The corresponding ones for $M^\prime$, calculated for $F^\prime(\lambda)$ with the same $\lambda$, should be equal.
 Their coefficients are given by the Newton polynomials. They start as $-Tr[X], (Tr[X]^2 - Tr[X . X])/2, 
-(Tr[X]^3 - 3 Tr[X] Tr[X . X] + 2 Tr[X . X . X] )/6,
 \ldots$
 It turns out to be enough to impose equality of  $Tr[X]$ and $Tr[X.X]$ to those calculated using the known $X=F_{0001}$ (their respective values are $15$ and $149$) and $X=F_{0010}$ (the values are $39$ and $1887$).
 
 This leads to linear and quadratic equations for the remaining parameters of $F^\prime_{0001}$ and $F^\prime_{0010}$ that are simple enough to be solved. The various solutions obtained turn out to be equivalent by permuting the vertices. It is actually enough to do it for $F^\prime_{0001}$ only since we have rational polynomials at our disposal to get $F^\prime_{0010}$ and the others.
  Plugging those solutions into the chiral intertwiners $V^L_f$, for $f=0001$ and $f=0010$, and following the route sketched previously, one eventually obtains all the intertwiners $V^L_m$, for all irreps $m$ of $A$, without any parameter left, and by restriction all the matrices $F_m^\prime$ describing the $A$-module structure of~$M^\prime$.

    \subsubsection{$M^\prime$ as a module over the fusion ring}
          \label{annularmatricesMprime}
    
We display the annular matrices (nimreps) for the four fundamental irreps of the quantum module $M^\prime$ and their associated adjacency graphs while discarding self-loops and multiedges.

\begin{center}
\begin{minipage}[t]{0.62\textwidth}
\vspace{0pt}
{\footnotesize
\[
\setlength{\arraycolsep}{2.5pt}
 F_{0001}^\prime=
 \left(
\begin{array}{cccccccccccccccc}
 . & 1 & . & . & . & . & . & . & . & . & . & . & . & . & . & . \\
 1 & 2 & . & 1 & 1 & 1 & 1 & . & . & . & 2 & 2 & 2 & 1 & 1 & 1 \\
 . & . & . & 1 & . & . & . & . & . & . & . & . & . & . & . & . \\
 . & 1 & 1 & 1 & . & . & . & 1 & 1 & 1 & 1 & 1 & 1 & 1 & 1 & 1 \\
 . & 1 & . & . & 1 & . & . & . & . & 1 & 1 & . & 1 & 1 & 1 & . \\
 . & 1 & . & . & . & 1 & . & . & 1 & . & . & 1 & 1 & . & 1 & 1 \\
 . & 1 & . & . & . & . & 1 & 1 & . & . & 1 & 1 & . & 1 & . & 1 \\
 . & . & . & 1 & . & . & 1 & . & . & . & 1 & 1 & . & . & . & . \\
 . & . & . & 1 & . & 1 & . & . & . & . & . & 1 & 1 & . & . & . \\
 . & . & . & 1 & 1 & . & . & . & . & . & 1 & . & 1 & . & . & . \\
 . & 2 & . & 1 & 1 & . & 1 & 1 & . & 1 & 2 & 1 & 1 & 1 & 1 & 1 \\
 . & 2 & . & 1 & . & 1 & 1 & 1 & 1 & . & 1 & 2 & 1 & 1 & 1 & 1 \\
 . & 2 & . & 1 & 1 & 1 & . & . & 1 & 1 & 1 & 1 & 2 & 1 & 1 & 1 \\
 . & 1 & . & 1 & 1 & . & 1 & . & . & . & 1 & 1 & 1 & 1 & . & . \\
 . & 1 & . & 1 & 1 & 1 & . & . & . & . & 1 & 1 & 1 & . & 1 & . \\
 . & 1 & . & 1 & . & 1 & 1 & . & . & . & 1 & 1 & 1 & . & . & 1 \\
\end{array}
\right)
\]
}
\end{minipage}
\hfill
\begin{minipage}[t]{0.32\textwidth}
\vspace{0pt}
\centering
\includegraphics[width=\textwidth]{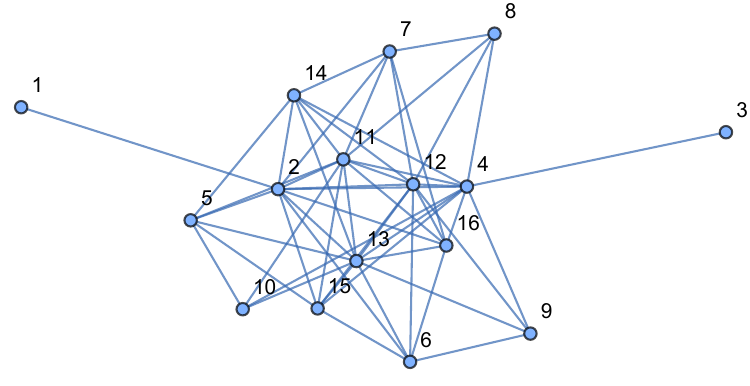}
\captionof{figure}{Adjacency graph of the matrix of the generator $0001$ of $(F_4,6)$ acting on the quantum module $M^\prime$.}
\label{mp0001}
\end{minipage}
\end{center}
\begin{center}
\begin{minipage}[t]{0.62\textwidth}
\vspace{0pt}
{\footnotesize
\[
\setlength{\arraycolsep}{2.5pt}
 F_{0010}^\prime=
 \left(
\begin{array}{cccccccccccccccc}
 . & . & . & . & . & . & . & . & . & . & 1 & 1 & 1 & 1 & 1 & 1 \\
 . & 9 & . & 6 & 4 & 4 & 4 & 2 & 2 & 2 & 7 & 7 & 7 & 4 & 4 & 4 \\
 . & . & . & . & . & . & . & . & . & . & 1 & 1 & 1 & . & . & . \\
 . & 6 & . & 3 & 2 & 2 & 2 & 2 & 2 & 2 & 4 & 4 & 4 & 3 & 3 & 3 \\
 . & 4 & . & 2 & 1 & 1 & 1 & 1 & 1 & 1 & 3 & 3 & 3 & 2 & 2 & 2 \\
 . & 4 & . & 2 & 1 & 1 & 1 & 1 & 1 & 1 & 3 & 3 & 3 & 2 & 2 & 2 \\
 . & 4 & . & 2 & 1 & 1 & 1 & 1 & 1 & 1 & 3 & 3 & 3 & 2 & 2 & 2 \\
 . & 2 & . & 2 & 1 & 1 & 1 & . & . & . & 2 & 2 & 2 & 1 & 1 & 1 \\
 . & 2 & . & 2 & 1 & 1 & 1 & . & . & . & 2 & 2 & 2 & 1 & 1 & 1 \\
 . & 2 & . & 2 & 1 & 1 & 1 & . & . & . & 2 & 2 & 2 & 1 & 1 & 1 \\
 1 & 7 & 1 & 4 & 3 & 3 & 3 & 2 & 2 & 2 & 6 & 6 & 6 & 4 & 4 & 4 \\
 1 & 7 & 1 & 4 & 3 & 3 & 3 & 2 & 2 & 2 & 6 & 6 & 6 & 4 & 4 & 4 \\
 1 & 7 & 1 & 4 & 3 & 3 & 3 & 2 & 2 & 2 & 6 & 6 & 6 & 4 & 4 & 4 \\
 1 & 4 & . & 3 & 2 & 2 & 2 & 1 & 1 & 1 & 4 & 4 & 4 & 2 & 2 & 2 \\
 1 & 4 & . & 3 & 2 & 2 & 2 & 1 & 1 & 1 & 4 & 4 & 4 & 2 & 2 & 2 \\
 1 & 4 & . & 3 & 2 & 2 & 2 & 1 & 1 & 1 & 4 & 4 & 4 & 2 & 2 & 2 \\
\end{array}
\right)
\]
}
\end{minipage}
\hfill
\begin{minipage}[t]{0.32\textwidth}
\vspace{0pt}
\centering
\includegraphics[width=\textwidth]{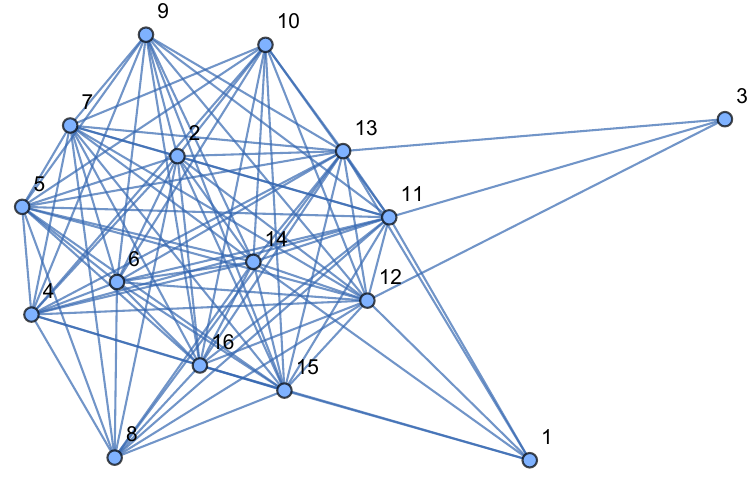}
\captionof{figure}{Adjacency graph of the matrix of the generator $0010$ of $(F_4,6)$ acting on the quantum module $M^\prime$.}
\label{mp0010}
\end{minipage}
\end{center}
\begin{center}
\begin{minipage}[t]{0.62\textwidth}
\vspace{0pt}
{\footnotesize
\[
\setlength{\arraycolsep}{2.5pt}
 F_{0100}^\prime=
 \left(
\begin{array}{cccccccccccccccc}
 . & 2 & . & 2 & 1 & 1 & 1 & . & . & . & 1 & 1 & 1 & . & . & . \\
 2 & 15 & 2 & 7 & 5 & 5 & 5 & 5 & 5 & 5 & 12 & 12 & 12 & 9 & 9 & 9 \\
 . & 2 & . & . & . & . & . & 1 & 1 & 1 & . & . & . & 1 & 1 & 1 \\
 2 & 7 & . & 8 & 5 & 5 & 5 & . & . & . & 9 & 9 & 9 & 3 & 3 & 3 \\
 1 & 5 & . & 5 & 2 & 4 & 4 & 1 & 1 & . & 5 & 7 & 5 & 2 & 2 & 3 \\
 1 & 5 & . & 5 & 4 & 2 & 4 & 1 & . & 1 & 7 & 5 & 5 & 3 & 2 & 2 \\
 1 & 5 & . & 5 & 4 & 4 & 2 & . & 1 & 1 & 5 & 5 & 7 & 2 & 3 & 2 \\
 . & 5 & 1 & . & 1 & 1 & . & 2 & 3 & 3 & 2 & 2 & 3 & 3 & 4 & 3 \\
 . & 5 & 1 & . & 1 & . & 1 & 3 & 2 & 3 & 3 & 2 & 2 & 4 & 3 & 3 \\
 . & 5 & 1 & . & . & 1 & 1 & 3 & 3 & 2 & 2 & 3 & 2 & 3 & 3 & 4 \\
 1 & 12 & . & 9 & 5 & 7 & 5 & 2 & 3 & 2 & 10 & 11 & 11 & 5 & 6 & 6 \\
 1 & 12 & . & 9 & 7 & 5 & 5 & 2 & 2 & 3 & 11 & 10 & 11 & 6 & 6 & 5 \\
 1 & 12 & . & 9 & 5 & 5 & 7 & 3 & 2 & 2 & 11 & 11 & 10 & 6 & 5 & 6 \\
 . & 9 & 1 & 3 & 2 & 3 & 2 & 3 & 4 & 3 & 5 & 6 & 6 & 5 & 5 & 5 \\
 . & 9 & 1 & 3 & 2 & 2 & 3 & 4 & 3 & 3 & 6 & 6 & 5 & 5 & 5 & 5 \\
 . & 9 & 1 & 3 & 3 & 2 & 2 & 3 & 3 & 4 & 6 & 5 & 6 & 5 & 5 & 5 \\
\end{array}
\right)
\]
}
\end{minipage}
\hfill
\begin{minipage}[t]{0.32\textwidth}
\vspace{0pt}
\centering
\includegraphics[width=\textwidth]{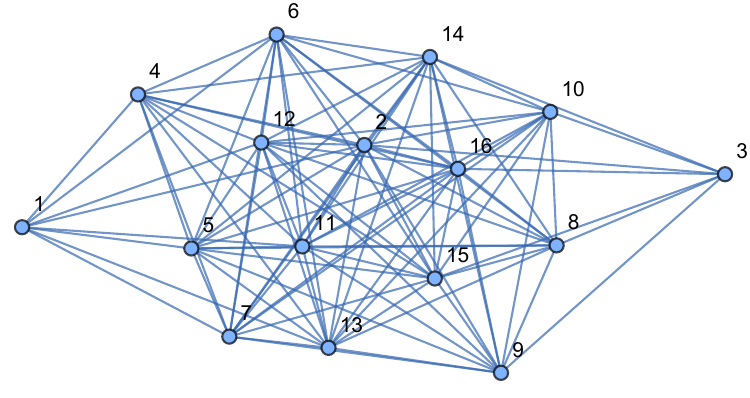}
\captionof{figure}{Adjacency graph of the matrix of the generator $0100$ of $(F_4,6)$ acting on the quantum module $M^\prime$.}
\label{mp0100}
\end{minipage}
\end{center}

\begin{center}
\begin{minipage}[t]{0.62\textwidth}
\vspace{0pt}
{\footnotesize
\[
\setlength{\arraycolsep}{2.5pt}
 F_{1000}^\prime=
 \left(
\begin{array}{cccccccccccccccc}
 . & . & . & . & 1 & 1 & 1 & . & . & . & . & . & . & . & . & . \\
 . & 3 & . & . & 1 & 1 & 1 & 1 & 1 & 1 & 2 & 2 & 2 & 2 & 2 & 2 \\
 . & . & . & . & . & . & . & 1 & 1 & 1 & . & . & . & . & . & . \\
 . & . & . & 3 & 1 & 1 & 1 & . & . & . & 2 & 2 & 2 & . & . & . \\
 1 & 1 & . & 1 & 1 & 1 & 1 & . & . & . & 1 & 1 & 1 & . & . & . \\
 1 & 1 & . & 1 & 1 & 1 & 1 & . & . & . & 1 & 1 & 1 & . & . & . \\
 1 & 1 & . & 1 & 1 & 1 & 1 & . & . & . & 1 & 1 & 1 & . & . & . \\
 . & 1 & 1 & . & . & . & . & 1 & 1 & 1 & . & . & . & 1 & 1 & 1 \\
 . & 1 & 1 & . & . & . & . & 1 & 1 & 1 & . & . & . & 1 & 1 & 1 \\
 . & 1 & 1 & . & . & . & . & 1 & 1 & 1 & . & . & . & 1 & 1 & 1 \\
 . & 2 & . & 2 & 1 & 1 & 1 & . & . & . & 2 & 2 & 2 & 1 & 1 & 1 \\
 . & 2 & . & 2 & 1 & 1 & 1 & . & . & . & 2 & 2 & 2 & 1 & 1 & 1 \\
 . & 2 & . & 2 & 1 & 1 & 1 & . & . & . & 2 & 2 & 2 & 1 & 1 & 1 \\
 . & 2 & . & . & . & . & . & 1 & 1 & 1 & 1 & 1 & 1 & 1 & 1 & 1 \\
 . & 2 & . & . & . & . & . & 1 & 1 & 1 & 1 & 1 & 1 & 1 & 1 & 1 \\
 . & 2 & . & . & . & . & . & 1 & 1 & 1 & 1 & 1 & 1 & 1 & 1 & 1 \\
\end{array}
\right)
\]
}
\end{minipage}
\hfill
\begin{minipage}[t]{0.32\textwidth}
\vspace{0pt}
\centering
\includegraphics[width=\textwidth]{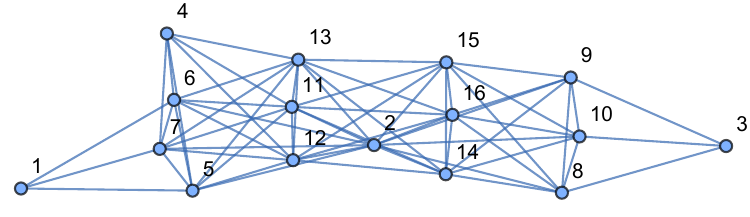}
\captionof{figure}{Adjacency graph of the matrix of the generator $1000$ of $(F_4,6)$ acting on the quantum module $M^\prime$.}
\label{mp1000}
\end{minipage}
\end{center}

      \subsection{Quantum dimensions for $M^\prime$ and other numerical results}
      
            \subsubsection{Perron-Frobenius and a normalization issue}
    \label{qprimenormalisation}       
    Using, for instance, $F_{0001}^\prime$, one finds that the Perron-Frobenius eigenvalue (dominant eigenvalue) is  equal to $3+\sqrt{5}+\sqrt{3 \left(5+2 \sqrt{5}\right)}$. The corresponding eigenvector is also well-defined, up to scale, and its components give unnormalised quantum dimensions for the simple objects of $M^\prime$.
    
Unlike $M$, there is no simple object of  $M^\prime$ that can play the role of a unit, and there is no special space of sections ${\Gamma^\prime_1}$ that can be used for normalization.
Therefore, from the study of $M^\prime$ alone, there is no marked vertex that could fix the scale of quantum dimensions. However, we can mark some arbitrary point and decide that its quantum dimension is $1$, for instance one of the two "extremal" vertices  labelled $3$ and $17$ on the Ocneanu graph of $V^L$ (see fig.~\ref{FL0001}),  or labelled  $1$ and $3$ on its restriction, the adjacency graph of $M^\prime$ (see fig.~\ref{mp0001}).
For definiteness, we choose the vertex labelled $1$ on the graph of $M^\prime$ (it gives the smallest numerical value for $q^\prime$), but one should remember that this $1$ is only a label; it does not refer to a unit.
 One could use another annular matrix, for instance $F^\prime(1,0,0,0)$, and obtain the same results.   
 {\footnotesize
\[
\begin{array}{@{}l@{\quad}l@{}}
q^\prime_1 =   1,&
   
q^\prime_2 =   3+\sqrt{5}+\sqrt{15+6 \sqrt{5}}, \\[3pt]
   
q^\prime_3 =   \frac{1}{2}
   \left(\sqrt{5}-1\right),&
   
q^\prime_4 =      1+\sqrt{5}+\sqrt{\frac{3}{2}
   \left(5+\sqrt{5}\right)}, \\[3pt]
   
q^\prime_5 =   \frac{1}{2} \left(3+\sqrt{5}+\sqrt{\frac{20}{3}+\frac{8
   \sqrt{5}}{3}}\right),&
   
q^\prime_6 =   \frac{1}{2} \left(3+\sqrt{5}+\sqrt{\frac{20}{3}+\frac{8
   \sqrt{5}}{3}}\right), \\[3pt]
   
q^\prime_7 =   \frac{1}{2} \left(3+\sqrt{5}+\sqrt{\frac{20}{3}+\frac{8
   \sqrt{5}}{3}}\right),&
   
q^\prime_8 =   \frac{1}{6} \left(3+3 \sqrt{5}+\sqrt{6
   \left(5+\sqrt{5}\right)}\right), \\[3pt]
   
q^\prime_9 =   \frac{1}{6} \left(3+3 \sqrt{5}+\sqrt{6
   \left(5+\sqrt{5}\right)}\right),&
   
q^\prime_{10}  =   \frac{1}{6} \left(3+3 \sqrt{5}+\sqrt{6
   \left(5+\sqrt{5}\right)}\right), \\[3pt]
   
q^\prime_{11}  =   2+\sqrt{5}+\sqrt{\frac{1}{6} \left(65+29
   \sqrt{5}\right)},&
   
q^\prime_{12}  =   2+\sqrt{5}+\sqrt{\frac{1}{6} \left(65+29
   \sqrt{5}\right)}, \\[3pt]
   
q^\prime_{13}  =   2+\sqrt{5}+\sqrt{\frac{1}{6} \left(65+29
   \sqrt{5}\right)},&
   
q^\prime_{14}  =   \frac{1}{6} \left(9+3 \sqrt{5}+\sqrt{150+66
   \sqrt{5}}\right), \\[3pt]
   
q^\prime_{15}  =    \frac{1}{6} \left(9+3 \sqrt{5}+\sqrt{150+66
   \sqrt{5}}\right),&
   
   q^\prime_{16}  = \frac{1}{6} \left(9+3 \sqrt{5}+\sqrt{150+66 \sqrt{5}}\right)
   
   \end{array}
   \]
   }

\subsubsection{$\vert{\mathcal M^\prime}\vert$ from $q_a^\prime$ (summing squares)}

We define the quantum mass of $\mathcal M^\prime$ as the sum of squares  of its quantum dimensions, normalized as above.
One finds  $\vert{\mathcal M^\prime}\vert=45 \left(87 \sqrt{5}+4 \sqrt{2046 \sqrt{5}+4575}+195\right)$.

       \subsection{Essential matrices for $M^\prime$}
       
    As for $M$, we could define rectangular essential matrices for $M^\prime$ by setting $(Ess^\prime_a)_{(n,b)}= (F^\prime_n)_{(a,b)}$ and,
    choosing some simple object $\star$ for $\mathcal M^\prime$ (that we mark),  we could  set  $\Gamma_a^\star = \sum _{n  \uparrow \Gamma_a^\star} \,  n$, where the sum runs over those simple objects $n$ of  ${\mathcal A}$ that appear in the column $a$ of the essential matrix $Ess^\prime_{{\star}}$.
    As for $M$, we could set  $\vert {\Gamma^\prime}_a^\star \vert = \sum _{n  \uparrow {\Gamma^\prime}_a^\star} \,  \text{qdim}(n)$, and calculate these quantities.
For $M$, we had $q_a = \vert \Gamma_a \vert  /    \vert \Gamma_1 \vert$, so that the list of $\vert \Gamma_a \vert$'s was  a scaled version of the list of quantum dimensions. This is also the case here, except that the unit $1$ is replaced by an arbitrarily chosen marked object $\star$.
    In other words, the list of these $\vert {\Gamma^\prime}_a^\star \vert$'s coincides, up to scale, with the list of the $q^\prime_{a}$'s given in     sec.~\ref{qprimenormalisation} and the list of $\vert \Gamma_a^\star \vert  /    \vert \Gamma_\star^\star \vert$'s is a list of quantum dimensions normalized in such a way that the chosen object $\star$ has q-dimension $1$. It does not seem that there is a way to remove this arbitrariness from the study of  $M^\prime$ alone, so we shall refrain from elaborating further on this matter.

 \section{Other results and suggestions for further study}
 \label{varia}
 
 \subsection{About toric matrices with two twists}
 In the previous sections we have obtained all the toric matrices $W_x$ with one twist (one defect line), i.e. all the $(V_{m,n})_{x,1}$ or $(V_{m,n})_{1,x}$, and we have described all the left and right intertwiners $(V_{m,1})_{x,y}$ and $(V_{1,m})_{x,y}$ but not the full family of the coefficients $(V_{m,n})_{x,y}$. This could be done and should not be too difficult.
 
 \subsection{Ring structure of $O$ and the generators of quantum symmetries $O_x$}
 We did not study this, but doing so requires solving the previous problem first.

 \subsection{Dual annular matrices}
 The action of $O$ on $M$ is described by the so-called dual annular matrices $S_x$. 
 This study requires knowledge of the $O_x$.

  \subsection{Bialgebra structure and generalized $6J$ symbols}

 One can associate a quantum groupoid  ${\mathcal B}$  to the module category ${\mathcal M}$.
  It is a finite-dimensional weak Hopf algebra.
  One can think of the algebra ${\mathcal B}$  as a direct sum of $r_A=39$ simple matrix components, and of its
dual, the algebra ${\mathcal B}^\star$, as a direct sum of $r_O=64$ simple matrix components. The dimensions $d_n$ and $d_x$ of these
blocks, called horizontal or vertical dimensions (or dimensions of the generalized spaces of essential paths, or spaces of admissible triples or generalized triangles, etc.), can be obtained from the annular matrices $F_n$ and from the dual annular matrices $S_x$.

Using labels $a,b,\ldots$ for the simple objects of ${\mathcal M}$, $m,n,\ldots$ for those of ${\mathcal A}$, and $x,y, \ldots$ for those of ${\mathcal O}$, one can build a particular basis of ${\mathcal B}$ indexed by double triangles $((a,n,b),(c,n,d))$ when $(F_n)_{a,b} \neq 0$ and $(F_n)_{c,d} \neq 0$, a particular basis of ${\mathcal B}^\star$ indexed by double triangles $((a,x,c),(b,x,d))$ when $(S_x)_{a,c} \neq 0$ and $(S_x)_{b,d} \neq 0$. The bialgebra pairing between these two bases is encoded by a collection of generalized $6J$ symbols $\left({}^{a,n,b}_{c,x,d}\right)$ that one can picture as tetrahedra with four edges of ${\mathcal M}$-type (labelled $a,b,c,d$), one edge of ${\mathcal A}$-type labelled $n$, and the opposite edge of ${\mathcal O}$-type labelled $x$. We only  give the dimensions of the horizontal blocks (the $n$-blocks), as well as the total dimension $d_{\mathcal B}$.

Horizontal dimensions: $d_n =$
{\footnotesize
$\{16, 95, 376, 363, 117, 819, 1590, 644, 739, 967, 2917, 1932, 2321, 1863, 1724, 374, \\
 609, 2580, 3304, 3161, 2541, 2273, 4170, 1516, 1348, 175, 925, 1802, 1036, 1629, 2580, 705, 1193, 3384, 2178, 2424,\\ 
 1569, 1407, 268\}.$}
The total horizontal dimension is $\sum_n d_n =  59634$ and the dimension of the bialgebra is $d_{\mathcal B} = \sum_n d_n^2 = 133686600$.
 Finding the list of vertical dimensions $d_x$ would require a detailed description of the multiplicative structure of $O$ (the structure coefficients given by the $O_x$ matrices), and a description of the dual annular matrices $S_x$. This would, in particular, allow one to check the equality $\sum_n d_n^2 =\sum_x d_x^2 $.

 \section*{Appendix}
 \label{appendix}
 
 Discussion of the self-fusion on $M$.
 This appendix should be considered as a continuation of sec.~\ref{selffusionofM}.
   
 The first step is to start with the matrix $g_{8+9}$ and replace its non-zero matrix elements by $28$ arbitrary parameters. This is an ansatz for the matrix $g_8$ (because of the non-negativity of the coefficients).
 Since $8$ and $9$ are ambichiral,  they should commute, so $8$ should also commute with the known sum $g_{8+9}$. This determines several parameters. The matrix $g_8$ should also commute with $g_6$ because they are both ambichiral; this determines another set of parameters. Then, one sets $g_{15} = g_6 . g_8$ (or $g_8. g_6$ since they commute). At this step, one can check that the ambichiral subalgebra is commutative: all six $g_a$ for $a\in \{1, 6, 8, 9, 15, 16\}$ commute.
 The $g_8$ matrix reads: 
  {\footnotesize $$
  \setlength{\arraycolsep}{2.pt}
  \renewcommand{\arraystretch}{0.8}
  \left(
\begin{array}{cccccccccccccccc}
 . & . & . & . & . & . & . & a & b & . & . & . & . & . & . & . \\
 . & . & . & . & . & . & . & . & . & c & d & . & . & . & . & . \\
 . & . & e & . & . & . & . & . & . & . & . & . & . & . & . & . \\
 . & . & . & f & . & . & . & . & . & . & . & . & . & . & . & . \\
 . & . & . & . & . & . & . & . & . & . & . & c & d & . & . & . \\
 . & . & . & . & . & . & . & . & . & . & . & . & . & . & a & b \\
 . & . & . & . & . & . & e & . & . & . & . & . & . & . & . & . \\
 b & . & . & . & . & . & . & . & a & . & . & . & . & . & . & . \\
 a & . & . & . & . & . & . & b & . & . & . & . & . & . & . & . \\
 . & d & . & . & . & . & . & . & . & . & c & . & . & . & . & . \\
 . & c & . & . & . & . & . & . & . & d & . & . & . & . & . & . \\
 . & . & . & . & d & . & . & . & . & . & . & . & c & . & . & . \\
 . & . & . & . & c & . & . & . & . & . & . & d & . & . & . & . \\
 . & . & . & . & . & . & . & . & . & . & . & . & . & f & . & . \\
 . & . & . & . & . & b & . & . & . & . & . & . & . & . & . & a \\
 . & . & . & . & . & a & . & . & . & . & . & . & . & . & b & . \\
\end{array}
\right)$$}
where the six (still unknown) parameters $a,b,c,d,e,f$ should be non-negative integers.
The rigidity condition implies, in particular, that $g_8^T$ and $g_{15}^T$ should correspond to simple objects. The only possible candidates are $g_8^T=g_9$  and  $g_{15}^T=g_{16}$.
These equalities (which imply $g_{16}=g_6 . g_9$) impose relations between the remaining coefficients: 
the following combinations of parameters should vanish:   ${1 - a - b, 1 - c - d, 2 - 2 e, 2 - 2 f}$, where $a,b,c,d,e,f$ are non-negative integers.

There are four solutions: $abcdef=101011$ or $abcdef = 010111$, which lead to two equivalent commutative solutions, and $abcdef = 100111$ or $abcdef = 011011$, which lead to two equivalent non-commutative solutions.
We find that the last two, which are quite interesting, should probably be discarded, because some rigidity constraints cannot be satisfied at the same time as the equations that describe the compatibility between module structure and self-fusion (at least we did not find a way to make it work). Although ultimately discarded, we shall describe this class of solutions in the next subsection.

  \subsubsection*{Self-fusion: The (questionable) non-commutative solution}  
    We take $abcdef = 100111$. The resulting $g_8$ does not commute with $g_2$.
    From the compatibility equation between the $g_a$ and $(F_n)_{1,b}$ one finds that $g_{10}$ (resp.  $g_{11}$) is obtained from $g_{8}$ (resp.  $g_{9}$) by multiplication by $g_2$. In a similar way, $g_{12}$ and $g_{13}$ are obtained, respectively, from $g_{15}$ and $g_{16}$ by multiplication by $g_2$.
     Since $g_2$ and $g_8$ (or $g_2$ and $g_{15}$) do not commute, one may try several possibilities. 
     
     If we set $g_{10} = g_8 . g_2$,  $g_{11} = g_9 . g_2$, and  $g_{12} = g_{15} . g_2$, $g_{13} = g_{16} . g_2$,
   we can check that $g_{10}+ g_{11}=g_{10+11}$, and that $g_{12}+ g_{13}=g_{12+13}$, as it should. However, we have a rigidity problem: the transposes of $g_{10}$, $g_{11}$, $g_{12}$, $g_{13}$ do not correspond to (self-fusion matrices of) any simple object of $M$. For instance $g_{10}^T = g_2. g_9$ (which is not $g_{11}$). This solution also leads to incompatibility between modular and self-fusion structure: the special constraint (\ref{compatibilityspecial})  holds, but not the general one  (\ref{compatibilitygeneral}); one may notice that a twisted version of the latter holds, namely \ref{compatibilitygeneraltwisted} where $g_a$ appears to the left of $F_n$.
\begin{equation}
g_a  . F_n = \sum_b (F_n)_{a,b} \, g_b.
\label{compatibilitygeneraltwisted}
\end{equation}  
    
We could  change the position of $g_2$  in the previous definitions, \ie set    
 $g_{10} = g_2.g_8$,  $g_{11} = g_2.g_9$, and  $g_{12} = g_2.g_{15}$, $g_{13} = g_2.g_{16}$. The obtained solution is then compatible with the usual version of the compatibility constraint (\ref{compatibilitygeneral}), but the rigidity problem is still present: the transpose of $g_{10}$, $g_{11}$, $g_{12}$, $g_{13}$ do not correspond to any simple object.
    
 Finally one could try a mixed version where one defines   $g_{10} = g_8. g_2$,  $g_{11} = g_2.g_9$,  and  $g_{12} = g_{15}.g_2$, $g_{13} = g_2.g_{16}$. This solves the rigidity problem: $g_{10}$ is now the transpose of $g_{11}$, and $g_{12}$ is the transpose of $g_{13}$. Unfortunately, none of the compatibility equations (the usual one or the twisted one) are satisfied, and, even worse, the sums  
 $g_{10}+ g_{11}$ and $g_{12}+ g_{13}$ do not reproduce the known $g_{10+11}$ and $g_{12+13}$.

The conclusion is that, despite several nice features (the algebra generated by the $g_a$'s, for the above solutions, is isomorphic to the direct sum of eight copies of the complex numbers and four copies of the algebra of $2\times2$ complex matrices), these solutions fail to describe the self-fusion ring associated with $M$.

  The reader can check that taking $abcdef = 010111$, instead of $100111$, leads to a very similar discussion and to the same negative conclusion.
    
 \subsubsection*{Self-fusion: The commutative solution} 
We take $abcdef = 101011$. The resulting matrix $g_8$, as well as its transpose $g_9$ and the matrices $g_{15}$ and $g_{16}$, now commutes with $g_2$.
One has then to set $g_{10} = g_2.g_8$,  $g_{11} = g_2.g_9$, and  $g_{12} = g_2.g_{15}$, $g_{13} = g_2.g_{16}$, where the order no longer matters, and it turns out that all the matrices $g_a$ commute.
The special version of the compatibility equations (\ref{compatibilityspecial}) holds, as well as the general version (\ref{compatibilitygeneral}). The rigidity constraints are obeyed: the matrices $g_1, g_2, g_3, g_4, g_5, g_6, g_7, g_{14}$ are symmetric (they correspond to self-dual objects), and the pairs $(g_8, g_9)$, $(g_{10}, g_{11})$, $(g_{12}, g_{13})$, $(g_{15}, g_{16})$ are pairs of mutually transposed matrices, meaning that $g_{10}^T= g_{11}$, etc.  These four pairs are not all on the same footing, since $(g_8, g_9)$ and $(g_{15}, g_{16})$ are pairs of ambichiral objects, whereas the members of the other two pairs are not.
This system of $16$ commuting matrices $g_a$ can be simultaneously diagonalised: 
There are eight self-dual objects  described by symmetric matrices with real eigenvalues, whereas the other eight matrices, namely the members of the four pairs just described, have complex eigenvalues.
Recall that the ambichiral objects are labelled $1, 6, 8, 9, 15, 16$, and that their self-fusion matrices generate a subring of the commutative ring $M$.

Taking $abcdef = 010111$, rather than $101011$, leads to a distinct commutative solution, in the sense that the matrices $g_8, g_9, g_{10}, g_{11}, g_{12}, g_{13}, g_{15}, g_{16}$ are different. However, one can check that the algebras generated by these two families of self-fusion matrices are isomorphic. In other words, we obtain two possible realisations of the same self-fusion structure.

We already gave explicitly the matrices $g_a$, for $a$ in $\{1,2,3,4,5,6,7,14\}$. We now give those obtained for $a$ belonging to  $\{8,9,10,11,12,13,15, 16\}$, for the first commutative solution previously discussed.

{\scriptsize
\setlength{\arraycolsep}{2.5pt}
 \renewcommand{\arraystretch}{0.6}
\begin{align*}
g_8=& 
\left(
\begin{array}{cccccccccccccccc}
 . & . & . & . & . & . & . & 1 & . & . & . & . & . & . & . & . \\
 . & . & . & . & . & . & . & . & . & 1 & . & . & . & . & . & . \\
 . & . & 1 & . & . & . & . & . & . & . & . & . & . & . & . & . \\
 . & . & . & 1 & . & . & . & . & . & . & . & . & . & . & . & . \\
 . & . & . & . & . & . & . & . & . & . & . & 1 & . & . & . & . \\
 . & . & . & . & . & . & . & . & . & . & . & . & . & . & 1 & . \\
 . & . & . & . & . & . & 1 & . & . & . & . & . & . & . & . & . \\
 . & . & . & . & . & . & . & . & 1 & . & . & . & . & . & . & . \\
 1 & . & . & . & . & . & . & . & . & . & . & . & . & . & . & . \\
 . & . & . & . & . & . & . & . & . & . & 1 & . & . & . & . & . \\
 . & 1 & . & . & . & . & . & . & . & . & . & . & . & . & . & . \\
 . & . & . & . & . & . & . & . & . & . & . & . & 1 & . & . & . \\
 . & . & . & . & 1 & . & . & . & . & . & . & . & . & . & . & . \\
 . & . & . & . & . & . & . & . & . & . & . & . & . & 1 & . & . \\
 . & . & . & . & . & . & . & . & . & . & . & . & . & . & . & 1 \\
 . & . & . & . & . & 1 & . & . & . & . & . & . & . & . & . & . \\
\end{array}
\right)
\qquad
g_9=&
\left(
\begin{array}{cccccccccccccccc}
 . & . & . & . & . & . & . & . & 1 & . & . & . & . & . & . & . \\
 . & . & . & . & . & . & . & . & . & . & 1 & . & . & . & . & . \\
 . & . & 1 & . & . & . & . & . & . & . & . & . & . & . & . & . \\
 . & . & . & 1 & . & . & . & . & . & . & . & . & . & . & . & . \\
 . & . & . & . & . & . & . & . & . & . & . & . & 1 & . & . & . \\
 . & . & . & . & . & . & . & . & . & . & . & . & . & . & . & 1 \\
 . & . & . & . & . & . & 1 & . & . & . & . & . & . & . & . & . \\
 1 & . & . & . & . & . & . & . & . & . & . & . & . & . & . & . \\
 . & . & . & . & . & . & . & 1 & . & . & . & . & . & . & . & . \\
 . & 1 & . & . & . & . & . & . & . & . & . & . & . & . & . & . \\
 . & . & . & . & . & . & . & . & . & 1 & . & . & . & . & . & . \\
 . & . & . & . & 1 & . & . & . & . & . & . & . & . & . & . & . \\
 . & . & . & . & . & . & . & . & . & . & . & 1 & . & . & . & . \\
 . & . & . & . & . & . & . & . & . & . & . & . & . & 1 & . & . \\
 . & . & . & . & . & 1 & . & . & . & . & . & . & . & . & . & . \\
 . & . & . & . & . & . & . & . & . & . & . & . & . & . & 1 & . \\
\end{array}
\right)
\\
\\
g_{10}=&
\left(
\begin{array}{cccccccccccccccc}
 . & . & . & . & . & . & . & . & . & 1 & . & . & . & . & . & . \\
 . & . & 1 & 1 & . & . & 1 & 1 & . & 1 & . & 1 & . & . & . & . \\
 . & 1 & 4 & 2 & 2 & . & 3 & . & . & 1 & 1 & 2 & 2 & 2 & . & . \\
 . & 1 & 2 & . & . & . & . & . & . & 1 & 1 & . & . & 1 & . & . \\
 . & . & 2 & . & . & . & 1 & . & . & 1 & . & 2 & . & 1 & 1 & . \\
 . & . & . & . & . & . & . & . & . & . & . & 1 & . & . & . & . \\
 . & 1 & 3 & . & 1 & . & 1 & . & . & 1 & 1 & 1 & 1 & 2 & . & . \\
 . & . & . & . & . & . & . & . & . & . & 1 & . & . & . & . & . \\
 . & 1 & . & . & . & . & . & . & . & . & . & . & . & . & . & . \\
 . & . & 1 & 1 & . & . & 1 & . & 1 & . & 1 & . & 1 & . & . & . \\
 1 & 1 & 1 & 1 & 1 & . & 1 & . & . & . & . & . & . & . & . & . \\
 . & . & 2 & . & . & . & 1 & . & . & . & 1 & . & 2 & 1 & . & 1 \\
 . & 1 & 2 & . & 2 & 1 & 1 & . & . & . & . & . & . & 1 & . & . \\
 . & . & 2 & 1 & 1 & . & 2 & . & . & . & . & 1 & 1 & 1 & . & . \\
 . & . & . & . & . & . & . & . & . & . & . & . & 1 & . & . & . \\
 . & . & . & . & 1 & . & . & . & . & . & . & . & . & . & . & . \\
\end{array}
\right)
\qquad
g_{11}=&
\left(
\begin{array}{cccccccccccccccc}
 . & . & . & . & . & . & . & . & . & . & 1 & . & . & . & . & . \\
 . & . & 1 & 1 & . & . & 1 & . & 1 & . & 1 & . & 1 & . & . & . \\
 . & 1 & 4 & 2 & 2 & . & 3 & . & . & 1 & 1 & 2 & 2 & 2 & . & . \\
 . & 1 & 2 & . & . & . & . & . & . & 1 & 1 & . & . & 1 & . & . \\
 . & . & 2 & . & . & . & 1 & . & . & . & 1 & . & 2 & 1 & . & 1 \\
 . & . & . & . & . & . & . & . & . & . & . & . & 1 & . & . & . \\
 . & 1 & 3 & . & 1 & . & 1 & . & . & 1 & 1 & 1 & 1 & 2 & . & . \\
 . & 1 & . & . & . & . & . & . & . & . & . & . & . & . & . & . \\
 . & . & . & . & . & . & . & . & . & 1 & . & . & . & . & . & . \\
 1 & 1 & 1 & 1 & 1 & . & 1 & . & . & . & . & . & . & . & . & . \\
 . & . & 1 & 1 & . & . & 1 & 1 & . & 1 & . & 1 & . & . & . & . \\
 . & 1 & 2 & . & 2 & 1 & 1 & . & . & . & . & . & . & 1 & . & . \\
 . & . & 2 & . & . & . & 1 & . & . & 1 & . & 2 & . & 1 & 1 & . \\
 . & . & 2 & 1 & 1 & . & 2 & . & . & . & . & 1 & 1 & 1 & . & . \\
 . & . & . & . & 1 & . & . & . & . & . & . & . & . & . & . & . \\
 . & . & . & . & . & . & . & . & . & . & . & 1 & . & . & . & . \\
\end{array}
\right)
\\
\\
g_{12}=&
\left(
\begin{array}{cccccccccccccccc}
 . & . & . & . & . & . & . & . & . & . & . & 1 & . & . & . & . \\
 . & . & 2 & . & . & . & 1 & . & . & 1 & . & 2 & . & 1 & 1 & . \\
 . & 2 & 7 & 2 & 3 & . & 4 & . & . & 2 & 2 & 3 & 3 & 4 & . & . \\
 . & . & 2 & 1 & 1 & . & 2 & . & . & . & . & 1 & 1 & 1 & . & . \\
 . & . & 3 & 1 & . & . & 2 & 1 & . & 2 & . & 3 & . & 1 & 1 & . \\
 . & . & . & . & . & . & . & . & . & 1 & . & 1 & . & . & . & . \\
 . & 1 & 4 & 2 & 2 & . & 3 & . & . & 1 & 1 & 2 & 2 & 2 & . & . \\
 . & . & . & . & . & . & . & . & . & . & . & . & 1 & . & . & . \\
 . & . & . & . & 1 & . & . & . & . & . & . & . & . & . & . & . \\
 . & . & 2 & . & . & . & 1 & . & . & . & 1 & . & 2 & 1 & . & 1 \\
 . & 1 & 2 & . & 2 & 1 & 1 & . & . & . & . & . & . & 1 & . & . \\
 . & . & 3 & 1 & . & . & 2 & . & 1 & . & 2 & . & 3 & 1 & . & 1 \\
 1 & 2 & 3 & 1 & 3 & 1 & 2 & . & . & . & . & . & . & 1 & . & . \\
 . & 1 & 4 & 1 & 1 & . & 2 & . & . & 1 & 1 & 1 & 1 & 2 & . & . \\
 . & . & . & . & . & . & . & . & . & . & 1 & . & 1 & . & . & . \\
 . & 1 & . & . & 1 & . & . & . & . & . & . & . & . & . & . & . \\
\end{array}
\right)
\qquad
g_{13}=&
\left(
\begin{array}{cccccccccccccccc}
 . & . & . & . & . & . & . & . & . & . & . & . & 1 & . & . & . \\
 . & . & 2 & . & . & . & 1 & . & . & . & 1 & . & 2 & 1 & . & 1 \\
 . & 2 & 7 & 2 & 3 & . & 4 & . & . & 2 & 2 & 3 & 3 & 4 & . & . \\
 . & . & 2 & 1 & 1 & . & 2 & . & . & . & . & 1 & 1 & 1 & . & . \\
 . & . & 3 & 1 & . & . & 2 & . & 1 & . & 2 & . & 3 & 1 & . & 1 \\
 . & . & . & . & . & . & . & . & . & . & 1 & . & 1 & . & . & . \\
 . & 1 & 4 & 2 & 2 & . & 3 & . & . & 1 & 1 & 2 & 2 & 2 & . & . \\
 . & . & . & . & 1 & . & . & . & . & . & . & . & . & . & . & . \\
 . & . & . & . & . & . & . & . & . & . & . & 1 & . & . & . & . \\
 . & 1 & 2 & . & 2 & 1 & 1 & . & . & . & . & . & . & 1 & . & . \\
 . & . & 2 & . & . & . & 1 & . & . & 1 & . & 2 & . & 1 & 1 & . \\
 1 & 2 & 3 & 1 & 3 & 1 & 2 & . & . & . & . & . & . & 1 & . & . \\
 . & . & 3 & 1 & . & . & 2 & 1 & . & 2 & . & 3 & . & 1 & 1 & . \\
 . & 1 & 4 & 1 & 1 & . & 2 & . & . & 1 & 1 & 1 & 1 & 2 & . & . \\
 . & 1 & . & . & 1 & . & . & . & . & . & . & . & . & . & . & . \\
 . & . & . & . & . & . & . & . & . & 1 & . & 1 & . & . & . & . \\
\end{array}
\right)
\\
\\
g_{15}=&
\left(
\begin{array}{cccccccccccccccc}
 . & . & . & . & . & . & . & . & . & . & . & . & . & . & 1 & . \\
 . & . & . & . & . & . & . & . & . & . & . & 1 & . & . & . & . \\
 . & . & 1 & . & . & . & 1 & . & . & . & . & . & . & . & . & . \\
 . & . & . & . & . & . & . & . & . & . & . & . & . & 1 & . & . \\
 . & . & . & . & . & . & . & . & . & 1 & . & 1 & . & . & . & . \\
 . & . & . & . & . & . & . & 1 & . & . & . & . & . & . & 1 & . \\
 . & . & 1 & . & . & . & . & . & . & . & . & . & . & . & . & . \\
 . & . & . & . & . & . & . & . & . & . & . & . & . & . & . & 1 \\
 . & . & . & . & . & 1 & . & . & . & . & . & . & . & . & . & . \\
 . & . & . & . & . & . & . & . & . & . & . & . & 1 & . & . & . \\
 . & . & . & . & 1 & . & . & . & . & . & . & . & . & . & . & . \\
 . & . & . & . & . & . & . & . & . & . & 1 & . & 1 & . & . & . \\
 . & 1 & . & . & 1 & . & . & . & . & . & . & . & . & . & . & . \\
 . & . & . & 1 & . & . & . & . & . & . & . & . & . & 1 & . & . \\
 . & . & . & . & . & . & . & . & 1 & . & . & . & . & . & . & 1 \\
 1 & . & . & . & . & 1 & . & . & . & . & . & . & . & . & . & . \\
\end{array}
\right)
\qquad
g_{16}=&
\left(
\begin{array}{cccccccccccccccc}
 . & . & . & . & . & . & . & . & . & . & . & . & . & . & . & 1 \\
 . & . & . & . & . & . & . & . & . & . & . & . & 1 & . & . & . \\
 . & . & 1 & . & . & . & 1 & . & . & . & . & . & . & . & . & . \\
 . & . & . & . & . & . & . & . & . & . & . & . & . & 1 & . & . \\
 . & . & . & . & . & . & . & . & . & . & 1 & . & 1 & . & . & . \\
 . & . & . & . & . & . & . & . & 1 & . & . & . & . & . & . & 1 \\
 . & . & 1 & . & . & . & . & . & . & . & . & . & . & . & . & . \\
 . & . & . & . & . & 1 & . & . & . & . & . & . & . & . & . & . \\
 . & . & . & . & . & . & . & . & . & . & . & . & . & . & 1 & . \\
 . & . & . & . & 1 & . & . & . & . & . & . & . & . & . & . & . \\
 . & . & . & . & . & . & . & . & . & . & . & 1 & . & . & . & . \\
 . & 1 & . & . & 1 & . & . & . & . & . & . & . & . & . & . & . \\
 . & . & . & . & . & . & . & . & . & 1 & . & 1 & . & . & . & . \\
 . & . & . & 1 & . & . & . & . & . & . & . & . & . & 1 & . & . \\
 1 & . & . & . & . & 1 & . & . & . & . & . & . & . & . & . & . \\
 . & . & . & . & . & . & . & 1 & . & . & . & . & . & . & 1 & . \\
\end{array}
\right)
 \end{align*}
 }

 \subsubsection*{Acknowledgements}

The techniques used here were applied (and sometimes developed) by the author more than 20 years ago,
but this return to the past would not have been initiated without a specific question that J.-B. Zuber asked the author in November 2025: ``What are the nimreps for this specific example?''.
Part of the present article was presented in a talk at the conference ``Jean-Bernard is eighty'', held at the Institut Henri Poincaré in Paris on September 24--25, 2026.
It is dedicated to Jean-Bernard on the occasion of his 80th birthday.

\subsubsection*{AI disclosure}

Most calculations performed here were carried out using programs developed by the author at the beginning of the 2000s, long before the rise of artificial intelligence.
However, during the preparation of this work, the author used ChatGPT to improve pre-existing programs and to draft exploratory Mathematica code, in particular the code leading to the determination of rational polynomials.
At a later stage, it (it seems that technology has no gender) proofread the manuscript and corrected numerous typos.
We also used Gemini for short discussions and for suggesting an interesting piece of code that, unfortunately, could not be used here.
Most importantly, the author thanks both of these AI systems for their company.

 \end{document}